\documentclass{article}

\usepackage[english]{babel}

\usepackage{cite}
\usepackage{appendix}
\usepackage{enumitem}
\usepackage{microtype}

\usepackage[letterpaper,top=2cm,bottom=2cm,left=3cm,right=3cm,marginparwidth=1.75cm]{geometry}
\usepackage{amssymb}
\usepackage{amsthm}
\usepackage{amsmath}
\makeatletter

\newcommand{\Rmnum}[1]{\expandafter\@slowromancap\romannumeral #1@}
\makeatother
\usepackage{mathrsfs}
\usepackage{indentfirst}
\usepackage{amsmath}
\usepackage{graphicx}
\usepackage[colorlinks=true, allcolors=blue]{hyperref}
\usepackage{autobreak}

\newtheorem{theorem}{Theorem}[section] 

\newtheorem{definition}{Definition}[section] 
\newtheorem{lemma}{Lemma}[section] 
\newtheorem{corollary}[theorem]{Corollary}
\newtheorem{example}{Example}[section]
\newtheorem{proposition}[theorem]{Proposition}
\newtheorem{remark}{Remark}[section] 
\numberwithin{equation}{section}
\usepackage{amsmath}
\allowdisplaybreaks[4]

\usepackage{enumerate}

\usepackage{amssymb}
\usepackage{authblk}

\newcommand{\enabstractname}{Abstract}
\newcommand{\cnabstractname}{摘要}
\newenvironment{enabstract}{%
	\par
	\noindent\mbox{}\hfill{\bfseries \enabstractname}\hfill\mbox{}\par
	\vskip 2.5ex}{\par\vskip 2.5ex}

\usepackage{footmisc}

\title{
Averaging Principle and Uniform Large Deviations for Time-Dependent Slow–Fast Systems via Young–Rough Integration}
\author{Wenting Xu\thanks{%
		Email: \href{mailto:qfnuxwt@163.com}{qfnuxwt@163.com}.\\
		Graduate School of Mathematics, Kyushu University, 
		744 Motooka, Nishi-ku, Fukuoka 819-0395, Japan.
}}

\date{\small}
\begin{document}
	\maketitle
	\begin{enabstract}
	The main aim of this work is to investigate the limiting behavior and asymptotic properties of time-dependent slow–fast stochastic dynamical systems. The slow component is driven by a stochastic process with H\"older-continuous paths of regularity greater than 1/2, while the fast component is driven by a Brownian rough path. We introduce a novel Young–rough integration framework tailored to such systems and establish well-posedness results, providing the analytical foundation for further asymptotic analysis. Building on this framework, we employ Khasminskii's time-discretization method to establish the strong convergence of the slow component to the solution of the corresponding averaged equation in the \(L^p\)-sense. Furthermore, when the driving noise in the slow equation is specified as a fractional Brownian motion, we apply the weak convergence approach to establish a uniform large deviation principle for the slow component, uniformly over compact sets of initial conditions, in the H\"older topology.
		
		
	\end{enabstract}
	
\tableofcontents

\section{Introduction}

Slow--fast systems arise naturally in various fields, including materials science, chemistry, fluid dynamics, biology, ecology, and climate dynamics \cite{bertram2017multi, harvey2011multiple, mastny2007two, kuehn2015multiple}. These systems consist of coupled subsystems evolving on different time scales and are often subject to stochastic perturbations. In particular, the driving noises associated with different components may exhibit distinct pathwise regularities and statistical properties, which pose additional challenges for the mathematical analysis of the resulting dynamics. The study of such systems therefore calls for appropriate analytical frameworks and integration techniques capable of accommodating driving signals of different regularities.

In this work, we mainly consider the time--dependent slow--fast system
\begin{align}\label{eq-1'}
	\begin{cases}
		dX^{\varepsilon,\delta}_t=b_1(t,X^{\varepsilon,\delta}_t,Y^{\varepsilon,\delta}_t)dt+\sqrt\varepsilon\sigma_1(t,X^{\varepsilon,\delta}_t)dB_t,\\
		dY^{\varepsilon,\delta}_t=\frac1\delta b_2(X^{\varepsilon,\delta}_t,Y^{\varepsilon,\delta}_t)dt+
		\frac{1}{\sqrt{\delta}}\sigma_2(X^{\varepsilon,\delta}_t,Y^{\varepsilon,\delta}_t)d^Iw_t,\\
		(X^{\varepsilon,\delta}_0,Y^{\varepsilon,\delta}_0)=(x_0,y_0)\in\mathbb{R}^{d_1+d_2},~ t\in[0,T],~ 0<\delta<\varepsilon\leq 1.
	\end{cases}
\end{align}
Here, $d^I$ denotes integration in the standard It\^o sense, \(B=(B_t)_{t\in[0,T]}\) is an \(\mathbb{R}^m\)-valued stochastic process with \(\alpha\)-H\"older continuous sample paths for \(\alpha\in(1/2,1)\), and \(w=(w_t)_{t\in[0,T]}\) is an \(r\)-dimensional Brownian motion (BM). 
For sufficiently small $\delta$,  $X^{\varepsilon,\delta}$ and $Y^{\varepsilon,\delta}$ are called the slow and fast component, respectively.

BM is a continuous semimartingale and thus admits stochastic integration in the It\^o and Stratonovich frameworks. For more general processes, such as fractional Brownian motion (FBM), which is not a semimartingale unless the Hurst parameter \(H=1/2\), several alternative integration theories have been developed. These include fractional calculus, the stochastic sewing lemma, and Young and rough path integration. Z\"ahle \cite{zahle1998integration} introduced a pathwise integration approach based on fractional calculus, while L\^e \cite{le2020stochastic} developed the stochastic sewing lemma, which constructs integrals from local increments by exploiting the regularity requirements on conditional expectations and random fluctuations. In contrast to the deterministic sewing lemma \cite{gubinelli2004controlling}, this approach accommodates weaker local error estimates and is therefore particularly useful in stochastic settings. More naturally, Young \cite{young1936inequality} established the existence of Riemann--Stieltjes integrals when the integrand and integrator have finite \(p\)- and \(q\)-variation, respectively, with \(\frac1p+\frac1q>1\). For rougher signals, Lyons \cite{lyons1994interpretation,lyons1994differential} developed rough path theory, which incorporates iterated-integral information of the driving signal and provides a framework for establishing well-posedness of differential equations driven by irregular paths. In particular, BM can be enhanced to a Brownian rough path by including its iterated integrals. Young integration can be viewed as the  simplest case of rough path integration, which extends pathwise integration to rougher signals through higher-order iterated integrals.

For slow--fast systems driven by standard BMs only, stochastic integration provides the standard framework for formulating and analyzing the system. Within the It\^o framework, Khasminskii \cite{khas1968averaging} developed a time--discretization approach to the averaging principle (AP), while Freidlin and Wentzell \cite{Wentzell1984} studied averaging with respect to the invariant measure of the fast dynamics. Such stochastic integration also provides a foundation for investigating the asymptotic behavior of slow--fast systems \cite{sun2021large,gasteratos2024importance,freidlin2012random,kifer2009large,liu2010strong,de2021averaging,liu2020averaging}.

Beyond the classical Brownian setting, there has been growing interest in understanding the behavior of stochastic systems driven by more general stochastic processes. Among these, FBM has attracted considerable attention as a prominent example of a non-Markovian Gaussian process. Since the analysis of long--time behavior typically relies on the existence of an invariant measure and the ergodic properties of the fast subsystem, existing studies have primarily focused on systems in which the slow equation is driven by FBM, while the fast equation is driven by standard BM. When the Hurst parameter satisfies $H>1/2$, Young integration and fractional calculus provide suitable tools for handling the corresponding integral terms. Guerra and Nualart \cite{guerra2008stochastic} combined fractional integration with classical It\^o stochastic calculus to establish well-posedness for equations driven by mixed FBM, paving the way for subsequent studies of AP and large deviation principle (LDP) \cite{pei2020stochastic, shen2022averaging,shen2024large, pei2023averaging, gailus2025large, inahama2023large}.
More recently, a few studies have employed the stochastic sewing lemma to investigate fully coupled slow–fast systems \cite{hairer2020averaging, li2025averaging, li2026large}.
For FBM with Hurst parameter $H\leq 1/2$, rough path theory \cite{friz2014course} provides a powerful framework for defining and analyzing the corresponding stochastic differential equations (SDEs). 
By lifting both the slow and fast equations to the rough path framework, Pei, Inahama and Xu \cite{pei2021averaging} first established well-posedness and AP for slow--fast systems using rough path theory via fractional calculus \cite{hu2009rough}. 
Subsequently, Inahama \cite{inahama2025averaging,inahama2025averaging2} established corresponding well-posedness results and AP within the controlled rough path framework under weaker assumptions, with \(L^p\)-convergence in H\"older topology, improving both the mode of convergence and the topology in which convergence is established.
Building on these results, Yang and Xu \cite{yang2025large} established an LDP for such systems.

A natural approach to dealing with a system driven by both a stochastic process with sufficiently regular H\"older continuous sample paths and BM is to interpret the corresponding integrals in the Young and It\^o senses, respectively. 
However, Young and It\^o integration are based on different mathematical constructions, which presents difficulties in analyzing the slow and fast components within a common analytical framework. 
Although the system can alternatively be formulated in terms of rough differential equations (RDEs), treating the more regular driving signal solely within rough path theory may fail to fully exploit its higher regularity and may require unnecessarily strong regularity assumptions on the coefficients.
To overcome these limitations, we extend the controlled rough path framework to incorporate Young integration while preserving the distinct regularity properties of the driving signals. This leads to a unified Young--rough integration framework in which the two types of integrals can be treated consistently and their interactions can be analyzed through refined estimates. 
The construction and fundamental properties of the generalized rough integral are presented in Section \ref{sec-2}.

We adopt a pathwise approach to the analysis of the slow--fast system \eqref{eq-1'}, exploiting the deterministic nature of Young and rough integration. As a first step, we consider the corresponding system without scale separation, obtained by setting both the noise intensity parameter $\varepsilon$ and the time-scale parameter $\delta$ equal to one.
The resulting system is given by \eqref{eq-2} in Section \ref{sec-3}.
In this setting, we interpret the first equation as a Young differential equation (YDE) and the second equation as a RDE within the generalized rough integration framework. By applying a fixed-point argument, we establish the global well-posedness of this deterministic Young--rough system and derive stability estimates with respect to the initial conditions and driving signals. The main well-posedness result is stated in Theorem \ref{th-1}.

The well-posedness and continuity of the resulting solution map provide a deterministic foundation for the pathwise analysis of stochastic systems driven by signals with the required regularity. In particular, by considering suitable sample paths of stochastic processes and their rough path lifts, the deterministic theory can be applied to stochastic slow--fast systems, subject to appropriate integrability and uniform estimates. A central objective in the study of slow--fast systems is to characterize the limiting and asymptotic behavior of the slow component as the time-scale separation becomes increasingly pronounced. Building on the deterministic well-posedness theory developed above, we investigate the AP and the uniform LDP for system \eqref{eq-1'} within the generalized Young--rough integration framework.
The Young--rough framework enables us to establish stronger results under less restrictive assumptions than those imposed in certain existing integration frameworks.

The AP provides a fundamental tool for analyzing slow--fast systems by approximating the slow component over finite time intervals through an effective averaged equation as the time-scale separation parameter tends to zero. In this work, we impose suitable regularity and dissipativity conditions on $b_2$ and $\sigma_2$ to ensure that, for every $x\in\mathbb{R}^{d_1}$, the following frozen SDE is exponentially ergodic:
\begin{align*}
	dY^{x,y}_t
	=b_2(x,Y^{x,y}_t)dt
	+\sigma_2(x,Y^{x,y}_t)d^Iw_t,
	\qquad
	Y^{x,y}_0=y\in\mathbb{R}^{d_2}.
\end{align*}
More precisely, the frozen equation admits a unique invariant probability measure $\mu^x$, and the distribution of its solution converges to $\mu^x$ at an exponential rate. This property implies that the frozen fast dynamics rapidly approach their statistical equilibrium, allowing their long-time behavior to be characterized by the invariant measure $\mu^x$. Consequently, the influence of the fast component on the slow dynamics can be captured through the averaged drift
$	\bar b_1(t,x) =\int_{\mathbb{R}^{d_2}} b_1(t,x,z)\mu^x(dz)$.

To investigate the averaging effect independently of the small-noise limit, we study the AP under fixed noise intensity. In this setting, the averaged equation retains the stochastic forcing and thus captures the effective stochastic dynamics of the slow component. Specifically, setting $\varepsilon=1$, we consider the following averaged equation:
\begin{align*}
	d\bar X_t
	=\bar b_1(t,\bar X_t)dt
	+\sigma_1(t,\bar X_t)dB_t,
	\qquad
	\bar X_0=x_0\in\mathbb{R}^{d_1}.
\end{align*}
We will show that for every $1\leq p<\infty$ and $\alpha^\prime\in(1/2, \alpha)$, we have 
$$
\lim_{\delta\to 0}\mathbb{E}[\|X^{1,\delta}-\bar{X}\|^p_{\alpha^\prime}]=0.
$$
Our main result of AP is given by Theorem \ref{th-2}.

Various integration frameworks, including It\^o integration, Young integration, fractional calculus, and rough path theory, have been employed to establish APs under different regularity assumptions \cite{liu2020averaging,pei2023averaging,pei2020stochastic, shen2022averaging, pei2021averaging,inahama2025averaging,inahama2025averaging2}.
Liu et al. \cite{liu2020averaging} established an AP in the uniform topology for time--dependent slow--fast systems driven by two BMs. Their results allow locally Lipschitz drift coefficients in the spatial variables but require additional time regularity, with stopping-time techniques used to overcome the difficulties arising from local Lipschitz continuity.
For time-dependent slow--fast systems driven by mixed FBM with Hurst parameter \(H>1/2\), Pei, Inahama and Xu \cite{pei2023averaging} established an AP with \(L^2\)-convergence using It\^o integration for BM and fractional calculus for FBM, under a global Lipschitz condition on the fast drift. 
In contrast, our controlled rough path framework, through its favorable topological properties and a priori estimates, allows us to obtain \(L^p\)-convergence in a H\"older topology. Moreover, we require only local Lipschitz continuity of the fast drift, with additional time regularity imposed solely on the diffusion coefficient of the slow equation.
Under a rather different set of assumptions, for fully coupled slow--fast systems in which the slow equation is driven by FBM, the stochastic sewing lemma provides a route to obtaining uniform estimates of the slow component. Hairer and Li \cite{hairer2020averaging} established an AP for fully coupled slow--fast systems driven by mixed FBM, proving convergence in probability via the stochastic sewing lemma. Next, Li, Gao and Qu \cite{li2025averaging} studied the AP in the infinite-dimensional setting, but with additive noise in the fast equation. 
Compared with the stochastic sewing approach, our Young--rough integration framework requires less restrictive regularity and boundedness assumptions on the fast-equation coefficients. It also yields stronger \(L^p\)-convergence in a H\"older topology and quantitative decay estimates for the error moments, while providing a more direct pathwise analysis that avoids elaborate probabilistic estimates. Moreover, our framework extends beyond FBM to more general driving signals with H\"older regularity greater than \(1/2\).
A detailed comparison of our results with those of preceding works is provided in Remark \ref{rem-3} below.

The AP and LDP address two distinct asymptotic mechanisms in stochastic systems: the former characterizes the effective dynamics arising from time-scale separation, whereas the latter quantifies the exponential decay of probabilities of rare deviations in the small-noise regime.
The theory of small-noise large deviations was developed in the 1960s by Freidlin, Wentzell, Varadhan, and others \cite{Wentzell1984,varadhan1984large}. Subsequent foundational works established a systematic framework for studying LDPs for stochastic differential equations driven by BM \cite{freidlin2012random,boue1998variational,budhiraja2000variational,budhiraja2008large}.

Classical LDPs typically describe the asymptotic behavior of stochastic systems with fixed initial conditions, with the corresponding rate functions depending on the initial states. However, in many applications, the initial conditions may vary within a prescribed subset of the state space, making pointwise large deviation estimates insufficient. This motivates the study of uniform LDPs, which provide a family of rate functions indexed by the initial conditions and ensure that the corresponding large deviation estimates hold uniformly over the prescribed set. Such uniformity is particularly important in problems involving varying initial states, such as the analysis of exit times of stochastic processes from a domain \cite{dembo2009large,freidlin2012random}.
Several notions of uniform LDP have been introduced in the literature. The first, due to Freidlin and Wentzell \cite{freidlin2012random}, will be referred to as the Freidlin--Wentzell uniform LDP. Another formulation, presented in \cite{dembo2009large}, will be called the Dembo--Zeitouni uniform LDP. A third formulation, known as the uniform Laplace principle, can be found in \cite{feng2006large}. These formulations are not equivalent in general.
Moreover, Salins \cite{salins2019equivalences} introduced the equicontinuous uniform Laplace principle, which strengthens the uniform Laplace principle by requiring uniform convergence over equibounded, equicontinuous families of test functions. Salins further established the equivalence between the equicontinuous uniform Laplace principle and the Freidlin--Wentzell uniform LDP without additional assumptions. An important advantage of the equicontinuous uniform Laplace principle is that it can be established using the variational methods developed by Budhiraja, Dupuis, and their collaborators \cite{boue1998variational,budhiraja2000variational,budhiraja2008large}.
In subsequent work, Salins \cite{salins2019uniform} developed general methods for establishing LDPs for small-noise SDEs taking values in separable Banach spaces, uniformly over bounded sets of initial conditions.

When studying small-noise large deviations for slow--fast systems with the limit of time-scale separation, one needs to consider the joint limit in which both the noise intensity $\varepsilon$ and the time-scale parameter $\delta$ tend to zero. The interplay between these two asymptotic effects introduces additional challenges to the large deviation analysis. In this setting, averaging techniques play a crucial role in identifying the effective dynamics and characterizing the associated rate functions.
For slow--fast systems driven by BM, Sun et al. \cite{sun2021large} investigated the LDP for stochastic Burgers equations with two time scales. For systems driven by FBM with Hurst parameter $H>1/2$, Gailus and Gasteratos \cite{gailus2025large} and Inahama et al. \cite{inahama2023large} established LDPs for slow--fast systems driven by mixed FBM. Shen et al. \cite{shen2024large} investigated distribution-dependent slow--fast systems, while Xu et al. \cite{xu2025large} extended the analysis to infinite-dimensional slow--fast systems. For FBM with $H<1/2$, Yang and Xu \cite{yang2025large,yang2025large2} recently established LDPs for slow--fast systems with $H\in(1/4,1/2)$ using rough path theory.
Concerning uniform LDPs, Li and Liu \cite{li2026uniform} recently established a uniform LDP for slow--fast McKean--Vlasov systems driven by standard BM. Nevertheless, uniform LDPs for slow--fast systems, particularly those driven by FBM, remain relatively unexplored. A detailed comparison of our results with previous works is provided in Remark \ref{rem-5} below.

In this work, we establish a uniform LDP for the slow component under the asymptotic regime \(0<\delta<\varepsilon<1\) and \(\delta/\varepsilon\to0\) as \(\varepsilon\to0\), restricting the driving noise \(B\) in the slow equation to be an FBM \(B^H\) with Hurst parameter \(H\in(1/2,1)\). Under suitable assumptions, we employ the weak convergence approach and the equivalence between the equicontinuous uniform Laplace principle and the Freidlin--Wentzell uniform LDP to establish the Freidlin--Wentzell uniform LDP for the slow component in the $(H-\kappa)$-H\"older topology  for any sufficiently small $\kappa >0$, uniformly with respect to its initial conditions. The result is stated precisely in Theorem \ref{th-3}.

The rest of the paper is organized as follows.
Section \ref{sec-2} extends certain classical results from rough path theory to accommodate diffusion coefficients with a mixed Young--rough structure. It develops a generalized controlled rough path framework and defines the associated rough integral. In particular, we show that the rough integral of a generalized controlled rough path is a controlled rough path in the classical sense.
Section \ref{sec-3} establishes the well-posedness of the deterministic system without time-scale separation and derives the corresponding stability estimates within the Young--rough framework.
Section \ref{sec-4} specifies the assumptions under which an AP is established, with convergence in \(L^p\) in the H\"older topology.
Section \ref{sec-5} imposes additional assumptions and establishes the Freidlin--Wentzell uniform LDP for the slow component in the H\"older topology, uniformly with respect to its initial conditions.
Section \ref{sec-A} provides the precise definitions of the equicontinuous uniform Laplace principle and the Freidlin--Wentzell uniform LDP, and presents the equivalence between these two principles.

Before concluding the introduction, we introduce the notation which will be used throughout the paper.
Let $T>0$ be arbitrary and we work on the time interval $[0,T]$ unless otherwise specified. Whenever there is no ambiguity concerning the time interval, we omit $[0,T]$ from the notation. 
Below, $\mathcal{V}$, $\mathcal{W}$, $\bar{\mathcal{V}}$, $\bar{\mathcal{W}}$ are Euclidean spaces. The set of bounded linear maps from $\mathcal{V}$ to $\mathcal{W}$ is denoted by $\mathcal{L}(\mathcal{V},\mathcal{W})$. When $\mathcal{V}=\mathbb{R}^n$ and $\mathcal{W}=\mathbb{R}^d$, $\mathcal{L}(\mathcal{V},\mathcal{W})$ coincides with the set of all real $d\times n$ matrices and is equipped with Hilbert-Schmidt norm instead of the operator norm. Moreover, for a function $\varphi:[0,T]\to \mathcal{V}$, we write $\varphi^1_{s,t}:=\varphi_t-\varphi_s$ for $0\leq s\leq t\leq T$.

For $\alpha\in(0,1]$, let $\mathcal{C}^{\alpha}=\mathcal{C}^{\alpha}(\mathcal{V})$ denote the space of $\alpha$-H\"older continuous paths taking values in $\mathcal{V}$. More precisely, we define
$$
\mathcal{C}^{\alpha}
:=\mathcal{C}^{\alpha}\big( \mathcal{V}\big)
:=\big\{
	\varphi\in \mathcal{C}( \mathcal{V}):
	\|\varphi\|_{\alpha}<\infty
	\big\},
$$
where
$$
\|\varphi\|_{\alpha}
:=\sup_{0\leq s<t\leq T}
\frac{|\varphi_t-\varphi_s|}{|t-s|^{\alpha}}
$$
is the $\alpha$-H\"older seminorm, and $\mathcal{C}( \mathcal{V})$ denotes the space of all continuous functions from $[0,T]$ to $\mathcal{V}$. 
For $\tilde{\alpha}\in(0,1)$, let $\tilde{\mathcal{C}}^{\tilde{\alpha}}=\tilde{\mathcal{C}}^{\tilde{\alpha}}(\mathcal{V})$ denote the little H\"older space, that is, the space consisting of all $\varphi\in\mathcal{C}^{\tilde{\alpha}}$ satisfying
$$\lim_{\delta\to0^+}
\sup_{\substack{0\leq s<t\leq T\\t-s\leq\delta}} \frac{|\varphi_t-\varphi_s|}{(t-s)^{\tilde{\alpha}}}=0.$$
Equipped with the norm $|\varphi_0|+\|\varphi\|_{\tilde{\alpha}}$, the space $\tilde{\mathcal{C}}^{\tilde{\alpha}}$ is a separable Banach space.
Moreover, $\tilde{\mathcal{C}}^{\tilde{\alpha}}$
coincides with the closure of 
$\bigcup_{0<\iota<1-{\tilde{\alpha}}}
	\mathcal{C}^{{\tilde{\alpha}}+\iota}(\mathcal{V})$ with respect to the above norm.
Furthermore, $\tilde{\mathcal{C}}^{\tilde{\alpha}}$ is continuously embedded in $\mathcal{C}^{\tilde{\alpha}}(\mathcal{V})$ (see \cite{ciesielski1960isomorphisms}).

Let $\Delta_{[0,T]}:=\big\{(s,t):0\leq s\leq t\leq T\big\}$
denote the $2$-simplex. For $\kappa\in(0,1)$, we define $\mathcal{C}_2^\kappa(\mathcal{V})$ as the space of continuous maps $G:\Delta_{[0,T]}\to\mathcal{V}$ satisfying $G_{t,t}=0$ for every $t\in[0,T]$ and
$$
\|G\|_{\kappa}
:=\sup_{0\leq s<t\leq T}
\frac{|G_{s,t}|}{|t-s|^\kappa}<\infty.
$$
Thus, $\mathcal{C}_2^{\kappa}(\mathcal{V})$ consists of two-parameter increments that are $\kappa$-H\"older continuous in the sense of this seminorm.

For $\kappa_1\in(0,1]$ and $\kappa_2>0$, we further define $\mathcal{C}_2^{\kappa_1,\kappa_2}(\mathcal{V})$ as the space of continuous maps $J:\Delta_{[0,T]}\to\mathcal{V}$ satisfying $J_{t,t}=0$ and
$$
\|J\|_{\kappa_1,\kappa_2}
:=\|J\|_{\kappa_1}
+\|\delta J\|_{\kappa_2}<\infty,
$$
where
$$
\|J\|_{\kappa_1}
:=\sup_{0\leq s<t\leq T}
\frac{|J_{s,t}|}{|t-s|^{\kappa_1}},
$$
and
$$
\delta J_{s,u,t}
:=J_{s,t}-J_{s,u}-J_{u,t},
\qquad 0\leq s\leq u\leq t\leq T.
$$
The seminorm of the three-parameter increment $\delta J$ is given by
$$
\|\delta J\|_{\kappa_2}
:=\sup_{0\leq s<u<t\leq T}
\frac{|\delta J_{s,u,t}|}{|t-s|^{\kappa_2}}.
$$
The quantity $\|J\|_{\kappa_1,\kappa_2}$ measures both the H\"older regularity of $J$ and the regularity of its three-point increment $\delta J$.

For the rest of this section, let $\beta\in(1/3,1/2]$. 
A continuous map
$$
\mathbb{X}=(X^1,X^2):
\Delta_{[0,T]}
\longrightarrow
\mathcal{V}\oplus(\mathcal{V}\otimes\mathcal{V})
$$
is called a $\mathcal{V}$-valued $\beta$-H\"older rough path if
$$
\|X^1\|_{\beta}+\|X^2\|_{2\beta}<\infty
$$
and Chen's relations hold:
$$
X^1_{s,t}=X^1_{s,u}+X^1_{u,t},
$$
and
$$
X^2_{s,t}
=X^2_{s,u}+X^2_{u,t}
+X^1_{s,u}\otimes X^1_{u,t},
\qquad 0\leq s\leq u\leq t\leq T.
$$
Here, $\|X^1\|_{\beta}$ and $\|X^2\|_{2\beta}$ denote the corresponding H\"older seminorms on $\Delta_{[0,T]}$. The set of all $\mathcal{V}$-valued $\beta$-H\"older rough paths is denoted by $\Omega^{\beta}(\mathcal{V})$. 
The $\beta$-H\"{o}lder homogeneous rough path norm is denoted by
$$\vert\!\vert\!\vert \mathbb{X} \vert\!\vert\!\vert_{\beta}:=\|X^1\|_\beta+\sqrt{\|X^2\|_{2\beta}}.$$

Finally, for a fixed path $X^1\in\mathcal{C}^{\beta}(\mathcal{V})$, we denote by $\hat{\mathscr{D}}^\beta_{X^1}(\mathcal{W})$ the space of $\mathcal{W}$-valued paths controlled by $X^1$. More precisely, its elements are pairs $(Y,Y^\dagger)$, where $Y\in\mathcal{C}^{\beta}(\mathcal{W})$ and $Y^\dagger\in\mathcal{C}^{\beta}(\mathcal{L}(\mathcal{V},\mathcal{W}))$ such that the remainder
$$
Y^\sharp_{s,t}:=Y_t-Y_s-Y^\dagger_sX^1_{s,t},
\qquad 0\leq s \leq t\leq T
$$
has finite $2\beta$-H\"older seminorm.

For a function $f:\mathcal{V}\to\mathcal{W}$ of class $C^k$, $k\geq 2$, we denote its first- and \(i\)-th-order gradients by \(\nabla f\) and \(\nabla^i f\), respectively, for \(2\leq i\leq k\).
Let $F:\mathcal{V}\times\mathcal{\bar{V}}\to\mathcal{W}$ be a $C^k$ function, where $k\geq 2$.
Then, $\nabla_x F$ and $\nabla_y F$ denote the partial gradients of $F$ with respect to the first and second variables, respectively, while $\nabla_{(x,y)} F$ denotes the gradient with respect to the pair $(x,y)$. More precisely, for any $(x,y),(h,l)\in\mathcal{V}\times\mathcal{\bar{V}}$, we denote
\begin{align*}
\nabla_{(x,y)}F(x,y)\langle h,l\rangle
:=&\nabla_xF(x,y)\langle h\rangle+
\nabla_yF(x,y)\langle l\rangle\\
:=&
\lim_{s\to 0}\frac{F(x+sh,y)-F(x,y)}{s}+\lim_{u\to 0}\frac{F(x,y+ul)-F(x,y)}{u}
.
\end{align*}
For simplicity, we write $\nabla$ in place of $\nabla_{(x,y)}$.
Higher-order gradients are denoted analogously. For instance, $\nabla_{xx}^2 F$, $\nabla_{yy}^2 F$, and $\nabla_{xy}^2 F$ represent the corresponding second-order gradient, while $\nabla^2F$ denotes the second-order gradient with respect to the pair $(x,y)$.

\section{Generalized Rough Integration}\label{sec-2}

Let $\alpha\in(1/2,1]$ and $\beta\in(1/3,1/2]$ satisfy $\alpha+\beta>1$. Throughout this section, $\mathcal{V}$, $\mathcal{W}$, etc. denote Euclidean spaces.
The norms of finite dimensional vector spaces are simply denoted by $|\cdot|$.

We now introduce generalized controlled paths adapted to the mixed Young--rough setting.
\begin{definition}[Generalized controlled path]
Given a path $X^1\in \mathcal{C}^\beta(\mathcal{V})$, we say that $Y\in \mathcal{C}^\beta(\mathcal{W})$ is controlled by $X^1$ if there exists $Y^\dagger\in \mathcal{C}^\beta \big(\mathcal{L}(\mathcal{V},\mathcal{W})\big)$ so that the remainder term $Y^\sharp$ given implicitly through the relation 
\begin{align}\label{eq-22}
	Y^1_{s,t}=Y^\dagger_sX^1_{s,t}+Y^\sharp_{s,t},\quad0\leq s\leq t\leq T,
\end{align}
satisfies $\|Y^\sharp\|_{\alpha\wedge 2\beta}<\infty$. This defines the space of generalized controlled paths,
$$(Y,Y^\dagger,Y^\sharp)\in\mathscr{D}^{\alpha,\beta}_{X^1}(\mathcal{W}).$$
For simplicity, we shall often write $(Y,Y^\dagger)$ instead of $(Y,Y^\dagger,Y^\sharp)$.
The space $\mathscr{D}^{\alpha,\beta}_{X^1}(\mathcal{W})$ is equipped with the seminorm 
$$
\|(Y,Y^\dagger)\|_{\mathscr{D}^{\alpha,\beta}_{X^1}}
:=
\|Y^\dagger\|_\beta+\|Y^\sharp\|_{\alpha\wedge 2\beta}
.
$$
Moreover, $\mathscr{D}^{\alpha,\beta}_{X^1}(\mathcal{W})$ becomes a Banach space when endowed with the norm
$$t\mapsto|Y_0|+|Y^\dagger_0|+\|(Y,Y^\dagger)\|_{\mathscr{D}^{\alpha,\beta}_{X^1}}.$$
In particular, when $\alpha=2\beta$, the space $\mathscr{D}^{\alpha,\beta}_{X^1}(\mathcal{W})$ coincides with the previously defined space $\hat{\mathscr{D}}^\beta_{X^1}(\mathcal{W})$, which is referred to as the space of normal controlled paths.
Note that $\hat{\mathscr{D}}^\beta_{X^1}(\mathcal{W})\subset\mathscr{D}^{\alpha,\beta}_{X^1}(\mathcal{W})$.

\end{definition}

\begin{example}
	Here are a few typical examples of generalized controlled paths for a given rough path  $\mathbb{Z}=(Z^1,Z^2)\in\Omega^\beta(\mathcal{V})$.
	
	\begin{itemize}
		\item[1.] For $\xi\in\mathcal{W}$ and $\varphi\in\mathcal{L}(\mathcal{V},\mathcal{W})$, $t\to (\xi+\varphi Z^1_{0,t},\varphi)$ belongs to $\hat{\mathscr{D}}^\beta_{Z^1}(\mathcal{W})\subset\mathscr{D}^{\alpha,\beta}_{Z^1}(\mathcal{W})$. Note that $\mathscr{D}^{\alpha,\beta}_{Z^1}$-seminorm of this element is zero.
		\item[2.] If $\psi\in\mathcal{C}^{\alpha\wedge2\beta}(\mathcal{W})$, then obviously $(\psi,0)\in\mathscr{D}^{\alpha,\beta}_{Z^1}(\mathcal{W})$ with $\|\psi\|_{\alpha\wedge2\beta}=\|\psi\|_{\mathscr{D}^{\alpha,\beta}_{Z^1}}$. In this way, we have a natural continuous embedding $\mathcal{C}^{\alpha\wedge2\beta}(\mathcal{W})\hookrightarrow \mathscr{D}^{\alpha,\beta}_{Z^1}(\mathcal{W})$.
		\item[3.] Suppose that $(Q,Q^\dagger)\in\mathscr{D}^{\alpha,\beta}_{Z^1}(\mathcal{W})$, and $g:~\mathcal{W}\to\mathcal{\bar{W}}$ is a $C^2$-function. Then, $\big(g(Q),g(Q)^\dagger\big)\in\mathscr{D}^{\alpha,\beta}_{Z^1}(\bar{\mathcal{W}})$ if we set $g(Q)_t:=g(Q_t)$ and $g(Q)^\dagger_t:=\nabla g(Q_t)Q^\dagger_t$, where the right hand side is the composition of the two linear maps $\nabla g(Q_t)\in\mathcal{L}(\mathcal{W},\bar{\mathcal{W}})$ and $Q^\dagger_t\in\mathcal{L}(\mathcal{V},\mathcal{W})$. Using Taylor expansion, we can verify this fact as follows:
		\begin{align*}
			g(Q)^1_{s,t}
			=&\nabla g(Q_s) Q^1_{s,t}+\int_0^1(1-\theta)\nabla^2g(Q_s+\theta Q^1_{s,t})d\theta\langle Q^1_{s,t},Q^1_{s,t}\rangle\\
			=& \nabla g(Q_s)Q^\dagger_{s}Z^1_{s,t}+g(Q)^\sharp_{s,t},
		\end{align*}
		with
		\begin{align*}
			g(Q)^\sharp_{s,t}:=\nabla g(Q_s) Q^\sharp_{s,t}+\int_0^1(1-\theta)\nabla^2g(Q_s+\theta Q^1_{s,t})d\theta\langle Q^1_{s,t},Q^1_{s,t}\rangle.
		\end{align*}
		From the assumptions that $Q\in\mathcal{C}^\beta(\mathcal{W})$ and $Q^\sharp\in\mathcal{C}^{\alpha\wedge2\beta}_2(\mathcal{W})$, it is easy to see that $g(Q)^\sharp\in\mathcal{C}^{\alpha\wedge2\beta}_2(\bar{\mathcal{W}})$.
		
		\item[4.]\label{lem-1}
		As a corollary of the first and third examples, we introduce a special class of generalized controlled paths.
		For later purposes, we replace the state space by $\mathcal{L}(\mathcal{V},\mathcal{W})$.

		Let $f\in C^2$ with $f: \bar{\mathcal{V}}\times \mathcal{W}\to \mathcal{L}(\mathcal{V},\mathcal{W})$, where $\mathcal{V}$, $\bar{\mathcal{V}}$ and $\mathcal{W}$ are Euclidean spaces. 
		Let $X\in\mathcal{C}^\alpha(\bar{\mathcal{V}})$, and let $(Y,Y^\dagger)\in \hat{\mathscr{D}}^\beta_{Z^1}(\mathcal{W})$, where $Z^1\in\mathcal{C}^\beta(\mathcal{V})$.   
		Then, the map $f(X,Y)\in\mathcal{C}^\beta\big(\mathcal{L}(\mathcal{V},\mathcal{W})\big)$ is controlled by $Z^1$. More precisely, there exists $f(X,Y)^\dagger=\nabla_yf(X_t,Y_t)Y^\dagger_t\in \mathcal{C}^\beta\big(\mathcal{L}\big(\mathcal{V},\mathcal{L}(\mathcal{V},\mathcal{W})\big)\big)$ such that $$f(X,Y)^1_{s,t}=f(X,Y)^\dagger_sZ^1_{s,t}+f(X,Y)^\sharp_{s,t},$$
		where the remainder satisfies $\|f(X,Y)^\sharp\|_{\alpha\wedge2\beta}<\infty$. 
		In particular,
		$$\big(f(X,Y),f(X,Y)^\dagger\big)\in\mathscr{D}^{\alpha,\beta}_{Z^1}\big(\mathcal{L}(\mathcal{V},\mathcal{W})\big).$$
		Moreover, by a Taylor expansion, we obtain
		\begin{align*}
			f(X_t,Y_t)-f(X_s,Y_s)=&\nabla_x f(X_s,Y_s)X^1_{s,t}+\nabla_y f(X_s,Y_s)Y^1_{s,t}\\
			&+\int_0^1 (1-\theta)\nabla^2_{xx}f(X_s+\theta X^1_{s,t},Y_s+\theta Y^1_{s,t}) \langle X^1_{s,t},X^1_{s,t},\cdot\rangle d\theta\\
			&+2\int_0^1 (1-\theta)\nabla^2_{xy}f(X_s+\theta X^1_{s,t},Y_s+\theta Y^1_{s,t}) \langle X^1_{s,t},Y^1_{s,t},\cdot\rangle d\theta\\
			&+\int_0^1 (1-\theta)\nabla^2_{yy}f(X_s+\theta X^1_{s,t},Y_s+\theta Y^1_{s,t}) \langle  Y^1_{s,t},Y^1_{s,t},\cdot\rangle d\theta\\
			=&\nabla_y f(X_s,Y_s)Y^\dagger_s Z^1_{s,t}+f(X,Y)^\sharp_{s,t},
		\end{align*}
		where $\cdot\in\mathcal{V}$, with
		\begin{align*}
			f(X,Y)^\sharp_{s,t}=&\nabla_x f(X_s,Y_s)X^1_{s,t}+\nabla_y f(X_s,Y_s)Y^\sharp_{s,t}\\
			&+\int_0^1 (1-\theta)\nabla^2_{xx}f(X_s+\theta X^1_{s,t},Y_s+\theta Y^1_{s,t}) \langle X^1_{s,t},X^1_{s,t},\cdot\rangle d\theta\\
			&+2\int_0^1 (1-\theta)\nabla^2_{xy}f(X_s+\theta X^1_{s,t},Y_s+\theta Y^1_{s,t}) \langle X^1_{s,t},Y^1_{s,t},\cdot\rangle d\theta\\
			&+\int_0^1 (1-\theta)\nabla^2_{yy}f(X_s+\theta X^1_{s,t},Y_s+\theta Y^1_{s,t}) \langle  Y^1_{s,t},Y^1_{s,t},\cdot\rangle d\theta
			.
		\end{align*}
		Since $(Y,Y^\dagger)\in \hat{\mathscr{D}}^\beta_{Z^1}(\mathcal{W})$, it is easy  to see that 
		$\|f(X,Y)^\sharp\|_{\alpha\wedge2\beta}<\infty$.

		\item[5.] Concatenation of two generalized controlled paths is also a generalized controlled path. Let $0\leq a<b<c\leq T$. For $(\tilde{Q},\tilde{Q}^\dagger)\in\mathscr{D}^{\alpha,\beta}_{Z^1}([a,b],\mathcal{W})$ and $(\bar{Q},\bar{Q}^\dagger)\in\mathscr{D}^{\alpha,\beta}_{Z^1}([b,c],\mathcal{W})$ with $(\tilde{Q}_b,\tilde{Q}^\dagger_b)=(\bar{Q}_b,\bar{Q}_b^\dagger)$, their concatenation $(Q,Q^\dagger):=(\tilde{Q}*\bar{Q},\tilde{Q}^\dagger*\bar{Q}^\dagger)$ can naturally be defined and belongs to 
		$\mathscr{D}^{\alpha,\beta}_{Z^1}([a,c],\mathcal{W})$. (Here, $*$ stands for the usual concatenation operation for two continuous paths.) It is clear that $Q^\dagger\in\mathcal{C}^\beta([a,c],\mathcal{L}(\mathcal{V},\mathcal{W}))$. To prove that $Q^\sharp\in\mathcal{C}^{\alpha\wedge2\beta}_2([a,c],\mathcal{W})$, it is sufficient to observe the following: For $a\leq s \leq b \leq t\leq c$, 
		\begin{align*}
			Q^\sharp_{s,t}
			=&Q^1_{s,t}-Q^\dagger_s Z^1_{s,t}\\
			=&(\bar{Q}^1_{b,t}-\bar{Q}^\dagger_bZ^1_{b,t})+(\tilde{Q}^1_{s,b}-\tilde{Q}^\dagger_sZ^1_{s,b})+\tilde{Q}^{\dagger,1}_{s,b}Z^1_{b,t}\\
			=& \bar{Q}^\sharp_{b,t}+\tilde{Q}^\sharp_{s,b}+\tilde{Q}^{\dagger,1}_{s,b}Z^1_{b,t}.
		\end{align*}
		The right hand side is clearly dominated by a constant multiple of $(t-s)^{\alpha\wedge2\beta}$. 
	\end{itemize}
\end{example}

Let
$
(Q,Q^\dagger)\in \mathscr{D}^{\alpha,\beta}_{Z^1}\big(\mathcal{L}(\mathcal{V},\mathcal{W})\big),
$
with $\mathbb{Z}=(Z^1,Z^2)\in\Omega^\beta(\mathcal{V})$.

Let $0\leq a<b\leq T$.
For 
$(Q,Q^\dagger)\in \mathscr{D}^{\alpha,\beta}_{Z^1}(\mathcal{L}(\mathcal{V},\mathcal{W}))$, define
$$J_{s,t}=Q_sZ^1_{s,t}+Q^\dagger_{s}Z^2_{s,t}, \quad (s,t)\in\Delta_{[a,b]}.$$

For $a\leq s\leq u\leq t\leq b$, we compute 
\begin{align*}
	\delta J_{s,u,t}:=&J_{s,t}-J_{s,u}-J_{u,t}\\
	=&Q_sZ^1_{s,t}+Q^\dagger_sZ^2_{s,t}
	-Q_sZ^1_{s,u}-Q^\dagger_sZ^2_{s,u}-Q_uZ^1_{u,t}-Q^\dagger_uZ^2_{u,t}\\
	=&-Q^1_{s,u}Z^1_{u,t}-Q^\dagger_sZ^2_{s,u}-Q^\dagger_uZ^2_{u,t}+Q^\dagger_s(Z^2_{s,u}+Z^2_{u,t}+Z^1_{s,u}\otimes Z^1_{u,t})\\
	=&Q^\dagger_sZ^1_{s,u}\otimes Z^1_{u,t}
	-(Q^\dagger_s Z^1_{s,u}+Q^\sharp_{s,u})Z^1_{u,t}-Q^{\dagger,1}_{s,u}Z^2_{u,t}\\
	=&-Q^{\dagger,1}_{s,u}Z^2_{u,t}-Q^\sharp_{s,u}Z^1_{u,t}, 
\end{align*}
where $Q^{\dagger,1}_{s,u}=Q^\dagger_u-Q^\dagger_s$.

Let $\mathcal{P}=\{s=t_0<t_1<\dots<t_N=t\}$ be a partition of $[s,t]\subset[a,b]$. Its mesh size is denoted by $|\mathcal{P}|$. Define $J_{s,t}(\mathcal{P}):=\sum_{i=1}^N J_{t_{i-1},t_i}$.
We define the integral by
\begin{align}\label{eq-1}
	\int_s^t Q_ud\mathbb{Z}_u:=\lim_{|\mathcal{P}|\downarrow 0}J_{s,t}(\mathcal{P}),\quad(s,t)\in\Delta_{[a,b]}
\end{align}
if it exists.


\begin{proposition}\label{pro-1}
Let $[a,b]\subset[0,T]$. 
Suppose that $\mathbb{Z}=(Z^1,Z^2)\in\Omega^\beta([a,b],\mathcal{V})$ and 
$(Q,Q^\dagger)\in\mathscr{D}^{\alpha,\beta}_{Z^1}\big([a,b],\mathcal{L}(\mathcal{V},\mathcal{W})\big)$.
Then, for all $(s,t)\in\Delta_{[a,b]}$, the limit in \eqref{eq-1} exists and defines a generalized rough integral.
Moreover, we have
\[
\Big(\int_a^\cdot Q_u\,d\mathbb{Z}_u,\, Q\Big)
\in \hat{\mathscr{D}}^\beta_{Z^1}\big([a,b],\mathcal{W}\big),
\]
and the estimate
\begin{align}\label{eq-5}
\lefteqn{\Big|
	\int_s^t Q_u\,d\mathbb{Z}_u
	- Q_s Z^1_{s,t}
	- Q^\dagger_s Z^2_{s,t}
	\Big|}\notag\\
&\leq
C
\big(
\|Q^\sharp\|_{\alpha\wedge2\beta;[a,b]}\|Z^1\|_{\beta;[a,b]}
+
\|Q^\dagger\|_{\beta;[a,b]}\|Z^2\|_{2\beta;[a,b]}
\big)
(t-s)^{(\alpha+\beta)\wedge 3\beta}
\end{align}
holds for all $(s,t)\in\Delta_{[a,b]}$. Here, the constant $C>0$ is independent of $s$, $t$, $(Q,Q^\dagger)$ and $\mathbb{Z}$.
\end{proposition}

\begin{proof}
\leavevmode\par

\medskip
\textbf{Step 1: Convergence of the Riemann sums.}
\medskip
	
Let $\mathcal{P}=\{s=t_0<\cdots<t_N=t\}$ be a partition of $[s,t]$.
Then there exists $i\in\{1,\dots,N-1\}$ such that
\[
t_{i+1}-t_{i-1} \le \frac{2(t-s)}{N-1}.
\]
Using the algebraic identity for $J$, we obtain
\begin{align*}
\big|J_{s,t}(\mathcal{P}) - J_{s,t}(\mathcal{P}\setminus\{t_i\})\big| 
&= |J_{t_{i-1},t_i} + J_{t_i,t_{i+1}} - J_{t_{i-1},t_{i+1}}| \\
&= \big| Q^\sharp_{t_{i-1},t_i} Z^1_{t_i,t_{i+1}} 
+ Q^{\dagger,1}_{t_{i-1},t_i} Z^2_{t_i,t_{i+1}} \big|.
\end{align*}
By the H\"{o}lder regularity of $Q^\dagger$, $Q^\sharp$, $Z^1$, and $Z^2$, we deduce
\begin{align*}
\lefteqn{\big|J_{s,t}(\mathcal{P}) - J_{s,t}(\mathcal{P}\setminus\{t_i\})\big|}\\
&\leq
\big(
\|Q^\sharp\|_{\alpha\wedge2\beta;[a,b]}\|Z^1\|_{\beta;[a,b]}
+
\|Q^\dagger\|_{\beta;[a,b]}\|Z^2\|_{2\beta;[a,b]}
\big) \Big(\frac{2(t-s)}{N-1}\Big)^{(\alpha+\beta)\wedge 3\beta}.
\end{align*}
Iterating this procedure and summing over all removed points yields
\begin{align*}
\lefteqn{|J_{s,t}(\mathcal{P}) - J_{s,t}|}\\
&\leq
2^{(\alpha+\beta)\wedge3\beta}\zeta\big((\alpha+\beta)\wedge3\beta\big)\big(
\|Q^\sharp\|_{\alpha\wedge2\beta;[a,b]}\|Z^1\|_{\beta;[a,b]}
+
\|Q^\dagger\|_{\beta;[a,b]}\|Z^2\|_{2\beta;[a,b]}
\big)
(t-s)^{(\alpha+\beta)\wedge 3\beta},
\end{align*}
which proves \eqref{eq-5}. Here, $\zeta$ is Riemann zeta function.
Since $(\alpha+\beta)\wedge 3\beta > 1$, it follows that 
$\big\{J_{s,t}(\mathcal{P})\big\}_{\mathcal{P}}$ is a Cauchy family as $|\mathcal{P}|\to 0$, and hence the limit exists.
	
\medskip
	
\textbf{Step 2: Controlled structure.}

\medskip
	
We show that the rough integral defines a controlled path.
From \eqref{eq-5}, we have
\[
\Big|\int_s^t Q_u\,d\mathbb{Z}_u\Big|
\le |Q_s Z^1_{s,t}| + |Q^\dagger_s Z^2_{s,t}| + C|t-s|^{(\alpha+\beta)\wedge 3\beta},
\]
where the constant $C>0$ depends only on $\alpha$ and $\beta$. Hence
\[
\int_a^\cdot Q_u\,d\mathbb{Z}_u \in \mathcal{C}^\beta\big([a,b],\mathcal{W}\big).
\]
Similarly, splitting up the left-hand side of \eqref{eq-5} after the second term yields a $2\beta$-H\"older type estimate for
$$
\int_s^t Q_u\,d\mathbb{Z}_u - Q_s Z^1_{s,t}
=:
\Big(\int_a^\cdot
Q_u\,d\mathbb{Z}_u\Big)^\sharp_{s,t},
$$
that is,
for the remainder term in the sense of \eqref{eq-22}. 
Thus the remainder is $(2\beta)$-H\"older continuous.
Then, we conclude that
$$
\Big(\int_a^\cdot Q_u\,d\mathbb{Z}_u,\, Q\Big)
\in \hat{\mathscr{D}}^\beta_{Z^1}\big([a,b],  \mathcal{W}\big).
$$
This completes the proof.
\end{proof}

\begin{corollary}\label{co-1}
	Suppose that  $\sigma:~\bar{\mathcal{V}}\times \mathcal{W}\to \mathcal{L}(\mathcal{V},\mathcal{W})$ is of $C^2$. Let $X\in\mathcal{C}^\alpha(\bar{\mathcal{V}})$, and $(Y,Y^\dagger)\in \hat{\mathscr{D}}^\beta_{Z^1}(\mathcal{W})$, where $\mathbb{Z}=(Z^1,Z^2)\in\Omega^\beta(\mathcal{V})$.   
	Then,
	$\big(\int_0^\cdot \sigma(X_s,Y_s)d\mathbb{Z}_s,\sigma(X,Y))\in\hat{\mathscr{D}}^\beta_{Z^1}(\mathcal{W}\big)$.
\end{corollary}
\begin{proof}
	By Example 4, we have  $\big(\sigma(X,Y),\sigma(X,Y)^\dagger\big)\in   \mathscr{D}^{\alpha,\beta}_{Z^1}\big(\mathcal{L}(\mathcal{V},\mathcal{W})\big)$. Therefore, Proposition \ref{pro-1} implies that $\big(\int_0^\cdot \sigma(X_s,Y_s)d\mathbb{Z}_s,\sigma(X,Y)\big)\in\hat{\mathscr{D}}^\beta_{Z^1}(\mathcal{W})$,
	which completes the proof.
\end{proof}

\section{Well-Posedness of Deterministic Young--Rough Systems}\label{sec-3}

Let $\alpha\in(1/2,1]$ and \(\beta\in(1/3,1/2]\) satisfy $\alpha+\beta>1$.
We introduce the following continuous maps:
\begin{align*}
\begin{cases}
b_1:[0,T]\times \mathbb{R}^{d_1}\times \mathbb{R}^{d_2}\to\mathbb{R}^{d_1},
\\
\sigma_1:[0,T]\times \mathbb{R}^{d_1}\times \mathbb{R}^{d_2}\to\mathcal{L}(\mathbb{R}^m, \mathbb{R}^{d_1}),\\
b_2: \mathbb{R}^{d_1}\times \mathbb{R}^{d_2}\to\mathbb{R}^{d_2},
\\
\sigma_2:\mathbb{R}^{d_1}\times \mathbb{R}^{d_2}\to\mathcal{L}(\mathbb{R}^r, \mathbb{R}^{d_2}).
\end{cases}
\end{align*}
For $(B,\mathbb{W})=(B,W^1,W^2)\in\mathcal{C}^\alpha([0,T],\mathbb{R}^m)\times \Omega^\beta([0,T],\mathbb{R}^r)$,
we consider the following Young-rough system:
\begin{align}\label{eq-2}
	\begin{cases}
		dX_t=b_1(t,X_t,Y_t)dt+ \sigma_1(t,X_t,Y_t)dB_t,\\  
		dY_t= b_2(X_t,Y_t)dt+ \sigma_2(X_t,Y_t)d\mathbb{W}_t,\\
		(X_0,Y_0)=(x_0,y_0)\in \mathbb{R}^{d_1}\times \mathbb{R}^{d_2},\\
		Y^\dagger_t=\sigma_2(X_t,Y_t),~ t\in[0,T].
	\end{cases}
\end{align}
Here, the first equation is understood in the Young sense, while the second equation is interpreted as a RDE.
We will solve this system in $\mathcal{C}^\alpha([0,T],\mathbb{R}^{d_1})\times \hat{\mathscr{D}}^\beta_{W^1}([0,T],\mathbb{R}^{d_2})$.
To establish the existence and uniqueness of solutions to system \eqref{eq-2}, we impose the following assumptions.

\paragraph{Assumptions.}
\begin{itemize}
\item[\textbf{(A1)}] The function $\sigma_1(t,\cdot,\cdot)\in C^2_b(\mathbb{R}^{d_1}\times \mathbb{R}^{d_2}, \mathcal{L}(\mathbb{R}^m,\mathbb{R}^{d_1}))$  in space, for any fixed $t\in[0,T]$, 
and is $\beta$-H\"{o}lder continuous in time. More precisely,
\begin{align}\label{eq-23}
S_k:=\sup_{t\in[0,T]}\Big\{\sum_{i=0}^{k}\|\nabla^i\sigma_1(t,\cdot,\cdot)\|_\infty\Big\}<\infty,
\end{align} 
for $k=0,1,2$,
and
there exist constants 
$L_1> 0$ such that for all $(x,y)\in \mathbb{R}^{d_1}\times\mathbb{R}^{d_2}$ and $s,t\in[0,T]$,
\begin{align*}
|\nabla\sigma_1(t,x,y) - \nabla\sigma_1(s,x,y)| 
+ |\sigma_1(t,x,y) - \sigma_1(s,x,y)| 
&\le L_1 |t-s|^{\beta}.
\end{align*}
No differentiability with respect to the time variable is assumed.

\item [\textbf{(A2)}]
The function $\sigma_2\in C^3_b(\mathbb{R}^{d_1}\times \mathbb{R}^{d_2}, \mathcal{L}(\mathbb{R}^r,\mathbb{R}^{d_2}))$.
	
\item [\textbf{(A3)}]
The drift coefficients $b_1$ and $b_2$ are bounded and 
Lipschitz continuous in the space variables uniformly in time. That is, 
there exist constants $L_2, ~L_3 > 0$ such that for all $(x_1,y_1),(x_2,y_2)\in \mathbb{R}^{d_1}\times\mathbb{R}^{d_2}$ and $t\in[0,T]$,
	\begin{align*}
		|b_1(t,x_1,y_1) - b_1(t,x_2,y_2)| 
		+ |b_2(x_1,y_1) - b_2(x_2,y_2)| 
		&\le L_2\big(|x_1-x_2| + |y_1-y_2|\big),\\
		|b_1(t,x_1,y_1)| + |b_2(x_1,y_1)| &\le L_3.
	\end{align*}
\end{itemize}

\subsection{Young-rough differential equations: The case of bounded and globally Lipschitz drift vector field}\label{sec-3.1}

\begin{theorem}\label{th-1}
Assume \textbf{(A1)-(A3)}. Then, for every initial condition
$
(x_0,y_0)\in\mathbb{R}^{d_1}\times\mathbb{R}^{d_2}
$,
there exists a unique global solution 
$$(X,Y,Y^\dagger)\in \mathcal{C}^\alpha([0,T],\mathbb{R}^{d_1})\times\hat{\mathscr{D}}^\beta_{W^1}([0,T],\mathbb{R}^{d_2})$$
of system \eqref{eq-2}. 		
\end{theorem}

\begin{proof}
Throughout, we take $\alpha^\prime$ and $\beta^\prime$ sufficiently close to $\alpha$ and $\beta$, respectively, such that $\frac13<\beta^\prime<\beta\leq \frac12<\alpha^\prime<\alpha<1$
and $\alpha^\prime+\beta^\prime>1$.
For simplicity, we assume $\alpha-\alpha^\prime=\beta-\beta^\prime$.
Define the unit ball
\begin{align*}
	\mathcal{B}^{x_0,y_0}_{[0,T]}
	=\Big\{&(X,Y,Y^\dagger)\in\mathcal{C}^{\alpha^\prime}([0,T],\mathbb{R}^{d_1})\times\hat{\mathscr{D}}^{\beta^\prime}_{W^1}([0,T],\mathbb{R}^{d_2})|\\
	&
	\|X\|_{\alpha^\prime}+\|(Y,Y^\dagger)\|_{\hat{\mathscr{D}}^{\beta^\prime}_{W^1};[0,\tau]}\leq 1,~ (X_0,Y_0)=(x_0,y_0),~
	Y^\dagger_0=\sigma_2(x_0,y_0)	\Big\}
\end{align*} 
If we restrict the seminorm to $\mathcal{B}^{x_0,y_0}_{[0,T]}$, it naturally induces a metric, and the resulting metric space is complete.
Since $\mathcal{C}^\alpha([0,T],\mathbb{R}^{d_1})$ and  $\hat{\mathscr{D}}^{\beta^\prime}_{W^1}([0,T],\mathbb{R}^{d_2})$ are Banach spaces with norms
\begin{align*}
X&\mapsto |X_0| +\|X\|_\alpha,\\
(Y,Y^\dagger)&\mapsto |Y_0|+|Y^\dagger_0|+\|(Y,Y^\dagger)\|_{\hat{\mathscr{D}}^{\beta^\prime}_{W^1}},
\end{align*} 
the set $\mathcal{B}^{x_0,y_0}_{[0,T]}$ is a closed ball centered at 
$$t\mapsto \big(x_0,y_0+\sigma_2(x_0,y_0)W^1_{0,t},\sigma_2(x_0,y_0)\big).$$ 
Define maps $$\mathcal{M}^{x_0,y_0}_{[0,\tau]}, \mathcal{M}^1_{[0,\tau]}, \mathcal{M}^2_{[0,\tau]}:\mathcal{C}^{\alpha}([0,T],\mathbb{R}^{d_1})\times \hat{\mathscr{D}}^{\beta}_{W^1}([0,T],\mathbb{R}^{d_2})\to \mathcal{C}^{\alpha}([0,T],\mathbb{R}^{d_1})\times \hat{\mathscr{D}}^{\beta}_{W^1}([0,T],\mathbb{R}^{d_2}),$$
where $\tau\in(0,1]$ and $(x_0,y_0)\in\mathbb{R}^{d_1}\times\mathbb{R}^{d_2}$,
\begin{align*}
	\mathcal{M}^1_{[0,\tau]}(X,Y,Y^\dagger)&=\Big(\int_0^\cdot b_1(s,X_s,Y_s)ds,\int_0^\cdot b_2(X_s,Y_s)ds,0\Big),\\
	\mathcal{M}^2_{[0,\tau]}(X,Y,Y^\dagger)&=\Big(\int_0^\cdot \sigma_1(s,X_s,Y_s)dB_s,\int_0^\cdot \sigma_2(X_s,Y_s)d\mathbb{W}_s,\sigma_2(X,Y)\Big),
\end{align*}
and 
$\mathcal{M}^{x_0,y_0}_{[0,\tau]}=(x_0,y_0,0)+\mathcal{M}^1_{[0,\tau]}+\mathcal{M}^2_{[0,\tau]}$.
Here, the item $\big(\int_0^\cdot \sigma_2(X_s,Y_s)d\mathbb{W}_s,\sigma_2(X,Y)\big)$ is a controlled path discussed in Corollary \ref{co-1}. 
We equip space $\mathcal{C}^{\alpha}([0,T],\mathbb{R}^{d_1})\times \hat{\mathscr{D}}^{\beta}_{W^1}([0,T],\mathbb{R}^{d_2})$ with the seminorm $\|(X,Y,Y^\dagger)\|_{\alpha,\beta,W^1}:=\|X\|_\alpha+\|(Y,Y^\dagger)\|_{\hat{\mathscr{D}}^\beta_{W^1}}$.
Strictly speaking, this would only yield a solution $(X,Y,Y^\dagger)$ in $\mathcal{C}^{\alpha^\prime}\times \hat{\mathscr{D}}^{\beta^\prime}_{W^1}$. But since $(B,\mathbb{W})\in\mathcal{C}^{\alpha}\times \Omega^{\beta}$, it turns out that solution is automatically an element of $\mathcal{C}^\alpha\times \hat{\mathscr{D}}^\beta_{W^1}$. Indeed, 
$$|X_{s,t}|\leq \|b_1\|_\infty|t-s|+C\|\sigma_1(\cdot,X_\cdot,Y_\cdot)\|_{\beta^\prime}\|B\|_\alpha|t-s|^\alpha
\leq \|b_1\|_\infty|t-s|+C\|B\|_\alpha|t-s|^\alpha,$$ 
and 
$$|Y_{s,t}|\leq \|Y^\dagger\|_\infty|W^1_{s,t}|+\|R^Y\|_{2\beta^\prime}|t-s|^{2\beta^\prime},$$
so that $X\in\mathcal{C}^\alpha$ and $Y\in\mathcal{C}^\beta$.
Here, the first inequality follows from the estimate for the Young integral in \eqref{eq-57}, which will be explained later.
It then follows that $Y^\dagger=\sigma_2(X,Y)\in\mathcal{C}^\beta$ and also $R^Y\in\mathcal{C}^{2\beta}_2$, since $W^2\in\mathcal{C}^{2\beta}_2$ and 
\begin{align*}
|R^Y_{s,t}|=&|Y_{s,t}-Y^\dagger_sW^1_{s,t}|=\Big|\int_s^t b_2(X_u,Y_u)du+
\int_s^t\sigma_2(X_u,Y_u)d\mathbb{W}_u-\sigma_2(X_s,Y_s)W^1_{s,t}\Big|\\
\leq &
\|b_2\|_\infty|t-s|+
\|
\sigma_2(X,Y)^\dagger\|_\infty|W^2_{s,t}|+O(|t-s|^{3\beta^\prime})\\
\leq &
\|b_2\|_\infty|t-s|+C\|Y^\dagger\|_\infty|W^2_{s,t}|+O(|t-s|^{3\beta^\prime}).
\end{align*}
Here, the second and the last inequalities follow from \eqref{eq-5} and Example 4, respectively.

To simplify the notation, we suppress the dependence on the interval $[0,\tau]$ whenever it is clear from the context. In particular,  
$\|\cdot\|_{\alpha^\prime}$
and$\|\cdot\|_{\infty}$ will denote the H\"older norm $\|\cdot\|_{\alpha^\prime;[0,\tau]}$
and  the uniform norm $\|\cdot\|_{\infty;[0,\tau]}$, respectively.

A fixed point of map $\mathcal{M}^{x_0,y_0}_{[0,\tau]}$ is a solution of system \eqref{eq-2} on $[0,\tau]$.

\medskip
\textbf{Step 1: Invariance:}
\medskip

By a straightforward calculation, we obtain
$$\|\mathcal{M}^{x_0,y_0}_{[0,\tau]}(X,Y,Y^\dagger)\|_{\alpha^\prime,\beta^\prime,W^1}
\leq \|\mathcal{M}^1_{[0,\tau]}(X,Y,Y^\dagger)\|_{\alpha^\prime,\beta^\prime,W^1}+\|\mathcal{M}^2_{[0,\tau]}(X,Y,Y^\dagger)\|_{\alpha^\prime,\beta^\prime,W^1}.$$

We first estimate the drift part. By the boundedness of $b_1$ and $b_2$, we have 
\begin{align*}
\|\mathcal{M}^1_{[0,\tau]}(X,Y,Y^\dagger)\|_{\alpha^\prime,\beta^\prime,W^1}
=&\Big\|\int_0^\cdot b_1(s,X_s,Y_s)ds\Big\|_{\alpha^\prime}+\Big\|\int_0^\cdot b_2(X_s,Y_s)ds\Big\|_{2\beta^\prime}\\
\leq& \sup_{0\leq s<t\leq \tau}\frac{|\int_s^tb_1(u,X_u,Y_u)du|}{|t-s|^{\alpha^\prime}}
+\sup_{0\leq s<t\leq \tau}\frac{|\int_s^tb_2(X_u,Y_u)du|}{|t-s|^{2\beta^\prime}}\\
\leq &\|b_1\|_\infty\tau^{1-\alpha^\prime}+\|b_2\|_\infty\tau^{1-2\beta^\prime}\\
\leq & L_3\tau^{\beta-\beta^\prime}.
\end{align*}
Here, we used the fact that $\alpha-\alpha^\prime=\beta-\beta^\prime$ is sufficiently small.

For $(X,Y,Y^\dagger)\in\mathcal{B}^{x_0,y_0}_{[0,\tau]}$, the following estimates hold:
\begin{align}\label{eq-3}
\|Y^\dagger\|_\infty
\leq |\sigma_2(x_0,y_0)|+\|Y^\dagger\|_{\beta^\prime}\tau^{\beta^\prime}
\leq \|\sigma_2\|_\infty+1, 
\end{align}
and
\begin{align}\label{eq-4}
|Y^1_{s,t}|
\leq& |Y^\dagger_sW^1_{s,t}|+|Y^\sharp_{s,t}|\notag\\
\leq &\big(\|\sigma_2\|_\infty+1\big)\|W^1\|_\beta|t-s|^\beta+\|Y^\sharp\|_{2\beta^\prime}|t-s|^{2\beta^\prime}\notag\\
\leq &\big(\|\sigma_2\|_\infty+1\big)\big(\|W^1\|_\beta+1\big)|t-s|^\beta.
\end{align}

Next, we estimate the Young-rough integrals part $\mathcal{M}^2_{[0,\tau]}$. By definition,
\begin{align*}
\|\mathcal{M}^2_{[0,\tau]}(X,Y,Y^\dagger)\|_{\alpha^\prime,\beta^\prime,W^1}
=&\Big\|\int_0^\cdot \sigma_1(s,X_s,Y_s)dB_s\Big\|_{\alpha^\prime}
+\Big\|\Big(\int_0^\cdot\sigma_2(X_s,Y_s)d\mathbb{W}_s\Big)^\sharp\Big\|_{2\beta^\prime}
+\|\sigma_2(X,Y)\|_{\beta^\prime}. 
\end{align*}

In order to estimate the Young integral term, recall the standard estimate for Young integrals. 
Define Young integral $Z:=\int_0^\cdot \tilde{Y}_sd\tilde{X}_s$ with $\tilde{X}\in\mathcal{C}^{\tilde{\alpha}}([0,\tau],\mathbb{R}^{d_1})$ and $\tilde{Y}\in\mathcal{C}^{\tilde{\beta}}([0,\tau],\mathcal{L}(\mathbb{R}^{d_1},\mathbb{R}^{d_2}))$, where $\tilde{\alpha}+\tilde{\beta}>1$, for $0\leq s \leq t\leq \tau$.
From \cite[P. 62 Eq. (4.3)]{friz2014course}, we have
\begin{align}\label{eq-14}
|Z_{s,t}|
\leq&
|\tilde{Y}_s||\tilde{X}^1_{s,t}|+C\|\tilde{Y}\|_{{\tilde{\beta}};[s,t]} \|\tilde{X}\|_{{\tilde{\alpha}};[s,t]}|t-s|^{{\tilde{\alpha}}+{\tilde{\beta}}}\notag\\
\leq &
\big(|\tilde{Y}_0|+\|\tilde{Y}\|_{{\tilde{\beta}};[0,\tau]}\big)\|\tilde{X}\|_{{\tilde{\alpha}};[0,\tau]}|t-s|^{\tilde{\alpha}}
+
C\|\tilde{Y}\|_{{\tilde{\beta}};[0,\tau]}\|\tilde{X}\|_{{\tilde{\alpha}};[0,\tau]}|t-s|^{\tilde{\alpha}}\notag\\
\leq &
C\big(|\tilde{Y}_0|+\|\tilde{Y}\|_{{\tilde{\beta}};[0,\tau]}\big)\|\tilde{X}\|_{{\tilde{\alpha}};[0,\tau]}|t-s|^{\tilde{\alpha}}.
\end{align}

By \eqref{eq-4} and \eqref{eq-14} together with Assumption \textbf{(A1)}, we obtain
\begin{align}\label{eq-57}
\Big\|\int_0^\cdot \sigma_1(s,X_s,Y_s)dB_s\Big\|_{\alpha^\prime}
\leq& C\big(|\sigma_1(0,x_0,y_0)|+\|\sigma_1(\cdot,X_\cdot,Y_\cdot)\|_{\beta^\prime}\big)\|B\|_{\alpha^\prime}\notag\\
\leq&C\big\{S_0 +
L_1|t-s|^{\beta-\beta^\prime}+S_1\big(\|X\|_{\beta^\prime}+\|Y\|_{\beta^\prime}\big)\big\}\|B\|_\alpha\tau^{\alpha-\alpha^\prime}\notag\\
\leq& C\big(S_1+L_1\big)\big(\|\sigma_2\|_\infty+1\big)\big(\|W^1\|_\beta+1\big)\|B\|_\alpha\tau^{\alpha-\alpha^\prime}.
\end{align}

We now turn to the estimate of $\|\sigma_2(X,Y)\|_{\beta^\prime}$. Using \eqref{eq-4} and the boundedness of $\nabla\sigma_2$, we have 
\begin{align*}
\|\sigma_2(X,Y)\|_{\beta^\prime}
\leq &\sup_{0\leq s<t\leq \tau}\frac{|\sigma_2(X_t,Y_t)-\sigma_2(X_s,Y_s)|}{|t-s|^{\beta^\prime}}\\
\leq & \|\nabla \sigma_2\|_\infty \sup_{0\leq s<t\leq \tau}\frac{|X_t-X_s|+|Y_t-Y_s|}{|t-s|^{\beta^\prime}}\\
\leq &\|\sigma_2\|_{C^1_b}\big(\|X\|_{\beta^\prime}+\|Y\|_{\beta^\prime}\big)\\
\leq &\|\sigma_2\|_{C^1_b}\tau^{\alpha^\prime-\beta^\prime}+\|\sigma_2\|_{C^1_b}\big(\|\sigma_2\|_\infty+1\big)\big(\|W^1\|_\beta+1\big)\tau^{\beta-\beta^\prime}\\
\leq &C\big(\|\sigma_2\|_{C^1_b}+1\big)^2\big(\|W^1\|_\beta+1\big)\tau^{\beta-\beta^\prime}.
\end{align*}

By Proposition \ref{pro-1} and the standard estimate for rough integrals, we have
\begin{align}\label{eq-13}
\Big\|\Big(\int_0^\cdot\sigma_2(X_s,Y_s)d\mathbb{W}_s\Big)^\sharp\Big\|_{2\beta^\prime}
=&\sup_{0\leq s<t\leq \tau}\frac{|\int_s^t\sigma_2(X_s,Y_s)d\mathbb{W}_s-\sigma_2(X_s,Y_s)W^1_{s,t}|}{|t-s|^{2\beta^\prime}}\notag\\
\leq &\|\sigma_2(X,Y)^\dagger\|_\infty\|W^2\|_{2\beta^\prime}\notag\\
&+C\big(\|\sigma_2(X,Y)^\dagger\|_{\beta^\prime}\|W^2\|_{2\beta^\prime}+\|\sigma_2(X,Y)^\sharp\|_{\alpha^\prime\wedge2\beta^\prime}\|W^1\|_{\beta^\prime}\big)\tau^{(\alpha^\prime-\beta^\prime)\wedge\beta^\prime}.
\end{align}
Therefore, it is sufficient to estimate  $\big\|\big(\sigma_2(X,Y),\sigma_2(X,Y)^\dagger\big)\big\|_{\mathscr{D}^{\alpha^\prime,\beta^\prime}_{W^1}}$.

Let $\gamma^\prime:=\alpha^\prime\wedge2\beta^\prime$.
From Lemma \ref{lem-1} and \eqref{eq-4}, we have
\begin{align*}
\|\sigma_2(X,Y)^\dagger\|_{\beta^\prime}
=& \|\nabla_y \sigma_2(X,Y)Y^\dagger\|_{\beta^\prime}\\
\leq & \|\nabla\sigma_2\|_\infty\|Y^\dagger\|_{\beta^\prime}
+\|Y^\dagger\|_\infty\|\nabla^2\sigma_2\|_\infty\big(\|X\|_{\beta^\prime}+\|Y\|_{\beta^\prime}\big)\\
\leq &\|\sigma_2\|_{C^2_b}\big(|Y^\dagger_0|+\|(Y,Y^\dagger)\|_{\hat{\mathscr{D}}^{\beta^\prime}_{W^1}}\big)\big(1+\|X\|_{\beta^\prime}+\|Y\|_{\beta^\prime}\big)\\
\leq & C
\|\sigma_2\|_{C^2_b}\big(\|\sigma_2\|_\infty+1\big)\big\{1+  
\big(\|\sigma_2\|_\infty+1\big)\big(\|W^1\|_\beta+1\big)\big\}\\
\leq &C\big(\|\sigma_2\|_{C^2_b}+1\big)^3\big(\|W^1\|_\beta+1\big).
\end{align*}

For the remainder term, using Taylor expansion, we obtain
\begin{align*}
\sigma_2(X,Y)^\sharp_{s,t}
=& \sigma_2(X_t,Y_t)-\sigma_2(X_s,Y_s)-\sigma_2(X,Y)^\dagger_sW^1_{s,t}\\
=&\sigma_2(X_t,Y_t)-\sigma_2(X_s,Y_s)-\nabla_y\sigma_2(X_s,Y_s)Y^\dagger_sW^1_{s,t}\\
=&\sigma_2(X_t,Y_t)-\sigma_2(X_s,Y_s)-\nabla_y\sigma_2(X_s,Y_s)\big(Y^1_{s,t}-Y^\sharp_{s,t}\big)\\
=& \nabla_x\sigma_2(X_s,Y_s)X^1_{s,t}+\int_0^1(1-\theta) \nabla^2_{xx}\sigma_2(Z^\theta_{s,t})\langle X^1_{s,t},X^1_{s,t}\rangle d\theta\\
&+\int_0^1(1-\theta) \nabla^2_{yy}\sigma_2(Z^\theta_{s,t})\langle Y^1_{s,t},Y^1_{s,t}\rangle d\theta
+2\int_0^1(1-\theta) \nabla^2_{xy}\sigma_2(Z^\theta_{s,t})\langle X^1_{s,t},Y^1_{s,t}\rangle d\theta\\
&+\nabla_y\sigma_2(X_s,Y_s)Y^\sharp_{s,t},
\end{align*}
where $Z^\theta_{s,t}:=(X_s+\theta X^1_{s,t},Y_s+\theta Y^1_{s,t})$.
Consequently, we have 
\begin{align*}
\|\sigma_2(X,Y)^\sharp\|_{\gamma^\prime}
\leq& \|\sigma_2\|_{C^1_b}\|X\|_{\gamma^\prime}
+\frac12\|\sigma_2\|_{C^2_b}\sup_{0\leq s<t\leq \tau}\frac{(X^1_{s,t}+Y^1_{s,t})^2}{|t-s|^{\gamma^\prime}}
+\|\sigma_2\|_{C^1_b}\|Y^\sharp\|_{\gamma^\prime}\\
\leq &\|\sigma_2\|_{C^1_b}\big(\|X\|_{\alpha^\prime}+\|Y^\sharp\|_{2\beta^\prime}\big)
+\frac12\|\sigma_2\|_{C^2_b}\big(\|X\|_{\frac{\gamma^\prime}2}+\|Y\|_{\frac{\gamma^\prime}2}\big)^2\\
\leq & \|\sigma_2\|_{C^1_b}
+\frac12\|\sigma_2\|_{C^2_b}\big(\|X\|_{\alpha^\prime}+\|Y\|_{\beta^\prime}\big)^2\\
\leq &\|\sigma_2\|_{C^1_b}
+\frac12\|\sigma_2\|_{C^2_b}\big\{1+  \big(\|\sigma_2\|_\infty+1\big)\big(\|W^1\|_\beta+1\big)\big\}^2\\
\leq &C  \big(\|\sigma_2\|_{C^2_b}+1\big)^3\big(\|W^1\|_\beta+1\big)^2.
\end{align*}

Hence, we have
\begin{align*}
\big\|\big(\sigma_2(X,Y),\sigma_2(X,Y)^\dagger\big)\big\|_{\mathscr{D}^{\alpha^\prime,\beta^\prime}_{W^1}}
\leq &C  \big(\|\sigma_2\|_{C^2_b}+1\big)^3\big(\|W^1\|_\beta+1\big)^2.
\end{align*}

Substituting the above estimates into the rough integral bound \eqref{eq-13}, we finally obtain
\begin{align}\label{eq-15}
\lefteqn{\Big\|\Big(\int_0^\cdot\sigma_2(X_s,Y_s)d\mathbb{W}_s\Big)^\sharp\Big\|_{2\beta^\prime}}\notag\\
&\leq \|\sigma_2(X,Y)^\dagger\|_\infty\|W^2\|_{2\beta^\prime}
+C\big(\|\sigma_2(X,Y)^\dagger\|_{\beta^\prime}\|W^2\|_{2\beta^\prime}+\|\sigma_2(X,Y)^\sharp\|_{\gamma^\prime}\|W^1\|_{\beta^\prime}\big)\notag\\
&\leq
\big(|\sigma_2(X,Y)^\dagger_0|+\|\sigma_2(X,Y)^\dagger\|_{\beta^\prime}\big)\|W^2\|_{2\beta^\prime}\notag\\
&\quad+C\big\|\big(\sigma_2(X,Y),\sigma_2(X,Y)^\dagger\big)\big\|_{\mathscr{D}^{\alpha^\prime,\beta^\prime}_{W^1}}\big(\|W^1\|_{\beta^\prime}+\|W^2\|_{2\beta^\prime}\big)\notag\\
&\leq C\big(|\nabla_y\sigma_2(x_0,y_0)\sigma_2(x_0,y_0)|+\big\|\big(\sigma_2(X,Y),\sigma_2(X,Y)^\dagger\big)\big\|_{\mathscr{D}^{\alpha^\prime,\beta^\prime}_{W^1}}\big)\big(\|W^1\|_\beta+\|W^2\|_{2\beta}\big)\tau^{\beta-\beta^\prime}\notag\\
&\leq C\big(\|\sigma_2\|_{C^2_b}+1\big)^3\big(\|W^1\|_\beta+1\big)^2\big(\|W^1\|_\beta+\|W^2\|_{2\beta}\big)\tau^{\beta-\beta^\prime}.
\end{align}
Here, $\nabla_y\sigma_2(x_0,y_0)\sigma_2(x_0,y_0)=\nabla_y\sigma_2(x_0,y_0)\langle\sigma_2(x_0,y_0)\star,\blacktriangle\rangle$, where $\star,\blacktriangle\in\mathbb{R}^r$,
and $$|\nabla_y\sigma_2(x_0,y_0)\sigma_2(x_0,y_0)|\leq |\nabla_y\sigma_2(x_0,y_0)|\times|\sigma_2(x_0,y_0)|.$$

Collecting all the above estimates, we conclude that there exits a constant $C>0$ such that
$$\mathcal{M}^{x_0,y_0}_{[0,\tau]}(X,Y,Y^\dagger)=C\tau^{\beta-\beta^\prime},$$
where $C$ is independent of the initial value $(x_0,y_0)\in\mathbb{R}^{d_1}\times\mathbb{R}^{d_2}$ and $(X,Y,Y^\dagger)\in\mathcal{B}^{x_0,y_0}_{[0,\tau]}$.

\medskip
\textbf{Step 2: Contraction:}
\medskip

Now, we will prove that $\mathcal{M}^{x_0,y_0}_{[0,\tau]}$ is a contraction on $\mathcal{B}^{x_0,y_0}_{[0,\tau]}$ when $\tau=\tau_0$ is choosen small enough.

Let
$$(X,Y,Y^\dagger)\in\mathcal{B}^{x_0,y_0}_{[0,\tau]},\quad(\tilde{X},\tilde{Y},\tilde{Y}^\dagger)\in\mathcal{B}^{\tilde{x}_0,\tilde{y}_0}_{[0,\tau]}.$$
We first derive estimates for the controlled paths $(Y,Y^\dagger)$ and $(\tilde{Y},\tilde{Y}^\dagger)$:
\begin{align}
\|Y^\dagger-\tilde{Y}^\dagger\|_\infty
\leq&
|Y^\dagger_0-\tilde{Y}^\dagger_0| 
+\|Y^\dagger-\tilde{Y}^\dagger\|_{\beta^\prime}\tau^{\beta^\prime}\notag\\
= &
|\sigma_2(x_0,y_0)-\sigma_2(\tilde{x}_0,\tilde{y}_0)|+\|Y^\dagger-\tilde{Y}^\dagger\|_{\beta^\prime}\tau^{\beta^\prime}\notag\\
\leq &
\|\sigma_2\|_{C^1_b}\big(|x_0-\tilde{x}_0|+|y_0-\tilde{y}_0|\big)+\|Y^\dagger-\tilde{Y}^\dagger\|_{\beta^\prime}\tau^{\beta^\prime},\label{eq-6}\\

|Y^1_{s,t}-\tilde{Y}^1_{s,t}|
\leq &
|Y^\dagger_s-\tilde{Y}^\dagger_s||W^1_{s,t}|+|Y^\sharp_{s,t}-\tilde{Y}^\sharp_{s,t}|\notag\\
\leq & 
\big\{\|\sigma_2\|_{C^1_b}\big(|x_0-\tilde{x}_0|+|y_0-\tilde{y}_0|\big)+\|Y^\dagger-\tilde{Y}^\dagger\|_{\beta^\prime}\tau^{\beta^\prime}\big\}\|W^1\|_\beta|t-s|^\beta\notag\\
&+\|Y^\sharp-\tilde{Y}^\sharp\|_{2\beta^\prime}|t-s|^{2\beta^\prime}
\notag\\
\leq &
\|\sigma_2\|_{C^1_b}\big(\|W^1\|_\beta+1\big)\big(|x_0-\tilde{x}_0|+|y_0-\tilde{y}_0|\big)|t-s|^\beta\notag\\
&+\big(\|W^1\|_\beta+1\big)\|(Y,Y^\dagger)-(\tilde{Y},\tilde{Y}^\dagger)\|_{\hat{\mathscr{D}}^{\beta^\prime}_{W^1}}|t-s|^\beta,\label{eq-7}\\

\|Y-\tilde{Y}\|_\infty
\leq &
|Y_0-\tilde{Y}_0|+\|Y^1_{0,\cdot}-\tilde{Y}^1_{0,\cdot}\|_\infty\notag\\
\leq &
(\|\sigma_2\|_{C^1_b}+1)\big(\|W^1\|_\beta+1\big)\big(|x_0-\tilde{x}_0|+|y_0-\tilde{y}_0|\big)\notag\\
&+\big(\|W^1\|_\beta+1\big)\|(Y,Y^\dagger)-(\tilde{Y},\tilde{Y}^\dagger)\|_{\hat{\mathscr{D}}^{\beta^\prime}_{W^1}}\tau^\beta.\label{eq-8}
\end{align}
For simplicity, let $Z=(X,Y)$ and $\tilde{Z}=(\tilde{X},\tilde{Y})$. Then, we will give some estimates.
From \eqref{eq-4}, we have
\begin{align}\label{eq-19}
	|Z^1_{s,t}|
	=&|X^1_{s,t}|+|Y^1_{s,t}|\notag\\
	\leq &
	\|X\|_{\alpha^\prime}|t-s|^{\alpha^\prime}
	+
	\big(\|\sigma_2\|_\infty+1\big)\big(\|W^1\|_\beta+1\big)|t-s|^\beta\notag\\
	\leq &
	C\big(\|\sigma_2\|_\infty+1\big)\big(\|W^1\|_\beta+1\big)|t-s|^\beta.
\end{align}
From \eqref{eq-7}, we have
\begin{align}\label{eq-12}
	|Z^1_{s,t}-\tilde{Z}^1_{s,t}|
	=&
	|X^1_{s,t}-\tilde{X}^1_{s,t}|+|Y^1_{s,t}-\tilde{Y}^1_{s,t}|\notag\\
	\leq &
	\|X-\tilde{X}\|_{\alpha^\prime}
	|t-s|^{\alpha^\prime}
	+
	\|\sigma_2\|_{C^1_b}\big(\|W^1\|_\beta+1\big)\big(|x_0-\tilde{x}_0|+|y_0-\tilde{y}_0|\big)|t-s|^\beta\notag\\
	&+\big(\|W^1\|_\beta+1\big)\|(Y,Y^\dagger)-(\tilde{Y},\tilde{Y}^\dagger)\|_{\hat{\mathscr{D}}^{\beta^\prime}_{W^1}}|t-s|^\beta\notag\\
	\leq &\|\sigma_2\|_{C^1_b}\big(\|W^1\|_\beta+1\big)\big(|x_0-\tilde{x}_0|+|y_0-\tilde{y}_0|\big)|t-s|^\beta\notag\\
	&+
	\big(\|W^1\|_\beta+1\big)\big(\|X-\tilde{X}\|_{\alpha^\prime}+
	\|(Y,Y^\dagger)-(\tilde{Y},\tilde{Y}^\dagger)\|_{\hat{\mathscr{D}}^{\beta^\prime}_{W^1}}\big)|t-s|^\beta.
\end{align}
From \eqref{eq-8}, we obtain
\begin{align}\label{eq-11}
	\|Z-\tilde{Z}\|_\infty
	=&
	\|X-\tilde{X}\|_\infty+\|Y-\tilde{Y}\|_\infty\notag\\
	\leq &
	|x_0-\tilde{x}_0|
	+
	\|X-\tilde{X}\|_{\alpha^\prime}\tau^{\alpha^\prime}\notag\\
	&+
	(\|\sigma_2\|_{C^1_b}+1)\big(\|W^1\|_\beta+1\big)\big(|x_0-\tilde{x}_0|+|y_0-\tilde{y}_0|\big)\notag\\
	&+\big(\|W^1\|_\beta+1\big)\|(Y,Y^\dagger)-(\tilde{Y},\tilde{Y}^\dagger)\|_{\hat{\mathscr{D}}^{\beta^\prime}_{W^1}}\tau^\beta\notag\\
	\leq &
	2(\|\sigma_2\|_{C^1_b}+1)\big(\|W^1\|_\beta+1\big)\big(|x_0-\tilde{x}_0|+|y_0-\tilde{y}_0|\big)\notag\\
	&+\big(\|W^1\|_\beta+1\big)\big(\|X-\tilde{X}\|_{\alpha^\prime}+
	\|(Y,Y^\dagger)-(\tilde{Y},\tilde{Y}^\dagger)\|_{\hat{\mathscr{D}}^{\beta^\prime}_{W^1}}\big)\tau^\beta.
\end{align}

From Assumption \textbf{(A2)} and \eqref{eq-8}, we have
\begin{align*}
&{\|\mathcal{M}^1_{[0,\tau]}(X,Y,Y^\dagger)-\mathcal{M}^1_{[0,\tau]}(\tilde{X},\tilde{Y},\tilde{Y}^\dagger)\|_{\alpha^\prime,\beta^\prime,W^1}}\\
&= \Big\|\int_0^\cdot\big(b_1(s,Z_s)-b_1(s,\tilde{Z}_s)\big)ds\Big\|_{\alpha^\prime}+\Big\|\Big(\int_0^\cdot\big(b_2(Z_s)-b_2(\tilde{Z}_s)\big)ds,0\Big)\Big\|_{\hat{\mathscr{D}}^{\beta^\prime}_{W^1}}\\
&=
\Big\|\int_0^\cdot\big(b_1(s,Z_s)-b_1(s,\tilde{Z}_s)\big)ds\Big\|_{\alpha^\prime}+\Big\|\int_0^\cdot\big(b_2(Z_s)-b_2(\tilde{Z}_s)\big)ds\Big\|_{2\beta^\prime}\\
&\leq 
L_2\big(\|X-\tilde{X}\|_\infty+\|Y-\tilde{Y}\|_\infty\big)\big(\tau^{1-\alpha^\prime}+\tau^{1-2\beta^\prime}\big)\\
&\leq 
CL_2\big(\|\sigma_2\|_{C^1_b}+1\big)\big(\|W^1\|_\beta+1\big)\big(|x_0-\tilde{x}_0|+|y_0-\tilde{y}_0|\big)\tau^{\beta-\beta^\prime}\notag\\
&\quad+L_2\big(\|W^1\|_\beta+1\big)\big(\|X-\tilde{X}\|_{\alpha^\prime}+\|(Y,Y^\dagger)-(\tilde{Y},\tilde{Y}^\dagger)\|_{\hat{\mathscr{D}}^{\beta^\prime}_{W^1}}\big)\tau^\beta.
\end{align*}

Define $\Delta_t:=\sigma_2(X_t,Y_t)-\sigma_2(\tilde{X}_t,\tilde{Y}_t)$.
From Lemma \ref{lem-1}, $(\Delta,\Delta^\dagger)\in\mathscr{D}^{\alpha^\prime,\beta^\prime}_{W^1}$.
Then, by \eqref{eq-15}, we have 
\begin{align}\label{eq-24}
&\|\mathcal{M}^2_{[0,\tau]}(X,Y,Y^\dagger)-\mathcal{M}^2_{[0,\tau]}(\tilde{X},\tilde{Y},\tilde{Y}^\dagger)\|_{\alpha^\prime,\beta^\prime,W^1}\notag\\
&=
\Big\|\Big(\int_0^\cdot\big(\sigma_1(s,Z_s)-\sigma_1(s,\tilde{Z}_s)\big)dB_s,
\int_0^\cdot\Delta_s d\mathbb{W}_s,
\Delta
\Big)\Big\|_{\alpha^\prime,\beta^\prime,W^1}\notag\\
&=
\Big\|\int_0^\cdot\big(\sigma_1(s,Z_s)-\sigma_1(s,\tilde{Z}_s)\big)dB_s\Big\|_{\alpha^\prime}
+\Big\|\Big(\int_0^\cdot\Delta_sd\mathbb{W}_s,\Delta\Big)\Big\|_{\hat{\mathscr{D}}^{\beta^\prime}_{W^1}}\notag\\
&=:A+\Big\|\Big(\int_0^\cdot\Delta_sd\mathbb{W}_s,\Delta\Big)\Big\|_{\hat{\mathscr{D}}^{\beta^\prime}_{W^1}}\notag\\
&\leq A+\|\Delta\|_{\beta^\prime}
+C\big(|\Delta^\dagger_0|+\|(\Delta,\Delta^\dagger)\|_{\mathscr{D}^{\alpha^\prime,\beta^\prime}_{W^1}}\big)
\big(\|W^1\|_\beta+\|W^2\|_{2\beta}\big)\tau^{\beta-\beta^\prime}.
\end{align}

For $A$,
by the standard estimate of Young integral \eqref{eq-14} and the fact that $\alpha^\prime+\beta^\prime>1$, we have
\begin{align}\label{eq-16}
A
\leq &
C\big(|\sigma_1(0,x_0,y_0)-\sigma_1(0,\tilde{x}_0,\tilde{y}_0)|
+\|\sigma_1(\cdot,Z_\cdot)-\sigma_1(\cdot,\tilde{Z}_\cdot)\|_{\beta^\prime}\big)\|B\|_{\alpha^\prime}\notag\\
\leq &
C\big\{S_1\big(|x_0-\tilde{x}_0|+|y_0-\tilde{y}_0|\big)
+\|\sigma_1(\cdot,Z_\cdot)-\sigma_1(\cdot,\tilde{Z}_\cdot)\|_{\beta^\prime}\big\}\|B\|_{\alpha^\prime}.
\end{align}

To estimate 
$
\|\sigma_1(\cdot,Z_\cdot)-\sigma_1(\cdot,\tilde{Z}_\cdot)\|_{\beta^\prime}
$,
we first establish the following lemma.
\begin{lemma}\label{lem-2}
Let $$\varphi:[0,T]\times \mathcal{W}\times \bar{\mathcal{W}}\to \mathcal{V},\quad0< T\leq 1,$$
satisfies the following conditions:
\begin{itemize}
\item[1.]
The map $\varphi(t,\cdot,\cdot)\in C^2_b(\mathcal{W}\times\bar{\mathcal{W}},\mathcal{V})$ uniformly in $t$, that is 
$$\tilde{S}_k:=\sup_{t\in[0,T]}\Big\{\sum_{i=0}^{k}\|\nabla^i\varphi(t,\cdot,\cdot)\|_\infty\Big\}<\infty,$$ 
for $k=0,1,2$.
No differentiability with respect to the time variable is assumed.

\item[2.]
The function $\varphi$ and $\nabla\varphi$ are $\beta$-H\"older continuous in time, uniformly in $(x,y)\in\mathcal{W}\times\bar{\mathcal{W}}$, that is, there exists a constant $C>0$ such that
$$|\varphi(t,x,y)-\varphi(s,x,y)|+|\nabla\varphi(t,x,y)-\nabla  \varphi(s,x,y)|\leq C|t-s|^\beta$$
for all $s,t\in[0,T]$ and $(x,y)\in\mathcal{W}\times\bar{\mathcal{W}}$.
\end{itemize}
Let $X,\tilde{X}\in\mathcal{C}^\alpha([0,T],\mathcal{W})$ and $Y,\tilde{Y}\in\mathcal{C}^\beta([0,T],\bar{\mathcal{W}})$. 
Assume that there exist constants $K_1$ and $K_2$ such that
$$\|X\|_{\alpha;[0,T]}+\|\tilde{X}\|_{\alpha;[0,T]}\leq K_1,\quad
\|Y\|_{\beta;[0,T]}+\|\tilde{Y}\|_{\beta;[0,T]}\leq K_2.$$
Then, there exists a constant $C>0$, depending only on $K_1$ and $K_2$, such that
$$\|\varphi(\cdot,X_\cdot,Y_\cdot)-\varphi(\cdot,\tilde{X}_\cdot,\tilde{Y}_\cdot)\|_{\beta;[0,T]}
\leq C\tilde{S}_2\big(|X_0-\tilde{X}_0|+|Y_0-\tilde{Y}_0|+\|X-\tilde{X}\|_\alpha+\|Y-\tilde{Y}\|_\beta\big).
$$
\end{lemma}

\begin{proof}
For $0\leq s<t\leq T$, we decompose the increment as
\begin{align*}
&\varphi(\cdot,X_\cdot,Y_\cdot)^1_{s,t}-\varphi(\cdot,\tilde{X}_\cdot,\tilde{Y}_\cdot)^1_{s,t}\\
=&
\varphi(t,X_t,Y_t)-\varphi(s,X_t,Y_t)
-\big(\varphi(t,\tilde{X}_t,\tilde{Y}_t)-\varphi(s,\tilde{X}_t,\tilde{Y}_t)\big)
\\
&+\varphi(s,X_t,Y_t)-\varphi(s,X_s,Y_s)
-\big(\varphi(s,\tilde{X}_t,\tilde{Y}_t)-\varphi(s,\tilde{X}_s,\tilde{Y}_s)\big).
\end{align*} 

We first estimate the time increment.
Define
$$f_{s,t}(x,y):=\varphi(t,x,y)-\varphi(s,x,y).$$
Then, it suffices to estimate 
$|f_{s,t}(X_t,Y_t)-f_{s,t}(\tilde{X}_t,\tilde{Y}_t)|.$

The first condition implies that $f_{s,t}\in C^2_b$.
Applying the division property, we have
\begin{align*}
f_{s,t}(X_t,Y_t)-f_{s,t}(\tilde{X}_t,\tilde{Y}_t)
=\int_0^1
\left[
\begin{array}{c|c}
	\nabla_x f_{s,t}(\gamma^\theta_t) & \nabla_y f_{s,t}(\gamma^\theta_t)
\end{array}
\right]
\begin{bmatrix}
	X_t-\tilde{X}_t\\
	Y_t-\tilde{Y}_t
\end{bmatrix}
d \theta,
\end{align*}
where $$\gamma^\theta_t:=\big(\tilde{X}_t+\theta(X_t-\tilde{X}_t),\tilde{Y}_t+\theta(Y_t-\tilde{Y}_t)\big).$$
From the second condition, we obtain 
$$\|\nabla f_{s,t}\|_\infty\leq C|t-s|^\beta,$$
where $\nabla f_{s,t}=\nabla\varphi(t,\cdot,\cdot)-\nabla\varphi(s,\cdot,\cdot)$.
Consequently,
\begin{align*}
|f_{s,t}(X_t,Y_t)-f_{s,t}(\tilde{X}_t,\tilde{Y}_t)|
\leq& C\big(\|X-\tilde{X}\|_\infty+\|Y-\tilde{Y}\|_\infty\big)|t-s|^\beta\\
\leq &
C\big(|X_0-\tilde{X}_0|+|Y_0-\tilde{Y}_0|+\|X-\tilde{X}\|_\alpha+\|Y-\tilde{Y}\|_\beta\big)|t-s|^\beta.
\end{align*}

Next, we estimate the spatial increment. By the mean value theorem, we have
\begin{align*}
E:=&\varphi(s,X_t,Y_t)-\varphi(s,X_s,Y_s)
-\big(\varphi(s,\tilde{X}_t,\tilde{Y}_t)-\varphi(s,\tilde{X}_s,\tilde{Y}_s)\big)\\
=&
\int_0^1 \nabla_x\varphi(s,\gamma^\theta_t)\langle X_t-\tilde{X}_t\rangle d\theta
+\int_0^1 \nabla_y\varphi(s,\gamma^\theta_t)\langle Y_t-\tilde{Y}_t\rangle d\theta\\
&-\int_0^1 \nabla_x\varphi(s,\gamma^\theta_s)\langle X_s-\tilde{X}_s\rangle d\theta
-\int_0^1 \nabla_y\varphi(s,\gamma^\theta_s)\langle Y_s-\tilde{Y}_s\rangle d\theta.
\end{align*}
Define 
$$g_s(x,\tilde{x},y,\tilde{y}):=
\int_0^1 \nabla\varphi(s,\tilde{x}+\theta(x-\tilde{x}),
\tilde{y}+\theta(y-\tilde{y}))d\theta,$$
and $$\Xi_t:=(X_t-\tilde{X}_t,Y_t-\tilde{Y}_t).$$
Then, we have 
\begin{align*}
E
=
g_s(X_t,\tilde{X}_t,Y_t,\tilde{Y}_t)\Xi_t
-g_s(X_s,\tilde{X}_s,Y_s,\tilde{Y}_s)\Xi_s.
\end{align*}

Since $\varphi(t,\cdot,\cdot)\in C^2_b$ for every $t\in[0,T]$, the function $g_s$ satisfies the following conditions: 
\begin{itemize}
\item[(1)]Boundedness: $$\sup_{s\in[0,T]}\|g_s\|_\infty\leq \tilde{S}_1,$$
\item[(2)]Lipschitz continuity: $$\sup_{s\in[0,T]}\|g_s\|_{\mathrm{Lip}}\leq C\tilde{S}_2.$$
\end{itemize} 
Therefore,
\begin{align*}
|E|
\leq 
\sup_{s\in[0,T]}\{\|g_s\|_\infty\}|\Xi_t-\Xi_s|
+\sup_{s\in[0,T]}\{\|g_s\|_{\mathrm{Lip}}\}\big|\big(X^1_{s,t},\tilde{X}^1_{s,t},Y^1_{s,t},\tilde{Y}^1_{s,t}\big)\big||\Xi_s|.
\end{align*}
Moreover,
\begin{align*}
|\Xi_s|
\leq
\|X-\tilde{X}\|_\infty+\|Y-\tilde{Y}\|_\infty
\leq \big(|X_0-\tilde{X}_0|+|Y_0-\tilde{Y}_0|+\|X-\tilde{X}\|_\alpha T^\alpha+\|Y-\tilde{Y}\|_\beta T^\beta\big),
\end{align*}
and
\begin{align*}
|\Xi_t-\Xi_s|
\leq 
\|X-\tilde{X}\|_\alpha|t-s|^\alpha+\|Y-\tilde{Y}\|_\beta|t-s|^\beta
\leq 
\big(\|X-\tilde{X}\|_\alpha+\|Y-\tilde{Y}\|_\beta\big)|t-s|^\beta.
\end{align*}
On the other hand,
\begin{align*}
	|X^1_{s,t}|
	\leq \|X\|_\alpha|t-s|^\alpha
	\leq K_1|t-s|^\alpha,
\end{align*}
and similarly,
\begin{align*}
	|\tilde{X}_{s,t}|
	&\leq
	K_1|t-s|^\alpha,\\
	|Y_{s,t}|
	&\leq
	K_2|t-s|^\beta,\\
	|\tilde{Y}_{s,t}|
	&\leq
	K_2|t-s|^\beta.
\end{align*}

Combining the above estimates, we conclude that
\begin{align*}
\|\varphi(\cdot,X_\cdot,Y_\cdot)-\varphi(\cdot,\tilde{X}_\cdot,\tilde{Y}_\cdot)\|_{\beta;[0,T]}
\leq C\tilde{S}_2\big(|X_0-\tilde{X}_0|+|Y_0-\tilde{Y}_0|+\|X-\tilde{X}\|_\alpha+\|Y-\tilde{Y}\|_\beta\big).
\end{align*}
This completes the proof of Lemma \ref{lem-2}.
\end{proof}

We now continue with the proof of Theorem \ref{th-1}.

By Lemma \ref{lem-2} and \eqref{eq-7}, we obtain
\begin{align*}
\lefteqn{\|\sigma_1(\cdot,Z_\cdot)-\sigma_1(\cdot,\tilde{Z}_\cdot)\|_{\beta^\prime}}\\
&\leq 
CS_2\big(|x_0-\tilde{x}_0|+|y_0-\tilde{y}_0|+\|X-\tilde{X}\|_{\alpha^\prime}+\|Y-\tilde{Y}\|_{\beta^\prime}\big)\\
&\leq 
CS_2\big\{|x_0-\tilde{x}_0|+|y_0-\tilde{y}_0|
+
\|\sigma_2\|_{C^1_b}\big(\|W^1\|_\beta+1\big)\big(|x_0-\tilde{x}_0|+|y_0-\tilde{y}_0|\big)\\
&\quad+\|X-\tilde{X}\|_{\alpha^\prime}+\big(\|W^1\|_\beta+1\big)\|(Y,Y^\dagger)-(\tilde{Y},\tilde{Y}^\dagger)\|_{\hat{\mathscr{D}}^{\beta^\prime}_{W^1}}
\big\}\\
&\leq 
CS_2
(\|\sigma_2\|_{C^1_b}+1)
\big(\|W^1\|_\beta+1\big)
\big(
|x_0-\tilde{x}_0|+|y_0-\tilde{y}_0|
\big)\\
&\quad+
CS_2
\big(\|W^1\|_\beta+1\big)
\big(
\|X-\tilde{X}\|_{\alpha^\prime}
+
\|(Y,Y^\dagger)-(\tilde{Y},\tilde{Y}^\dagger)\|_{\hat{\mathscr{D}}^{\beta^\prime}_{W^1}}
\big).
\end{align*}
Plugging the above inequality into \eqref{eq-16}, we obtain
\begin{align}\label{eq-70}
A
\leq &
CS_2
(\|\sigma_2\|_{C^1_b}+1)
\big(\|W^1\|_\beta+1\big)
\big(|x_0-\tilde{x}_0|+|y_0-\tilde{y}_0|\big)
\|B\|_\alpha\tau^{\alpha-\alpha^\prime}\notag\\
&+
CS_2
\big(\|W^1\|_\beta+1\big)
\big(
\|X-\tilde{X}\|_{\alpha^\prime}
+
\|(Y,Y^\dagger)-(\tilde{Y},\tilde{Y}^\dagger)\|_{\hat{\mathscr{D}}^{\beta^\prime}_{W^1}}
\big)
\|B\|_\alpha\tau^{\alpha-\alpha^\prime}
.
\end{align}

For the term $\|\Delta\|_{\beta^\prime}$, since $\Delta\in\mathscr{D}^{\alpha^\prime,\beta^\prime}_{W^1}$, similar to the calculation for $\|Y\|_{\beta^\prime}$ \eqref{eq-4}, we have 
\begin{align}\label{eq-25}
\|\Delta\|_{\beta^\prime}
\leq&
\|\Delta^\dagger\|_\infty\|W^1\|_{\beta^\prime} +\|\Delta^\sharp\|_{\beta^\prime}\notag\\
\leq &
\big(|\Delta^\dagger_0|\|W^1\|_\beta
+\|\Delta^\dagger\|_{\beta^\prime}\|W^1\|_\beta+\|\Delta^\sharp\|_{\gamma^\prime}\big)\tau^{\beta-\beta^\prime}\notag\\
\leq &
\big(\|W^1\|_\beta+1\big)\big(|\Delta^\dagger_0|+\|(\Delta,\Delta^\dagger)\|_{\mathscr{D}^{\alpha^\prime,\beta^\prime}_{W^1}}\big)\tau^{\beta-\beta^\prime},
\end{align}
where $\gamma^\prime=\alpha^\prime\wedge2\beta^\prime$.

For the term $|\Delta^\dagger_0|$,
\begin{align}\label{eq-26}
|\Delta^\dagger_0|
=&\big|\nabla_y\sigma_2(x_0,y_0)\sigma_2(x_0,y_0)-\nabla_y\sigma_2(\tilde{x}_0,\tilde{y}_0)\sigma_2(\tilde{x}_0,\tilde{y}_0)\big|\notag\\
\leq &
\big|\nabla_y\sigma_2(x_0,y_0)\big\{\sigma_2(x_0,y_0)-\sigma_2(\tilde{x}_0,\tilde{y}_0)\big\}\big|
+\big|\big\{\nabla_y\sigma_2(x_0,y_0)-\nabla_y\sigma_2(\tilde{x}_0,\tilde{y}_0)\big\}\sigma_2(\tilde{x}_0,\tilde{y}_0)\big|
\notag\\
\leq &
C\|\sigma_2\|^2_{C^2_b}\big(|x_0-\tilde{x}_0|+|y_0-\tilde{y}_0|\big).
\end{align}

For the term
$\|(\Delta,\Delta^\dagger)\|_{\mathscr{D}^{\alpha^\prime,\beta^\prime}_{W^1}}$,
we divide the estimate into two parts, i.e. $\|\Delta^\dagger\|_{\beta^\prime}$ and $\|\Delta^\sharp\|_{\gamma^\prime}$.
For the first part $\|\Delta^\dagger\|_{\beta^\prime}$, we first derive an auxiliary estimate
\begin{align}\label{eq-9}
\sigma_2(Z)^{\dagger,1}_{s,t}
=&\nabla_y\sigma_2(Z_t)Y^\dagger_t-\nabla_y\sigma_2(Z_s)Y^\dagger_s\notag\\
= &
\big\{\nabla_y\sigma_2(Z_t)-\nabla_y\sigma_2(Z_s)\big\}Y^\dagger_t
+\nabla_y\sigma_2(Z_s)Y^{\dagger,1}_{s,t}\notag\\
=&
\int_0^1\nabla\big\{\nabla_y\sigma_2(Z_s+\theta Z^1_{s,t})\big\}
\langle Z^1_{s,t},Y^\dagger_t\bullet,\star\rangle d\theta 
+
\nabla_y\sigma_2(Z_s)Y^{\dagger,1}_{s,t}\notag\\
=&
\int_0^1\big\{\nabla^2_{xy}\sigma_2(Z_s+\theta Z^1_{s,t})\langle X^1_{s,t},Y^\dagger_t\bullet,\star\rangle +\nabla^2_{yy}\sigma_2(Z_s+\theta Z^1_{s,t})\langle Y^1_{s,t},Y^\dagger_t\bullet,\star\rangle\big\}d\theta\notag\\
&+
\nabla_y\sigma_2(Z_s)Y^{\dagger,1}_{s,t},
\end{align}
where $\bullet,\star\in\mathbb{R}^r$.
In the sequel, unless otherwise specified, the symbols
$\bullet$ and $\star$ always denote arbitrary elements of $\mathbb{R}^r$.
Then, from \eqref{eq-9}, we have
\begin{align*}
|\Delta^{\dagger,1}_{s,t}|&=|\sigma_2(Z)^{\dagger,1}_{s,t}-\sigma_2(\tilde{Z})^{\dagger,1}_{s,t}|\\
&\leq \Big|
\int_0^1\nabla^2_{xy}\sigma_2(Z_s+\theta Z^1_{s,t})\langle X^1_{s,t},Y^\dagger_t\bullet,\star\rangle d\theta
-
\int_0^1\nabla^2_{xy}\sigma_2(\tilde{Z}_s+\theta \tilde{Z}^1_{s,t})\langle \tilde{X}^1_{s,t},\tilde{Y}^\dagger_t\bullet,\star\rangle d\theta
\Big|\\
&\quad+
\Big|
\int_0^1\nabla^2_{yy}\sigma_2(Z_s+\theta Z^1_{s,t})\langle Y^1_{s,t},Y^\dagger_t\bullet,\star\rangle d\theta
-
\int_0^1\nabla^2_{yy}\sigma_2(\tilde{Z}_s+\theta \tilde{Z}^1_{s,t})\langle \tilde{Y}^1_{s,t},\tilde{Y}^\dagger_t\bullet,\star\rangle d\theta
\Big|\\
&\quad+
\big|\nabla_y\sigma_2(Z_s)\big(Y^{\dagger,1}_{s,t}
-
\tilde{Y}^{\dagger,1}_{s,t}\big)
\big|\\
&\quad
+
\big|\big\{\nabla_y\sigma_2(Z_s)
-
\nabla_y\sigma_2(\tilde{Z}_s)\big\} \tilde{Y}^{\dagger,1}_{s,t}
\big|
\\
&=:B_1+B_2+B_3+B_4.
\end{align*}


From \eqref{eq-3} and \eqref{eq-6}, we have 
\begin{align}\label{eq-17}
B_1
\leq &
\Big|
\int_0^1\nabla^2_{xy}\sigma_2(Z_s+\theta Z^1_{s,t})\langle X^1_{s,t}-\tilde{X}^1_{s,t},Y^\dagger_t\bullet,\star\rangle d\theta
\Big|\notag\\
&+
\Big|
\int_0^1\nabla^2_{xy}\sigma_2(\tilde{Z}_s+\theta \tilde{Z}^1_{s,t})\langle \tilde{X}^1_{s,t},(Y^\dagger_t-\tilde{Y}^\dagger_t)\bullet,\star\rangle d\theta
\Big|\notag\\
&+
\Big|
\int_0^1\big\{\nabla^2_{xy}\sigma_2(Z_s+\theta Z^1_{s,t})-\nabla^2_{xy}\sigma_2(\tilde{Z}_s+\theta \tilde{Z}^1_{s,t})\big\}\langle \tilde{X}^1_{s,t},\tilde{Y}^\dagger_t\bullet,\star\rangle d\theta
\Big|\notag\\
\leq &
C\|\sigma_2\|_{C^2_b}\|X-\tilde{X}\|_{\alpha^\prime}\|Y^\dagger\|_\infty|t-s|^{\alpha^\prime}
+
C\|\sigma_2\|_{C^2_b}\|\tilde{X}\|_{\alpha^\prime}\|Y^\dagger-\tilde{Y}^\dagger\|_\infty|t-s|^{\alpha^\prime}\notag\\
&+
C\|\sigma_2\|_{C^3_b}\|\tilde{X}\|_{\alpha^\prime}\|\tilde{Y}^\dagger\|_\infty|t-s|^{\alpha^\prime}
\int_0^1 |Z_s+\theta Z^1_{s,t}-\tilde{Z}_s-\theta \tilde{Z}^1_{s,t}|d\theta\notag\\
\leq &
C\|\sigma_2\|_{C^2_b}\big(\|\sigma_2\|_\infty+1\big)\|X-\tilde{X}\|_{\alpha^\prime}|t-s|^{\alpha^\prime}\notag\\
&+
C\|\sigma_2\|_{C^2_b}
\big\{\|\sigma_2\|_{C^1_b}\big(|x_0-\tilde{x}_0|+|y_0-\tilde{y}_0|\big)+\|Y^\dagger-\tilde{Y}^\dagger\|_{\beta^\prime}\tau^{\beta^\prime}\big\}
|t-s|^{\alpha^\prime}\notag\\
&+
C\|\sigma_2\|_{C^3_b}\big(\|\sigma_2\|_\infty+1\big)
\|Z-\tilde{Z}\|_\infty
|t-s|^{\alpha^\prime}
. 
\end{align}
Plugging \eqref{eq-11}  into \eqref{eq-17}, we have
\begin{align*}
B_1
\leq &
C\big(\|\sigma_2\|_{C^3_b}+1\big)^3
\big(\|W^1\|_\beta+1\big)
\big(|x_0-\tilde{x}_0|+|y_0-\tilde{y}_0|\big)|t-s|^{\alpha^\prime}\\
&+
C\big(\|\sigma_2\|_{C^3_b}+1\big)^2
\big(\|W^1\|_\beta+1\big)\big(\|X-\tilde{X}\|_{\alpha^\prime}+\|(Y,Y^\dagger)-(\tilde{Y},\tilde{Y}^\dagger)\|_{\hat{\mathscr{D}}^{\beta^\prime}_{W^1}}\big)|t-s|^{\alpha^\prime}.
\end{align*}

For the term $B_2$, by \eqref{eq-3}, \eqref{eq-4}, \eqref{eq-6}, \eqref{eq-7} and \eqref{eq-11}, we have
\begin{align*}
B_2
\leq &
\Big|
\int_0^1\nabla^2_{yy}\sigma_2(Z_s+\theta Z^1_{s,t})\langle Y^1_{s,t}-\tilde{Y}^1_{s,t},Y^\dagger_t\bullet,\star\rangle d\theta
\Big|\\
&+
\Big|
\int_0^1\nabla^2_{yy}\sigma_2(\tilde{Z}_s+\theta \tilde{Z}^1_{s,t})\langle \tilde{Y}^1_{s,t},(Y^\dagger_t-\tilde{Y}^\dagger_t)\bullet,\star\rangle d\theta
\Big|\\
&+
\Big|
\int_0^1\big\{\nabla^2_{yy}\sigma_2(Z_s+\theta Z^1_{s,t})-\nabla^2_{yy}\sigma_2(\tilde{Z}_s+\theta \tilde{Z}^1_{s,t})\big\}\langle \tilde{Y}^1_{s,t},\tilde{Y}^\dagger_t\bullet,\star\rangle d\theta
\Big|\\
\leq&
\|\sigma_2\|_{C^2_b}|Y^1_{s,t}-\tilde{Y}^1_{s,t}|\|Y^\dagger\|_\infty\\
&+
\|\sigma_2\|_{C^2_b}|\tilde{Y}^1_{s,t}|\|Y^\dagger-\tilde{Y}^\dagger\|_\infty\\
&+
C\|\sigma_2\|_{C^3_b}|\tilde{Y}^1_{s,t}|\|\tilde{Y}^\dagger\|_\infty
\|Z-\tilde{Z}\|_\infty
\\
\leq &
\|\sigma_2\|_{C^2_b}\big(\|\sigma_2\|_\infty+1\big)\big(\|W^1\|_\beta+1\big)\big\{\|\sigma_2\|_{C^1_b}\big(|x_0-\tilde{x}_0|+|y_0-\tilde{y}_0|\big)+\|(Y,Y^\dagger)-(\tilde{Y},\tilde{Y}^\dagger)\|_{\hat{\mathscr{D}}^{\beta^\prime}_{W^1}}\big\}|t-s|^\beta\\
&+
\|\sigma_2\|_{C^2_b}\big(\|\sigma_2\|_\infty+1\big)\big(\|W^1\|_\beta+1\big)
\big\{\|\sigma_2\|_{C^1_b}\big(|x_0-\tilde{x}_0|+|y_0-\tilde{y}_0|\big)+\|Y^\dagger-\tilde{Y}^\dagger\|_{\beta^\prime}\big\}|t-s|^\beta\\
&+C\|\sigma_2\|_{C^3_b}
\big(\|\sigma_2\|_\infty+1\big)^2\big(\|W^1\|_\beta+1\big)|t-s|^\beta\\
&\times
\Big\{\big(\|\sigma_2\|_{C^1_b}+1\big)\big(\|W^1\|_\beta+1\big)\big(|x_0-\tilde{x}_0|+|y_0-\tilde{y}_0|\big)
\notag\\
&+
\big(\|W^1\|_\beta+1\big)\big(\|X-\tilde{X}\|_{\alpha^\prime}+\|(Y,Y^\dagger)-(\tilde{Y},\tilde{Y}^\dagger)\|_{\hat{\mathscr{D}}^{\beta^\prime}_{W^1}}\big)
\Big\}\notag\\
\leq &
C \big(\|\sigma_2\|_{C^3_b}+1\big)^4
\big(\|W^1\|_\beta+1\big)^2
\big(|x_0-\tilde{x}_0|+|y_0-\tilde{y}_0|\big)|t-s|^{\beta}\\
&+
C \big(\|\sigma_2\|_{C^3_b}+1\big)^3
\big(\|W^1\|_\beta+1\big)^2
\big(\|X-\tilde{X}\|_{\alpha^\prime}+\|(Y,Y^\dagger)-(\tilde{Y},\tilde{Y}^\dagger)\|_{\hat{\mathscr{D}}^{\beta^\prime}_{W^1}}\big)
|t-s|^{\beta}.
\end{align*}

For the term $B_3$,
\begin{align*}
B_3
\leq
C\|\sigma_2\|_{C^1_b}\|Y^\dagger-\tilde{Y}^\dagger\|_{\beta^\prime}|t-s|^{\beta^\prime}. 
\end{align*}

For the term $B_4$, by \eqref{eq-11}, we have
\begin{align*}
B_4
\leq &
C\|\sigma_2\|_{C^2_b}\|Z-\tilde{Z}\|_\infty\|\tilde{Y}^\dagger\|_{\beta^\prime}|t-s|^{\beta^\prime}\\
\leq &
C\|\sigma_2\|_{C^2_b}(\|\sigma_2\|_{C^1_b}+1)
\big(\|W^1\|_\beta+1\big)
\big(|x_0-\tilde{x}_0|+|y_0-\tilde{y}_0|\big)|t-s|^{\beta^\prime}\notag\\
&+
C\|\sigma_2\|_{C^2_b}\big(\|W^1\|_\beta+1\big)\big(\|X-\tilde{X}\|_{\alpha^\prime}+
\|(Y,Y^\dagger)-(\tilde{Y},\tilde{Y}^\dagger)\|_{\hat{\mathscr{D}}^{\beta^\prime}_{W^1}}\big)
|t-s|^{\beta^\prime}.
\end{align*}

In summary, it can be conclude that
\begin{align*}
|\Delta^{\dagger,1}_{s,t}|
&\leq 
C \big(\|\sigma_2\|_{C^3_b}+1\big)^4
\big(\|W^1\|_\beta+1\big)^2
\big(|x_0-\tilde{x}_0|+|y_0-\tilde{y}_0|\big)|t-s|^{\beta^\prime}\\
&\quad+
C \big(\|\sigma_2\|_{C^3_b}+1\big)^3
\big(\|W^1\|_\beta+1\big)^2
\big(\|X-\tilde{X}\|_{\alpha^\prime}+\|(Y,Y^\dagger)-(\tilde{Y},\tilde{Y}^\dagger)\|_{\hat{\mathscr{D}}^{\beta^\prime}_{W^1}}\big)
|t-s|^{\beta^\prime}.
\end{align*}
Consequently, we obtain
\begin{align}\label{eq-20}
\|\Delta^\dagger\|_{\beta^\prime}=&\|\sigma_2(Z)^\dagger-\sigma_2(\tilde{Z})^\dagger\|_{\beta^\prime}\notag\\
\leq &
C \big(\|\sigma_2\|_{C^3_b}+1\big)^4
\big(\|W^1\|_\beta+1\big)^2
\big(|x_0-\tilde{x}_0|+|y_0-\tilde{y}_0|\big)\notag\\
&+
C \big(\|\sigma_2\|_{C^3_b}+1\big)^3
\big(\|W^1\|_\beta+1\big)^2
\big(\|X-\tilde{X}\|_{\alpha^\prime}+\|(Y,Y^\dagger)-(\tilde{Y},\tilde{Y}^\dagger)\|_{\hat{\mathscr{D}}^{\beta^\prime}_{W^1}}\big).
\end{align}

Before estimating the second term$\|\Delta^\sharp\|_{\gamma^\prime}$, we first derive a Taylor expansion of
\(\sigma_2(Z)^\sharp_{s,t}\) as follows:
\begin{align}\label{eq-10}
	\sigma_2(Z)^\sharp_{s,t}
	=&
	\sigma_2(Z)^1_{s,t}-\sigma_2(Z)^\dagger_s W^1_{s,t}\notag\\
	=&
	\nabla_x\sigma_2(Z_s)\langle X^1_{s,t},\bullet \rangle 
	+
	\nabla_y\sigma_2(Z_s)\langle Y^\sharp_{s,t},\bullet \rangle
	+ 
	\int_0^1(1-\theta)\nabla^2\sigma_2(Z_s+\theta Z^1_{s,t})
	\langle Z^1_{s,t},Z^1_{s,t},\bullet \rangle d\theta,
\end{align}
where we used the identities $Y^1_{s,t}=Y^\dagger_s W^1_{s,t}=\sigma_2(Z)^\dagger_sW^1_{s,t}$ and $\nabla_y\sigma_2(Z)^\dagger_s=\nabla_y\sigma_2(Z_s)Y^\dagger_s$.
Then, by \eqref{eq-10}, we have
\begin{align*}
|\sigma_2(Z)^\sharp_{s,t}-\sigma_2(\tilde{Z})^\sharp_{s,t}|
\leq &
\big|\nabla_x\sigma_2(Z_s)\langle X^1_{s,t},\bullet\rangle-
\nabla_x\sigma_2(\tilde{Z}_s)\langle \tilde{X}^1_{s,t},\bullet\rangle\big|\\
&+
\big|\nabla_y\sigma_2(Z_s)\langle Y^\sharp_{s,t},\bullet\rangle-
\nabla_y\sigma_2(\tilde{Z}_s)\langle \tilde{Y}^\sharp_{s,t},\bullet\rangle\big|\\
&+
\Big|
\int_0^1(1-\theta)\nabla^2\sigma_2(Z_s+\theta Z^1_{s,t})\langle Z^1_{s,t},Z^1_{s,t},\bullet\rangle d\theta\\
&\qquad-
\int_0^1(1-\theta)\nabla^2\sigma_2(\tilde{Z}_s+\theta \tilde{Z}^1_{s,t})\langle \tilde{Z}^1_{s,t},\tilde{Z}^1_{s,t},\bullet\rangle d\theta
\Big|\\
=:& E_1+E_2+E_3.
\end{align*}
For the term $E_1$, by \eqref{eq-11}, we have
\begin{align*}
E_1
\leq &
\|\sigma_2\|_{C^1_b}\|X-\tilde{X}\|_{\alpha^\prime}|t-s|^{\alpha^\prime}
+\|\sigma_2\|_{C^2_b}\|Z-\tilde{Z}\|_\infty\|\tilde{X}\|_{\alpha^\prime}|t-s|^{\alpha^\prime}\\
\leq &
\|\sigma_2\|_{C^1_b}\|X-\tilde{X}\|_{\alpha^\prime}|t-s|^{\alpha^\prime}\\
&+
C\|\sigma_2\|_{C^2_b}
(\|\sigma_2\|_{C^1_b}+1)\big(\|W^1\|_\beta+1\big)\big(|x_0-\tilde{x}_0|+|y_0-\tilde{y}_0|\big)
|t-s|^{\alpha^\prime}
\notag\\
&+
\|\sigma_2\|_{C^2_b}
\big(\|W^1\|_\beta+1\big)\big(\|X-\tilde{X}\|_{\alpha^\prime}+
\|(Y,Y^\dagger)-(\tilde{Y},\tilde{Y}^\dagger)\|_{\hat{\mathscr{D}}^{\beta^\prime}_{W^1}}\big)
|t-s|^{\alpha^\prime}
\notag\\
\leq &
C\|\sigma_2\|_{C^2_b}
(\|\sigma_2\|_{C^1_b}+1)\big(\|W^1\|_\beta+1\big)\big(|x_0-\tilde{x}_0|+|y_0-\tilde{y}_0|\big)
|t-s|^{\alpha^\prime}\\
&+C\|\sigma_2\|_{C^2_b}\big(\|W^1\|_\beta+1\big)\big(\|X-\tilde{X}\|_{\alpha^\prime}+
\|(Y,Y^\dagger)-(\tilde{Y},\tilde{Y}^\dagger)\|_{\hat{\mathscr{D}}^{\beta^\prime}_{W^1}}\big)|t-s|^{\alpha^\prime}.
\end{align*}
For the term $E_2$, by \eqref{eq-11}, we have
\begin{align*}
E_2
\leq &
\big|\nabla_y\sigma_2(Z_s)\langle Y^\sharp_{s,t}-\tilde{Y}^\sharp_{s,t},\bullet\rangle\big|
+
\big|\big\{\nabla_y\sigma_2(Z_s)-\nabla_y\sigma_2(\tilde{Z}_s)\big\}\langle \tilde{Y}^\sharp_{s,t},\bullet\rangle\big|\\
\leq &
\|\sigma_2\|_{C^1_b}\|Y^\sharp-\tilde{Y}^\sharp\|_{2\beta^\prime}|t-s|^{2\beta^\prime}
+
\|\sigma_2\|_{C^2_b}\|Z-\tilde{Z}\|_\infty\|\tilde{Y}^\sharp\|_{2\beta^\prime}|t-s|^{2\beta^\prime}\\
\leq &
C\|\sigma_2\|_{C^2_b}
(\|\sigma_2\|_{C^1_b}+1)\big(\|W^1\|_\beta+1\big)\big(|x_0-\tilde{x}_0|+|y_0-\tilde{y}_0|\big)|t-s|^{2\beta^\prime}
\\
&+
C\|\sigma_2\|_{C^2_b}\big(\|W^1\|_\beta+1\big)\big(\|X-\tilde{X}\|_{\alpha^\prime}+
\|(Y,Y^\dagger)-(\tilde{Y},\tilde{Y}^\dagger)\|_{\hat{\mathscr{D}}^{\beta^\prime}_{W^1}}\big)
|t-s|^{2\beta^\prime}.
\end{align*}
For the term $E_3$, by \eqref{eq-19}, \eqref{eq-12} and \eqref{eq-11}, we have
\begin{align*}
E_3
\leq &
\int_0^1\big|\nabla^2\sigma_2(Z_s+\theta Z^1_{s,t})\langle Z^1_{s,t}-\tilde{Z}^1_{s,t},Z^1_{s,t},\bullet\rangle\big|d\theta\\
&+
\int_0^1\big|\nabla^2\sigma_2(Z_s+\theta Z^1_{s,t})\langle \tilde{Z}^1_{s,t},Z^1_{s,t}-\tilde{Z}^1_{s,t},\bullet\rangle\big|d\theta\\
&+
\int_0^1\big|\big\{\nabla^2\sigma_2(Z_s+\theta Z^1_{s,t})-\nabla^2\sigma_2(\tilde{Z}_s+\theta \tilde{Z}^1_{s,t})
\big\}\langle \tilde{Z}^1_{s,t},\tilde{Z}^1_{s,t},\bullet\rangle\big|d\theta\\
\leq &
\|\sigma_2\|_{C^2_b}|Z^1_{s,t}-\tilde{Z}^1_{s,t}||Z^1_{s,t}|
+
\|\sigma_2\|_{C^2_b}|\tilde{Z}^1_{s,t}||Z^1_{s,t}-\tilde{Z}^1_{s,t}|
+
C\|\sigma_2\|_{C^3_b}\|Z-\tilde{Z}\|_\infty|\tilde{Z}^1_{s,t}|^2\\
\leq &
C\|\sigma_2\|_{C^2_b}\|\sigma_2\|_{C^1_b}\big(\|\sigma_2\|_\infty+1\big)\big(\|W^1\|_\beta+1\big)^2\big(|x_0-\tilde{x}_0|+|y_0-\tilde{y}_0|\big)|t-s|^{2\beta}\notag\\
&+
C\|\sigma_2\|_{C^2_b}\big(\|\sigma_2\|_\infty+1\big)\big(\|W^1\|_\beta+1\big)^2\big(\|X-\tilde{X}\|_{\alpha^\prime}+
\|(Y,Y^\dagger)-(\tilde{Y},\tilde{Y}^\dagger)\|_{\hat{\mathscr{D}}^{\beta^\prime}_{W^1}}\big)|t-s|^{2\beta}\\
&+
C\|\sigma_2\|_{C^3_b}\big(\|\sigma_2\|_{C^1_b}+1\big)^3
\big(\|W^1\|_\beta+1\big)^3
\big(|x_0-\tilde{x}_0|+|y_0-\tilde{y}_0|\big)|t-s|^{2\beta}
\\
&+
C\|\sigma_2\|_{C^3_b}
\big(\|\sigma_2\|_\infty+1\big)^2
\big(\|W^1\|_\beta+1\big)^3\big(\|X-\tilde{X}\|_{\alpha^\prime}+\|(Y,Y^\dagger)-(\tilde{Y},\tilde{Y}^\dagger)\|_{\hat{\mathscr{D}}^{\beta^\prime}_{W^1}}\big)|t-s|^{2\beta}\\
\leq 
&C\|\sigma_2\|_{C^3_b}\big(\|\sigma_2\|_{C^1_b}+1\big)^3
\big(\|W^1\|_\beta+1\big)^3
\big(|x_0-\tilde{x}_0|+|y_0-\tilde{y}_0|\big)|t-s|^{2\beta}
\\
&+
C\|\sigma_2\|_{C^3_b}
\big(\|\sigma_2\|_\infty+1\big)^2
\big(\|W^1\|_\beta+1\big)^3\big(\|X-\tilde{X}\|_{\alpha^\prime}+\|(Y,Y^\dagger)-(\tilde{Y},\tilde{Y}^\dagger)\|_{\hat{\mathscr{D}}^{\beta^\prime}_{W^1}}\big)|t-s|^{2\beta}.
\end{align*}
Combining the above estimates, we obtain
\begin{align}\label{eq-21}
\|\Delta^\sharp\|_{\gamma^\prime}
=&\|\sigma_2(Z)^\sharp-\sigma_2(\tilde{Z})^\sharp\|_{\gamma^\prime}\notag\\
\leq &
C\big(\|\sigma_2\|_{C^3_b}+1\big)^4
\big(\|W^1\|_\beta+1\big)^3
\big(|x_0-\tilde{x}_0|+|y_0-\tilde{y}_0|\big)\notag\\
&+
C
\big(\|\sigma_2\|_{C^3_b}+1\big)^3
\big(\|W^1\|_\beta+1\big)^3
\big(\|X-\tilde{X}\|_{\alpha^\prime}+\|(Y,Y^\dagger)-(\tilde{Y},\tilde{Y}^\dagger)\|_{\hat{\mathscr{D}}^{\beta^\prime}_{W^1}}\big)
.
\end{align}
Therefore, $
(\Delta,\Delta^\dagger)\in\mathscr{D}^{\alpha^\prime,\beta^\prime}_{W^1}$, and 
\begin{align}\label{eq-27}
\|(\Delta,\Delta^\dagger)\|_{\mathscr{D}^{\alpha^\prime,\beta^\prime}_{W^1}}
&=
\|\Delta^\dagger\|_{\beta^\prime}
+
\|\Delta^\sharp\|_{\gamma^\prime}\notag\\
&\leq 
C\big(\|\sigma_2\|_{C^3_b}+1\big)^4
\big(\|W^1\|_\beta+1\big)^3
\big(|x_0-\tilde{x}_0|+|y_0-\tilde{y}_0|\big)\notag\\
&\quad+
C
\big(\|\sigma_2\|_{C^3_b}+1\big)^3
\big(\|W^1\|_\beta+1\big)^3
\big(\|X-\tilde{X}\|_{\alpha^\prime}+\|(Y,Y^\dagger)-(\tilde{Y},\tilde{Y}^\dagger)\|_{\hat{\mathscr{D}}^{\beta^\prime}_{W^1}}\big)
.
\end{align}

Plugging \eqref{eq-70}, \eqref{eq-25}, \eqref{eq-26} and \eqref{eq-27} into \eqref{eq-24}, we obtain
\begin{align*}
\lefteqn{\|\mathcal{M}^2_{[0,\tau]}(X,Y,Y^\dagger)-\mathcal{M}^2_{[0,\tau]}(\tilde{X},\tilde{Y},\tilde{Y}^\dagger)\|_{\alpha^\prime,\beta^\prime,W^1}}\\
&\leq 
A+\|\Delta\|_{\beta^\prime}
+C\big(|\Delta^\dagger_0|+\|(\Delta,\Delta^\dagger)\|_{\mathscr{D}^{\alpha^\prime,\beta^\prime}_{W^1}}\big)
\big(\|W^1\|_\beta+\|W^2\|_{2\beta}\big)\tau^{\beta-\beta^\prime}\\
&\leq 
CS_2
(\|\sigma_2\|_{C^1_b}+1)
\big(\|W^1\|_\beta+1\big)
\big(|x_0-\tilde{x}_0|+|y_0-\tilde{y}_0|\big)
\|B\|_\alpha\tau^{\alpha-\alpha^\prime}\\
&\quad+
CS_2
\big(\|W^1\|_\beta+1\big)
\big(
\|X-\tilde{X}\|_{\alpha^\prime}
+
\|(Y,Y^\dagger)-(\tilde{Y},\tilde{Y}^\dagger)\|_{\hat{\mathscr{D}}^{\beta^\prime}_{W^1}}
\big)
\|B\|_\alpha\tau^{\alpha-\alpha^\prime}\\
&\quad+
C
\big(\|\sigma_2\|_{C^3_b}+1\big)^4
\big(\|W^1\|_\beta+\|W^2\|_{2\beta}+1\big)^4
\big(|x_0-\tilde{x}_0|+|y_0-\tilde{y}_0|\big)\tau^{\beta-\beta^\prime}\\
&\quad+
C
\big(\|\sigma_2\|_{C^3_b}+1\big)^3
\big(\|W^1\|_\beta+\|W^2\|_{2\beta}+1\big)^4
\big(\|X-\tilde{X}\|_{\alpha^\prime}+\|(Y,Y^\dagger)-(\tilde{Y},\tilde{Y}^\dagger)\|_{\hat{\mathscr{D}}^{\beta^\prime}_{W^1}}\big)
\tau^{\beta-\beta^\prime}\\
&\leq 
C
\big(S_2+1\big)\big(\|B\|_\alpha+1\big)
\big(\|\sigma_2\|_{C^3_b}+1\big)^4
\big(\|W^1\|_\beta+\|W^2\|_{2\beta}+1\big)^4 
\big(|x_0-\tilde{x}_0|+|y_0-\tilde{y}_0|\big)\tau^{\beta-\beta^\prime}\\
&\quad+
C
\big(S_2+1\big)\big(\|B\|_\alpha+1\big)
\big(\|\sigma_2\|_{C^3_b}+1\big)^3
\big(\|W^1\|_\beta+\|W^2\|_{2\beta}+1\big)^4\\
&\qquad\qquad\times
\big(\|X-\tilde{X}\|_{\alpha^\prime}+\|(Y,Y^\dagger)-(\tilde{Y},\tilde{Y}^\dagger)\|_{\hat{\mathscr{D}}^{\beta^\prime}_{W^1}}\big)
\tau^{\beta-\beta^\prime},
\end{align*}
where $S_2$ is defined in \eqref{eq-23}.

Consequently, we obtain
\begin{align*}
\lefteqn{\|\mathcal{M}^{x_0,y_0}_{[0,\tau]}(X,Y,Y^\dagger)-\mathcal{M}^{\tilde{x}_0,\tilde{y}_0}_{[0,\tau]}(\tilde{X},\tilde{Y},\tilde{Y}^\dagger)\|_{\alpha^\prime,\beta^\prime,W^1}}\\
&\leq 
CL_2\big(\|\sigma_2\|_{C^1_b}+1\big)\big(\|W^1\|_\beta+1\big)\big(|x_0-\tilde{x}_0|+|y_0-\tilde{y}_0|\big)\tau^{\beta-\beta^\prime}\notag\\
&\quad+L_2\big(\|W^1\|_\beta+1\big)\big(\|X-\tilde{X}\|_{\alpha^\prime}+\|(Y,Y^\dagger)-(\tilde{Y},\tilde{Y}^\dagger)\|_{\hat{\mathscr{D}}^{\beta^\prime}_{W^1}}\big)\tau^\beta\\
&\quad+
C
\big(S_2+1\big)\big(\|B\|_\alpha+1\big)
\big(\|\sigma_2\|_{C^3_b}+1\big)^4
\big(\|W^1\|_\beta+\|W^2\|_{2\beta}+1\big)^4 
\big(|x_0-\tilde{x}_0|+|y_0-\tilde{y}_0|\big)\tau^{\beta-\beta^\prime}\\
&\quad+
C
\big(S_2+1\big)\big(\|B\|_\alpha+1\big)
\big(\|\sigma_2\|_{C^3_b}+1\big)^3
\big(\|W^1\|_\beta+\|W^2\|_{2\beta}+1\big)^4\\
&\qquad\qquad\times
\big(\|X-\tilde{X}\|_{\alpha^\prime}+\|(Y,Y^\dagger)-(\tilde{Y},\tilde{Y}^\dagger)\|_{\hat{\mathscr{D}}^{\beta^\prime}_{W^1}}\big)
\tau^{\beta-\beta^\prime}\\
&\leq 
C (L_2+1)
\big(S_2+1\big)\big(\|B\|_\alpha+1\big)
\big(\|\sigma_2\|_{C^3_b}+1\big)^4
\big(\|W^1\|_\beta+\|W^2\|_{2\beta}+1\big)^4 
\big(|x_0-\tilde{x}_0|+|y_0-\tilde{y}_0|\big)\tau^{\beta-\beta^\prime}\\
&\quad+
C (L_2+1)
\big(S_2+1\big)\big(\|B\|_\alpha+1\big)
\big(\|\sigma_2\|_{C^3_b}+1\big)^3
\big(\|W^1\|_\beta+\|W^2\|_{2\beta}+1\big)^4\\
&\qquad\qquad\times
\big(\|X-\tilde{X}\|_{\alpha^\prime}+\|(Y,Y^\dagger)-(\tilde{Y},\tilde{Y}^\dagger)\|_{\hat{\mathscr{D}}^{\beta^\prime}_{W^1}}\big)
\tau^{\beta-\beta^\prime}.
\end{align*}
Let $\tau=\tau_1$ be small enough, we have 
\begin{align*}
\lefteqn{\big\|\mathcal{M}^{x_0,y_0}_{[0,\tau_1]}(X,Y,Y^\dagger)-\mathcal{M}^{\tilde{x}_0,\tilde{y}_0}_{[0,\tau_1]}(\tilde{X},\tilde{Y},\tilde{Y}^\dagger)\big\|_{\alpha^\prime,\beta^\prime,W^1}}\\
&\leq
\frac12|(x_0,y_0)-(\tilde{x}_0,\tilde{y}_0)|
+
\frac12
\big\|(X,Y,Y^\dagger)-(\tilde{X},\tilde{Y},\tilde{Y}^\dagger)\big\|_{\alpha^\prime,\beta^\prime,W^1}.
\end{align*}
Note that $\tau_1$ is determined by $L_2$, $\|W^1\|_\beta$, $\|W^2\|_{2\beta}$, $\|B\|_\alpha$, $S_2$ and $\|\sigma_2\|_{C^3_b}$, but is choosen independent of $(x_0,y_0)$.

If
$(x_0,y_0)=(\tilde{x}_0,\tilde{y}_0)$
in particular, this estimate implies that $\mathcal{M}^{x_0,y_0}$ is a contraction on $\mathcal{B}^{x_0,y_0}_{[0,\tau_1]}$ and has a unique fixed point in this set, which is a local solution of the system \eqref{eq-2}.
Thus, we have obtained a solution on $[0,\tau_1]$.

Next, we can choose $\tau_2=2\tau_1$ and do the same argument on the second interval $[\tau_1,\tau_2]$ with the initial condition $(x_0,y_0)$ at the time 0 being replaced by $(X_{\tau_1},Y_{\tau_1})$ at the time $\tau_1$. 
Since all the estimates above are independent of $(x_0,y_0)$, $(X_t,Y_t,Y^\dagger_t)_{t\in[\tau_1,\tau_2]}$ satisfies the same estimates as those for $(X_t,Y_t,Y^\dagger_t)_{t\in[0,\tau_1]}$. 
Concatenating them as in Example 4, we obtain a solution on $[0,2\tau_1]$.
We can continue this procedure to obtain a global solution $(X_t,Y_t,Y^\dagger_t)_{t\in[0,T]}$.

The uniqueness is a local-in-time property. Therefore, it suffices to show that any two solutions $(X,Y,Y^\dagger$) and $(\tilde{X},\tilde{Y},\tilde{Y}^\dagger)$ of system \eqref{eq-2} coincide on a sufficiently small interval near $t = 0$.
Take any $\alpha^{\prime\prime}$ and $\beta^{\prime\prime}$ such that 
$\frac13<\beta^{\prime\prime}<\beta^\prime<\frac12<\alpha^{\prime\prime}<\alpha^\prime<1$,
and we work in $(\alpha^{\prime\prime}, \beta^{\prime\prime})$-H\"older topology instead of $(\alpha^\prime, \beta^\prime)$-H\"older topology. Then, 
$$
\|(X,Y,Y^\dagger)\|_{\alpha^{\prime\prime},\beta^{\prime\prime},W^1}
\leq\|(X,Y,Y^\dagger)\|_{\alpha^\prime,\beta^\prime,W^1} 
\tau^{(\alpha^\prime-\alpha^{\prime\prime})\wedge(\beta^\prime-\beta^{\prime\prime})}.
$$
Hence, if $\tau>0$ is sufficiently small, both $(X,Y,Y^\dagger$) and $(\tilde{X},\tilde{Y},\tilde{Y}^\dagger)$ restricted to $[0,\tau]$ belong to 
$\mathcal{B}^{x_0,y_0}_{[0,\tau]}$.
Hence, $(X,Y,Y^\dagger$) and $(\tilde{X},\tilde{Y},\tilde{Y}^\dagger)$ must be identically equal on $[0,\tau]$. Since the choice of $\tau$ is uniform, the above argument can be iterated finitely many times to conclude that the two solutions coincide on the whole interval $[0,T]$. 

This completes the proof of Theorem \ref{th-1}.
\end{proof}

\begin{remark}\label{rem-1}
We make assumptions as follows: $\sigma_1(t,\cdot,\cdot)$ is of $C^2$ for any fixed $t\in[0,T]$ and $\sigma_2$ is of $C^3$. Moreover, 
for any $N>0$, there exists a positive constant $C_N$ such that 
$$
\sup_{\substack{t\in[0,T]\\|x|+|y|\leq N}}|\nabla^i\sigma_1(t,x,y)|\leq C_N,
$$
and
$$
|\sigma_1(t,x,y)-\sigma_1(s,x,y)|+|\nabla\sigma_1(t,x,y)-\nabla\sigma_1(s,x,y)|\leq C_N|t-s|^\beta, 
$$
for all $s,t\in[0,T]$ and all $(x,y)\in\mathbb{R}^{d_1+d_2}$ with $|x|+|y|\leq N$, where $C_N$ may depend on $N$.
For any $N>0$, there exists a positive constant $L_N$ such that 
$$
|b_1(t,x_1,y_1)-b_1(t,x_2,y_2)|+|b_2(x_1,y_1)-b_2(x_2,y_2)|
\leq L_N(|x_1-x_2|+|y_1-y_2|),
$$
for all $t\in[0,T]$ and all $(x_i,y_i)\in\mathbb{R}^{d_1+d_2}$ with $|x_i|+|y_i|\leq N$, $i=1,2$, where $L_N$ may depend on $N$.

Then, by standard cut-off argument, it immediately follows that under the above conditions,
the system \ref{eq-2} has a unique local solution for any given initial value. Hence, a unique solution exists up to either the explosion time or the time horizon $T$.

\end{remark}

From this point onward, we write
$\sigma_1(t,x,y)=\sigma_1(t,x),$
that is, we consider a diffusion coefficient
$\sigma_1:[0,T]\times\mathbb{R}^{d_1}
\to\mathcal{L}(\mathbb{R}^m,\mathbb{R}^{d_1})$
that is independent of the variable $y$.
For a continuous path $\varphi:~[0,T]\to \mathbb{R}^{d_2}$,
we consider the following single-time-scale YDE.
\begin{align}\label{eq-28}
dZ_t=b_1(t,Z_t,\varphi_t)dt+\sigma_1(t,Z_t)dB_t,\quad
Z_0=z_0\in\mathbb{R}^{d_1},\quad t\in[0,T], 
\end{align} 
where $B\in\mathcal{C}^\alpha([0,T],\mathbb{R}^m)$.
From this point onward, we assume that \(\sigma_1\) satisfies the following condition:
\begin{itemize}
	\item[\textbf{(A1)$^\prime$}] The function $\sigma_1(t,\cdot)\in C^2_b(\mathbb{R}^{d_1}, \mathcal{L}(\mathbb{R}^m,\mathbb{R}^{d_1}))$, for any fixed $t\in[0,T]$.
	Moreover, $\sigma_1$ is $\alpha$-H\"{o}lder continuous in time, uniformly in space. More precisely,
	\begin{align}\label{eq-56}
		S_k:=\sup_{t\in[0,T]}\Big\{\sum_{i=0}^{k}\|\nabla^i\sigma_1(t,\cdot)\|_\infty\Big\}<\infty,
	\end{align} 
	for $k=0,1,2$,
	and
	there exist constants 
	$L_1> 0$ such that, for all $x\in \mathbb{R}^{d_1}$ and $s,t\in[0,T]$,
	\begin{align*}
		|\nabla\sigma_1(t,x) - \nabla\sigma_1(s,x)| 
		+ |\sigma_1(t,x) - \sigma_1(s,x)| 
		&\le L_1 |t-s|^\beta.
	\end{align*}
	No differentiability with respect to the time variable is assumed.
\end{itemize}

\begin{lemma}\label{lem-3}
Assume \textbf{(A1)$^\prime$} and \textbf{(A3)}. Then, for every $B\in\mathcal{C}^\alpha([0,T],\mathbb{R}^m)$, $z_0\in\mathbb{R}^{d_1}$ and $\varphi$, there exists a unique global solution $Z\in\mathcal{C}^\alpha([0,T],\mathbb{R}^{d_1})$ of the YDE \eqref{eq-28}. Moreover, it satisfies the following estimate:
there exist positive constants $C$ and $l$ independent of $B$, $z_0$, $\varphi$, $b_1$ and $\sigma_1$, such that, for any $\alpha^\prime\in(1/2,\alpha)$,
\begin{align*}
\|Z\|_{\alpha^\prime}\leq C\big\{
(K+1)(\|B\|_\alpha+1)
\big\}^l,
\end{align*}
where $K:=\max\{S_2,L_1,L_2,\|b_1\|_\infty\}$.
\end{lemma}
\begin{proof}
Following the same argument as in the proof of Theorem \ref{th-1}, with considerably simpler estimates, we obtain the existence and uniqueness of a global solution to the YDE \eqref{eq-28}. It remains to establish the desired priori estimate.
For a fixed $\tau\in(0,T]$,  set
\begin{align*}
&\mathcal{M}^{z_0}_{[0,\tau]}(Z):= z_0+
\int_0^\cdot b_1(t,Z_t,\varphi_t)dt
+\int_0^\cdot\sigma_1(t,Z_t)dB_t,\\
&\mathcal{B}^{z_0}_{[0,\tau]}:=\big\{
Z\in\mathcal{C}^{\alpha^\prime}([0,T],\mathbb{R}^{d_1})~:~\|Z\|_{\alpha^\prime;[0,\tau]}\leq 1,Z_0=z_0
\big\}.
\end{align*}
By the estimates established in the proof of Theorem \ref{th-1}, there exist positive constants $C_1$ and $C_2$, independent of $B$, $z_0$, $\varphi$, $b_1$ and $\sigma_1$, such that, for any $Z\in\mathcal{B}^{z_0}_{[0,\tau]}$ and $\tilde{Z}\in\mathcal{B}^{\tilde{z}_0}_{[0,\tau]}$, the following estimates hold:
\begin{align*}
\|\mathcal{M}^{z_0}_{[0,\tau]}(Z)\|_{\alpha^\prime;[0,\tau]}
\leq C_1(K+1)(\|B\|_\alpha+1)\tau^{\alpha-\alpha^\prime}
\end{align*}
and
\begin{align*}
\|\mathcal{M}^{z_0}_{[0,\tau]}(Z)-\mathcal{M}^{\tilde{z}_0}_{[0,\tau]}(\tilde{Z})\|_{\alpha^\prime;[0,\tau]}
\leq
C_2(K+1)(|z_0-\tilde{z}_0|+\|Z-\tilde{Z}\|_{\alpha^\prime;[0,\tau]})(\|B\|_\alpha+1)\tau^{\alpha-\alpha^\prime}.
\end{align*}
Choose $\lambda>0$ such that
\begin{align*}
\lambda=\big\{
2(C_1+C_2)(K+1)(\|B\|_\alpha+1)
\big\}^{\frac{-1}{\alpha-\alpha^\prime}}.
\end{align*}
Then, for every $\tau\leq \lambda$, the above estimates imply that
\begin{align*}
\|\mathcal{M}^{z_0}_{[0,\tau]}(Z)\|_{\alpha^\prime;[0,\tau]}\leq \frac12,
\end{align*}
and 
\begin{align*}
\|\mathcal{M}^{z_0}_{[0,\tau]}(Z)-\mathcal{M}^{\tilde{z}_0}_{[0,\tau]}(\tilde{Z})\|_{\alpha^\prime;[0,\tau]}\leq \frac12\big(|z_0-\tilde{z}_0|+\|Z-\tilde{Z}\|_{\alpha^\prime;[0,\tau]}\big),
\end{align*}
that is 
$\mathcal{M}^{z_0}_{[0,\tau]}$
maps $\mathcal{B}^{z_0}_{[0,\tau]}$
into itself and is a contraction. Consequently, on every time interval of length at most $\lambda$, the solution satisfies a uniform $\alpha^\prime$-H\"older estimate.
By H\"older's inequality for finite sums, we have
\begin{align*}
\|Z\|_{\alpha^\prime;[0,T]}
\leq
\|Z\|_{\alpha^\prime;[0,\tau]}\Big\lfloor\frac T\lambda\Big\rfloor^{1-\alpha^\prime}
\leq  
C\big\{
(K+1)(\|B\|_\alpha+1)
\big\}^{\frac{1-\alpha^\prime}{\alpha-\alpha^\prime}}.
\end{align*}
This completes the proof.
\end{proof}

\begin{remark}\label{rem-2}
To establish the existence and uniqueness of a global solution of \eqref{eq-28} for any given continuous path $\varphi$ and initial value $z_0$, it suffices to assume that $\sigma_1$ satisfies Assumption \textbf{(A1)$^\prime$} and that $b_1$ satisfies
\begin{align}\label{eq-55}
\sup_{t\in[0,T],x\in\mathbb{R}^{d_1}}|b_1(t,x,\varphi_t)|
+
\sup_{t\in[0,T],x,x^\prime\in\mathbb{R}^{d_1},x\neq x^\prime}\frac{|b_1(t,x,\varphi_t)-b_1(t,x^\prime,\varphi_t)|}{|x-x^\prime|}<\infty.
\end{align}

\end{remark}

Let $\tilde{\varphi}:[0,T]\to \mathbb{R}^{d_2}$ be another continuous path.
In addition to \eqref{eq-28}, we consider the following YDE with the same $B$ and $z_0$:
\begin{align}\label{eq-39}
d\tilde{Z}_t=\tilde{b}_1(t,\tilde{Z}_t,\tilde{\varphi}_t)dt+\sigma_1(t,\tilde{Z}_t)dB_t,\quad
\tilde{Z}_0=z_0\in\mathbb{R}^{d_1},\quad t\in[0,T].
\end{align}
We assume that $\tilde{b}_1:[0,T]\times \mathbb{R}^{d_1}\times\mathbb{R}^{d_2}\to \mathbb{R}^{d_1}$ also satisfies \eqref{eq-55}, with $b_1$ and $\varphi$ replaced by $\tilde{b}_1$ and $\tilde{\varphi}$, respectively. In other words, $\tilde{b}_1$ is bounded and 
Lipschitz continuous in its first spatial variable, uniformly in time and in its second spatial variable. 
\begin{proposition}\label{pro-6}
Let $b_1$, $\tilde{b}_1$ and $\sigma_1$ satisfy the assumptions stated above. Let $B\in\mathcal{C}^\alpha([0,T],\mathbb{R}^m)$, $z_0\in\mathbb{R}^{d_1}$ and let $\varphi,\tilde\varphi:[0,T]\to\mathbb{R}^{d_2}$ be continuous paths. Let $Z$ and $\tilde{Z}$ denote the unique solutions to the YDEs \eqref{eq-28} and \eqref{eq-39} on $[0,T]$, respectively. Let $g:[0,T]\times \mathbb{R}^{d_1}\to \mathbb{R}^{d_1}$ be bounded and globally Lipschitz continuous in its spatial variable, uniformly in time. 
Define
$$M_t:=(Z_t-\tilde{Z}_t)-\int_0^t\big\{g(s,Z_s)-g(s,\tilde{Z}_s)\big\}ds-\int_0^t\big\{\sigma_1(s,Z_s)-\sigma_1(s,\tilde{Z}_s)\big\}dB_s,\quad t\in[0,T].$$
Then, $M\in\mathcal{C}^1([0,T],\mathbb{R}^{d_1})$ (i.e. Lipschitz continuous).
Moreover, for every $\alpha^\prime\in(1/2,\alpha)$, there exist constants $C>0$ and $l>0$, depending only on $\alpha$ and $\alpha^\prime$, such that 
$$\|Z-\tilde{Z}\|_{\alpha^\prime}\leq C \exp{\big[C\big\{(K^\prime+1)(\|B\|_\alpha+1)\big\}^l\big]}\|M\|_{\alpha^\prime},$$
where  $K^\prime:=\max\{S_2,L_1,L_2,\|b_1\|_\infty, L_{\tilde{b}_1},\|\tilde{b}_1\|_\infty,\|g\|_\infty,L_g\}.$
Here, $L_{\tilde{b}_1},L_g> 0$ are constants such that, for all $x_1,x_2\in \mathbb{R}^{d_1}$, $t\in[0,T]$ and continuous path $\varphi:[0,T]\to\mathbb{R}^{d_2}$,
\begin{align*}
	|\tilde{b}_1(t,x_1,\varphi_t) - \tilde{b}_1(t,x_2,\varphi_t)| 
	&\le L_{\tilde{b}_1}|x_1-x_2| ,\\
	|g(t,x_1) -g(t,x_2)| &\le L_g |x_1-x_2|.
\end{align*}
\end{proposition}
\begin{proof}
By the definition of $M$, we have $$M_t=\int_0^t\big\{b(s,Z_s,\varphi_s)-\tilde{b}_1(s,\tilde{Z}_s,\tilde{\varphi}_s)\big\}ds-\int_0^t\big\{g(s,Z_s)-g(s,\tilde{Z}_s)\big\}ds,$$
Since the integrands are bounded and continuous, $M\in\mathcal{C}^1([0,T],\mathbb{R}^{d_1})$.

Without loss of generality, we  assume that $T=1$. Mimicking the proof of Lemma \ref{lem-3}, we set 
$$\lambda^\prime:=\big\{2(C_1+C_2)(K^\prime+1)(\|B\|_\alpha+1)\big\}^{\frac{-1}{\alpha-\alpha^\prime}},$$
where $C_1$ and $C_2$ are the constants appearing in the proof of Lemma \ref{lem-3}, and
$$K^\prime:=\max\{S_2,L_1,L_2,\|b_1\|_\infty, L_{\tilde{b}_1},\|\tilde{b}_1\|_\infty,  \|g\|_\infty,L_g\}.$$ 
Set $s_j:=j\lambda^\prime$ for $0<j\leq \lfloor\frac1{\lambda^\prime}\rfloor$ and $s_N:=1$, with $N:= \lfloor\frac1{\lambda^\prime}\rfloor+1$. 
Then, by the same argument as in the proof of Lemma \ref{lem-3}, the choice of $\lambda^\prime$ ensures that, on each subinterval $[s_{j-1},s_j]$, for
$Z\in \mathcal{B}^{Z_{s_{j-1}}}_{[s_{j-1},s_j]}$ and $\tilde{Z}\in \mathcal{B}^{\tilde{Z}_{s_{j-1}}}_{[s_{j-1},s_j]}, $
we have
\begin{align*}
\|\mathcal{M}^{Z_{s_{j-1}}}_{[s_{j-1},s_j]}(Z)\|_{\alpha^\prime;[s_{j-1},s_j]}&\leq \frac12,\\
\|\mathcal{M}^{Z_{s_{j-1}}}_{[s_{j-1},s_j]}(Z)-\mathcal{M}^{\tilde{Z}_{s_{j-1}}}_{[s_{j-1},s_j]}(\tilde{Z})\|_{\alpha^\prime;[s_{j-1},s_j]}&\leq \frac12\big(|Z_{s_{j-1}}-\tilde{Z}_{s_{j-1}}|+\|Z-\tilde{Z}\|_{\alpha^\prime;[s_{j-1},s_j]}\big).
\end{align*}
Consequently, by a similar argument, for every $0\leq j\leq \lfloor\frac1{\lambda^\prime}\rfloor$, we obtain
\begin{align*}
&\Big\|\int_{s_{j-1}}^\cdot\big\{g(s,Z_s)-g(s,\tilde{Z}_s)\big\}ds
+\int_{s_{j-1}}^\cdot\big\{\sigma_1(s,Z_s)-\sigma_1(s,\tilde{Z}_s)\big\}dB_s\Big\|_{\alpha^\prime;[s_{j-1},s_j]}\\
\leq &\frac12 \big(|Z_{s_{j-1}}-\tilde{Z}_{s_{j-1}}|+\|Z-\tilde{Z}\|_{\alpha^\prime;[s_{j-1},s_j]}\big).
\end{align*}
Therefore,
\begin{align*}
\|Z-\tilde{Z}\|_{\alpha^\prime;[s_{j-1},s_j]}
\leq & \|M\|_{\alpha^\prime}
+ \frac12 \big(|Z_{s_{j-1}}-\tilde{Z}_{s_{j-1}}|+\|Z-\tilde{Z}\|_{\alpha^\prime;[s_{j-1},s_j]}\big).
\end{align*}
Rearranging the above inequality yields
\begin{align}\label{eq-40}
\|Z-\tilde{Z}\|_{\alpha^\prime;[s_{j-1},s_j]}\leq 2\|M\|_{\alpha^\prime}+|Z_{s_{j-1}}-\tilde{Z}_{s_{j-1}}|.
\end{align}
We next estimate the difference at the endpoints of the subintervals. By the triangle inequality and \eqref{eq-40}, since $\lambda^\prime\leq 1$, we obtain 
\begin{align*}
|Z_{s_j}-\tilde{Z}_{s_j}|-|Z_{s_{j-1}}-\tilde{Z}_{s_{j-1}}|
\leq &
|Z^1_{s_{j-1},s_j}-\tilde{Z}^1_{s_{j-1},s_j}|\\
\leq &
\big(2\|M\|_{\alpha^\prime}+|Z_{s_{j-1}}-\tilde{Z}_{s_{j-1}}|\big)(\lambda^\prime)^{\alpha^\prime}\\
\leq &
2\|M\|_{\alpha^\prime}
+|Z_{s_{j-1}}-\tilde{Z}_{s_{j-1}}|.
\end{align*}
Iterating this inequality and using $Z_0=\tilde{Z}_0=z_0$, we obtain
\begin{align*}
|Z_{s_j}-\tilde{Z}_{s_j}|
\leq 
2^j|z_0-z_0|+2\sum^{j-1}_{k=0}2^k\|M\|_{\alpha^\prime} 
=2(2^j-1)\|M\|_{\alpha^\prime}.
\end{align*}
Substituting this estimate into \eqref{eq-40}, we have, for every $0\leq j\leq \lfloor\frac1{\lambda^\prime}\rfloor$,
$$\|Z-\tilde{Z}\|_{\alpha^\prime;[s_{j-1},s_j]}
\leq 2\|M\|_{\alpha^\prime}+(2^N-2)\|M\|_{\alpha^\prime}=2^N\|M\|_{\alpha^\prime}.$$ 
Finally, combining the estimates on all subintervals and using H\"older's inequality for finite sums, we obtain
\begin{align*}
\|Z-\tilde{Z}\|_{\alpha^\prime,[0,1]}
\leq &
\|Z-\tilde{Z}\|_{\alpha^\prime;[s_{j-1},s_j]}\lfloor\frac1{\lambda^\prime}\rfloor^{1-\alpha^\prime}\\
\leq &
2^NN^{1-\alpha^\prime}\|M\|_{\alpha^\prime}\\
\leq &
2^{2N}\|M\|_{\alpha^\prime}\\
\leq &
C\exp{[C\{(K^\prime+1)(\|B\|_\alpha+1)\}^{\frac1{\alpha-\alpha^\prime}}]}\|M\|_{\alpha^\prime}.
\end{align*}
This completes the proof. 
\end{proof}

\subsection{Young differential equations: The case of bounded and locally Lipschitz drift vector field}\label{sec-3.2}

For a continuous path $\varphi:[0,T]\to \mathbb{R}^{d_2}$,
we consider the YDE \eqref{eq-28} with a locally Lipschitz drift vector field.
\begin{proposition}\label{pro-4}
Suppose that the diffusion coefficient $\sigma_1$ satisfies \textbf{(A1)$^\prime$} and that the drift coefficient $b_1$ satisfies the following conditions.
\begin{itemize}
\item[(1)] $b_1$ is bounded, that is,  $\|b_1\|_\infty:=\sup_{t\in[0,T], (x,y)\in\mathbb{R}^{d_1+d_2}}|b_1(t,x,y)|<\infty$.
\item[(2)] $b_1$ is locally Lipschitz continuous with respect to the first spatial variable, uniformly in  $t\in[0,T]$, in the sense that
$$
\sup_{\substack{t\in[0,T], |x_1|\vee|x_2|\leq N, \\x_1\neq x_2}}\frac{|b_1(t,x_1,\varphi_t)-b_1(t,x_2,\varphi_t)|}{|x_1-x_2|}<\infty,  
$$ 
for every $N\in\mathbb{N}$.
\end{itemize}
Then, for every $B\in\mathcal{C}^\alpha([0,T],\mathbb{R}^m)$, $z_0\in\mathbb{R}^{d_1}$ and admissible $\varphi$, the YDE \eqref{eq-28} has a unique global solution $Z\in\mathcal{C}^\alpha([0,T],\mathbb{R}^{d_1})$. Moreover, there exist constants $C,l>0$ such that 
\begin{align}\label{eq-35}
	\|Z\|_{\alpha^\prime}\leq& C\big\{\|b_1\|_\infty+(S_0+L_1)\|B\|_\alpha\big\}\big\{(S_1\|B\|_\alpha)^{\frac{1-\alpha^\prime}{\alpha-\alpha^\prime}}\vee 1+1\big\}\notag\\
	\leq& C\big\{\|b_1\|_\infty+S_0+L_1\big\}\big(\|B\|_\alpha+1\big)^l
	.
\end{align}
Here, $C$ and $l$ are independent of $B$, $z_0$, $\varphi$, $b_1$ and  $\sigma_1$. 
\end{proposition}
\begin{proof}
Without loss of generality, we assume that $T=1$.
By Remarks \ref{rem-1} and \ref{rem-2}, there exists a  unique local solution, denoted by $(Z_t)_{t\in[0,T_0]}$, for some $T_0\in(0,1]$. 
It suffices to show that $\|Z\|_{\alpha^\prime,[0,T_0]}$  is dominated by the right-hand side  of \eqref{eq-35}, with a bound independent of $T_0$. Indeed, this would imply that $Z$ does not explode at $T_0$, and hence the local solution can be extended to the whole interval $[0,1]$.

For $h\in(0,T_0]$, we introduce the notation $\|\cdot\|_{\alpha^\prime;h}:=\sup_{0<t-s\leq h}\|\cdot\|_{\alpha^\prime;[s,t]}$.
By the standard estimate for Young integrals \eqref{eq-57}, together with the condition \(\alpha^\prime+\beta^\prime>1\), we obtain
\begin{align*}
\|Z\|_{\alpha^\prime;h}
\leq &\|b_1\|_\infty \sup_{0<t-s\leq h}|t-s|^{1-\alpha^\prime}
+C_1\big(\sup_{0\leq s\leq T}|\sigma_1(s,Z_s)|+\|\sigma_1(\cdot,Z_\cdot)\|_{\beta^\prime;h}\big)\|B\|_{\alpha^\prime;h}\\
\leq &
\|b_1\|_\infty h^{1-\alpha^\prime}\\
&+C_1\Big(S_0+\sup_{0<t-s\leq h}\frac{|\sigma_1(t,Z_t)-\sigma_1(s,Z_t)|}{|t-s|^{\beta^\prime}}+\sup_{0<t-s\leq h}\frac{|\sigma_1(s,Z_t)-\sigma_1(s,Z_s)|}{|t-s|^{\beta^\prime}}\Big)\|B\|_{\alpha;h}h^{\alpha-\alpha^\prime}\\
\\
\leq&
\|b_1\|_\infty h^{1-\alpha^\prime}+C_1(S_0+L_1+S_1\|Z\|_{\alpha^\prime;h})\|B\|_{\alpha;h}h^{\alpha-\alpha^\prime}\\
\leq &
\|b_1\|_\infty h^{1-\alpha^\prime}+C_1(S_0+L_1)\|B\|_{\alpha;h}h^{\alpha-\alpha^\prime}+C_1S_1\|B\|_{\alpha;h}h^{\alpha-\alpha^\prime}\|Z\|_{\alpha^\prime;h}.
\end{align*}

We first consider the general case where $\|B\|_\alpha\neq0$ and $S_1\neq0$.
Choose $h$ sufficiently small such that $C_1S_1\|B\|_\alpha h^{\alpha-\alpha^\prime}\leq \frac12$.
Equivalently, $h\leq C_2(S_1\|B\|_\alpha)^{\frac{-1}{\alpha-\alpha^\prime}}$ with a suitable constant $C_2>0$.
Then, for every $h\leq \Lambda$, where $\Lambda:=C_2(S_1\|B\|_\alpha)^{\frac{-1}{\alpha-\alpha^\prime}}\wedge T_0\leq 1$,
we have 
\begin{align}\label{eq-36}
	\|Z\|_{\alpha^\prime;h}\leq 2\|b_1\|_\infty h^{1-\alpha^\prime}+2C_1(S_0+L_1)\|B\|_\alpha h^{\alpha-\alpha^\prime}.
\end{align}
Set $N:=\lfloor\frac{T_0}{\Lambda}\rfloor+1
\leq C_3
(S_1\|B\|_\alpha)^{\frac{1}{\alpha-\alpha^\prime}}\vee 1+1$
and $s_j:=j\Lambda$, $0\le j\leq \lfloor\frac{T_0}{\Lambda}\rfloor$, $s_N:=T_0$.
Since \eqref{eq-36} is valid on each subinterval $[s_{j-1},s_j]$, H\"older's inequality for finite sums yields 
\begin{align*}
	\|Z\|_{\alpha^\prime;[0,T_0]}
	\leq &
	\|Z\|_{\alpha^\prime;\Lambda}N^{1-\alpha^\prime}\\
	\leq &
	\big\{2\|b_1\|_\infty\Lambda^{1-\alpha^\prime}+2C_1(S_0+L_1)\|B\|_\alpha\Lambda^{\alpha-\alpha^\prime}
	\big\}N^{1-\alpha^\prime}
	\\
	\leq &
	2\big\{\|b_1\|_\infty+ C_1(S_0+L_1)\|B\|_\alpha
	\big\}\big\{C_3
	(S_1\|B\|_\alpha)^{\frac{1}{\alpha-\alpha^\prime}}\vee 1+1\big\}^{1-\alpha^\prime}
	\\
	\leq & 
	C_4\big\{\|b_1\|_\infty+(S_0+L_1)\|B\|_\alpha\big\}\big\{(S_1\|B\|_\alpha)^{\frac{1-\alpha^\prime}{\alpha-\alpha^\prime}}\vee 1+1\big\}
	.
\end{align*}
The right-hand side is independent of $T_0$. Hence, $Z$ does not explode at $T_0$, and the local solution extends globally to $[0,1]$.

We next consider the case $\|B\|_\alpha=0$. In this case,
$$
\|Z\|_{\alpha^\prime;h}\leq \|b_1\|_\infty h^{1-\alpha^\prime}.
$$
Taking $h=h_0=\frac{T_0}{2}<1$, we have $N=\frac{T_0}{h_0}=2$.
Hence, by H\"older's inequality for finite sums, 
$$
\|Z\|_{\alpha^\prime;[0,T_0]}\leq \|Z\|_{\alpha^\prime;h_0}N^{1-\alpha^\prime} \leq C\|b_1\|_\infty<\infty.
$$
Finally, suppose that $S_1=0$. Then, by \textbf{(A1)$^\prime$}, we also have $S_0=0$, and therefore
$$
\|Z\|_{\alpha^\prime;h}\leq \|b_1\|_\infty h^{1-\alpha^\prime}+C_1L_1\|B\|_{\alpha;h} h^{\alpha-\alpha^\prime}.
$$
Again, taking $h=h_0=\frac{T_0}{2}<1$, we have
\begin{align*}
\|Z\|_{\alpha^\prime;[0,T_0]}\leq \|Z\|_{\alpha^\prime;h_0}N^{1-\alpha^\prime}
\leq 
\|b_1\|_\infty +C_5L_1\|B\|_\alpha 
.
\end{align*}
This proves the proposition.
\end{proof}

We next extend Proposition \ref{pro-6} to the case in which the drift vector fields are bounded and locally Lipschitz continuous.

Let $M$ be defined as in Proposition \ref{pro-6}, with
$g:[0,T]\times\mathbb{R}^{d_1}\to \mathbb{R}^{d_1}$
assumed to be continuous.
In contrast to Proposition \ref{pro-6}, we do not assume that $b_1$, $\tilde{b}_1$ and $g$ are globally Lipschitz continuous.
For each $N\in\mathbb{N}$, define
\begin{align*}
L^N_{b_1}:=\sup_{\substack{t\in[0,T],|x|\vee|x^\prime|\leq N,\\
		x\neq x^\prime}}\frac{|b_1(t,x,\varphi_t)-b_1(t,x^\prime,\varphi_t)|}{|x-x^\prime|}.
\end{align*}
Define $L^N_{\tilde{b}_1}$ and $L^N_g$ analogously. We assume that $L^N_{b_1}$, $L^N_{\tilde{b}_1}$ and $L^N_g$ are finite.
\begin{proposition}\label{pro-10}
Let $B\in\mathcal{C}^\alpha([0,T],\mathbb{R}^m)$, $z_0\in\mathbb{R}^{d_1}$ and let $\varphi,\tilde\varphi:[0,T]\to\mathbb{R}^{d_2}$ be continuous paths. Let $Z$ and $\tilde{Z}$ denote the unique solutions to the YDEs \eqref{eq-28} and \eqref{eq-39} on $[0,T]$, respectively.
Assume that $\sigma_1$ satisfies  \textbf{(A1)$^\prime$}, 
and that $b_1$, $\tilde{b}_1$, $g$ are all bounded and satisfy 
$$L^N_{b_1}+L^N_{\tilde{b}_1}+L^N_g=O(N^{d})\quad \text{as }N\to \infty,$$ 
for some $d>0$. Here, $O$ stands for Landau's large $O$.
Then, $M\in\mathcal{C}^1([0,T],\mathbb{R}^{d_1})$ and there exist constants $C>0$ and $l>0$ such that, for every $\alpha^\prime\in(1/2,\alpha)$, 
$$\|Z-\tilde{Z}\|_{\alpha^\prime}\leq 
C\exp{\big\{C(\|B\|^l_\alpha+1)\big\}}\|M\|_{\alpha^\prime}.$$
Here, $C$ and $l$ are independent of $B$, $\varphi$, $\tilde{\varphi}$, $M$.
\end{proposition}
\begin{proof}
In this proof, $C_i$'s are positive constants independent of $B$, $\varphi$, $\tilde{\varphi}$, $M$.
By Proposition \ref{pro-4}, the equations \eqref{eq-28} and \eqref{eq-39} admit unique global solutions. Moreover, there exist constants $C_1\geq 1$ and $v>0$ such that
$$\|Z\|_\infty\vee\|\tilde{Z}\|_\infty\leq |z_0|+\big(\|Z\|_{\alpha^\prime}\vee\|\tilde{Z}\|_{\alpha^\prime}\big)
T^{\alpha^\prime}\leq C_1(\|B\|_\alpha^v+1)=:R_B\geq 1.$$
Choose a Lipschitz continuous function $\chi:\mathbb{R}^{d_1}\to [0,1]$ such that 
\begin{itemize}
\item [1.]  $\chi(x)=1$ if $|x|\leq R_B$;
\item [2.]  $\chi(x)=0$ if $|x|\geq 2R_B$;
\item [3.] the Lipschitz constant of $\chi$ is at most 1.
\end{itemize}
We claim that the truncated drift $\chi b_1$ is globally Lipschitz continuous in the $x$-variable. Indeed, for $x,x^\prime\in\mathbb{R}^{d_1}$,
\begin{align*}
|\chi(x)b_1(t,x,\varphi_t)-\chi(x^\prime)b_1(t,x^\prime,\varphi_t)|
\leq &
|\chi(x)-\chi(x^\prime)||b_1(t,x,\varphi_t)|+|\chi(x^\prime)||b_1(t,x,\varphi_t)-b_1(t,x^\prime,\varphi_t)|\\
\leq &
C_2\{\|b_1\|_\infty+(R_B)^{d}\}|x-x^\prime|\\
\leq &
C_3 \{(R_B)^{d}+1\}|x-x^\prime|.
\end{align*}
By the same argument, $\chi \tilde{b}_1$ and $\chi g$ also satisfy the same property.
Since $\|Z\|_\infty\vee\|\tilde{Z}\|_\infty
\leq R_B,$
noting that $Z$ (resp. $\tilde{Z}$) solves \eqref{eq-28} (resp. \eqref{eq-39}) whose drift vector field is replaced by $\chi b_1$ (resp. $\chi \tilde{b}_1$).

Hence, applying Proposition \ref{pro-6} with
$$
K^\prime\leq C_3\{R_B^d+1\}
\leq C_4(\|B\|_\alpha^{vd}+1),
$$
we obtain
$$
\|Z-\tilde{Z}\|_{\alpha^\prime}
\leq C\exp\big\{C(\|B\|_\alpha^{vd+1}+1)\big\}.
$$
The desired estimate follows by setting $l:=vd+1$.
\end{proof}

Next, we consider the following single-time-scale YDE.
\begin{align}\label{eq-58}
	dZ_t=b(t,Z_t)dt+\sigma_1(t,Z_t)dB_t,\quad
	Z_0=z_0\in\mathbb{R}^{d_1},\quad t\in[0,T], 
\end{align} 
where $B\in\mathcal{C}^\alpha([0,T],\mathbb{R}^m)$. This equation is a special case of \eqref{eq-28}. By Proposition \ref{pro-4}, it admits a unique solution under the assumptions of Proposition \ref{pro-11} below.
\begin{proposition}\label{pro-11}
Suppose that the diffusion coefficient $\sigma_1$ satisfies \textbf{(A1)$^\prime$} and that the drift coefficient $b$ satisfies the following conditions:
\begin{itemize}
	\item[(1)] $b$ is bounded, that is,  $\|b\|_\infty:=\sup_{t\in[0,T], x\in\mathbb{R}^{d_1}}|b(t,x)|<\infty$.
	\item[(2)] $b$ is locally Lipschitz continuous with respect to the spatial variable, uniformly in time, in the sense that
	$$
	\sup_{\substack{t\in[0,T], |x_1|\vee|x_2|\leq N, \\x_1\neq x_2}}\frac{|b(t,x_1)-b(t,x_2)|}{|x_1-x_2|}<\infty,  
	$$ 
	for every $N\in\mathbb{N}$.
\end{itemize}
	Let $Z\in\mathcal{C}^\alpha([0,T],\mathbb{R}^{d_1})$ be the unique solution to YDE
	\eqref{eq-58}. Similarly, let $\tilde{Z}$ denote the unique solution to the corresponding YDE driven by
	$\tilde{B}\in\mathcal{C}^{\alpha}([0,T],\mathbb{R}^m)$ and with initial condition $\tilde{z}_0\in\mathbb{R}^{d_1}$.
	Assume that
	\[
	\|B\|_\alpha\vee\|\tilde{B}\|_\alpha\leq M<\infty,
	\]
	and that $z_0$ and $\tilde{z}_0$ belong to a compact subset of $\mathbb{R}^{d_1}$.
	Then, for every $\alpha^\prime\in(1/2,\alpha)$, there exists a constant $C>0$, depending on $M$ and the coefficients but independent of $z_0,\tilde{z}_0,B,\tilde{B}$, such that
	\[
	\|Z-\tilde{Z}\|_{\alpha^\prime}
	\leq
	C\big(
	|z_0-\tilde{z}_0|
	+
	\|B-\tilde{B}\|_\alpha
	\big).
	\]
\end{proposition}

\begin{proof}
	Throughout this proof, subscripts on constants indicate their dependence on the corresponding parameters, and the values of the constants may change from line to line. 
	
	By Proposition \ref{pro-4} and the assumption
	$\|B\|_\alpha\vee\|\tilde{B}\|_\alpha\leq M$, there exists a constant
	$C_M>0$ such that
	\[
	\|Z\|_{\alpha^\prime}
	+
	\|\tilde{Z}\|_{\alpha^\prime}
	\leq C_M.
	\]
Since the initial values belong to a compact set,
	we also have
	\[
	\|Z\|_{\infty}
	\leq |z_0|+C_M T^{\alpha^\prime}\leq \tilde{C}_M,
	\qquad
	\|\tilde{Z}\|_{\infty}
	\leq |\tilde z_0|+C_M T^{\alpha^\prime}\leq \tilde{C}_M.
	\]
Moreover, for every $M>0$, the coefficient $b$ is Lipschitz continuous in the spatial variable on bounded subsets of $\mathbb{R}^{d_1}$, for any fixed $t\in[0,T]$. More precisely, there exists a constant $L_M>0$, depending only on $M$, such that
$$
|b(t,x)-b(t,y)|\leq L_M|x-y|,
\qquad
t\in[0,T],\quad x,y\in B_M,
$$
where $B_M:=\{x\in\mathbb{R}^{d_1}:|x|\leq \tilde{C}_M\}$.
	
	For $t\in[0,T]$, we have
	\begin{align*}
		Z_t-\tilde{Z}_t
		={}&
		z_0-\tilde{z}_0
		+\int_0^t
		\big\{
		b(s,Z_s)
		-b(s,\tilde{Z}_s)
		\big\}\,ds
		+
		\int_0^t
		\big\{
		\sigma_1(s,Z_s)
		-\sigma_1(s,\tilde{Z}_s)
		\big\}\,dB_s\\
		&+
		\int_0^t
		\sigma_1(s,\tilde{Z}_s)\,d(B_s-\tilde{B}_s).
	\end{align*}
	Let $\tau\in(0,1)$ be small. Then, we have
	\begin{align}\label{eq-59}
\Big\|
\int_0^\cdot
\big\{
b(s,Z_s)
-b(s,\tilde Z_s)
\big\}\,ds
\Big\|_{\alpha^\prime;[0,\tau]}
		\leq&
		L_M
		\|Z-\tilde Z\|_{\infty;[0,\tau]}
		\tau^{1-\alpha^\prime}
		\notag\\
		\leq&
		L_M
		\big(|z_0-\tilde{z}_0|+\|Z-\tilde{Z}\|_{\alpha^\prime;[0,\tau]}\big)\tau^{1-\alpha^\prime}.
	\end{align}
Let $0<\beta^\prime<\beta<\frac12<\alpha^\prime<\alpha<1$, and $\alpha^\prime+\beta^\prime>1$.
	Next, by the standard estimate for Young integrals \eqref{eq-57} and Lemma \ref{lem-2}, there exists a constant $C_M>0$ such that
	\begin{align}\label{eq-60}
		&\Big\|
		\int_0^\cdot
		\big\{
		\sigma_1(s,Z_s)
		-\sigma_1(s,\tilde Z_s)
		\big\}\,dB_s
		\Big\|_{\alpha^\prime;[0,\tau]}\notag\\
		&\leq
C_M		\big(|\sigma_1(0,z_0)-\sigma_1(0,\tilde{z}_0)|+\|\sigma_1(\cdot,Z_\cdot)-\sigma_1(\cdot,\tilde{Z}_\cdot)\|_{\beta^\prime;[0,\tau]}\big)\|B\|_{\alpha^\prime}\notag\\
&\leq C_M		\big(|z_0-\tilde{z}_0|+\|Z-\tilde{Z}\|_{\alpha^\prime;[0,\tau]}\big)\|B\|_\alpha\tau^{\alpha-\alpha^\prime}
.
	\end{align}
	For the last integral, the standard Young estimate \eqref{eq-57} gives
	\begin{align}\label{eq-61}
		\Big\|
		\int_0^\cdot
		\sigma_1(s,\tilde Z_s)
		\,d(B_s-\tilde B_s)
		\Big\|_{\alpha^\prime;[0,\tau]}
		\leq&
		C_M
		\big(|\sigma_1(0,\tilde z_0)|+\|\sigma_1(\cdot,\tilde Z_\cdot)\|_{\beta^\prime;[0,\tau]}\big)
		\|B-\tilde B\|_\alpha
		\tau^{\alpha-\alpha^\prime}\notag\\
		\leq &
		C_M(1+\|\tilde{Z}\|_{\alpha^\prime;[0,\tau]})	\|B-\tilde B\|_\alpha
		\tau^{\alpha-\alpha^\prime}\notag\\
		\leq &
		C_M	\|B-\tilde B\|_\alpha.
	\end{align}
Combining \eqref{eq-59}--\eqref{eq-61}, we obtain
	\begin{align*}
		\|Z-\tilde Z\|_{\alpha^\prime;[0,\tau]}
		\leq{}&
		\bar{C}_M
		\big(|z_0-\tilde z_0|+\|Z-\tilde Z\|_{\alpha^\prime;[0,\tau]}\big)
		\tau^{\alpha-\alpha^\prime}
		+
		\bar{C}_M\|B-\tilde B\|_\alpha.
	\end{align*}
Taking $\Lambda:= (2\bar{C}_M)^{\frac{-1}{\alpha-\alpha^\prime}}\wedge 1$. 
Then, for every $\tau\leq\Lambda$, satisfying $1-\bar{C}_M\tau^{\alpha-\alpha^\prime}>0$, we have 
	\begin{align*}
		\|Z-\tilde Z\|_{\alpha^\prime;[0,\tau]}
		\leq
		C_M
		\big(
		|z_0-\tilde z_0|
		+
		\|B-\tilde B\|_\alpha
		\big).
	\end{align*}
	
	We now extend the above estimate to the whole interval $[0,T]$.
	Set
$
	N:=\left\lfloor\frac{T}{\Lambda}\right\rfloor+1,
$
	and let
$
	s_j:=j\Lambda$,
	$0\leq j\leq
	\left\lfloor\frac{T}{\Lambda}\right\rfloor,
$
	with $s_N:=T$.
	
	For each $j=1,\ldots,N$, we apply the preceding argument on the
	interval $[s_{j-1},s_j]$. In this case, the initial difference is
	$Z_{s_{j-1}}-\tilde Z_{s_{j-1}}$. Thus,
	\begin{align}\label{eq-62}
		\|Z-\tilde Z\|_{\alpha^\prime;[s_{j-1},s_j]}
		\leq
		C_M
		\Big(
		|Z_{s_{j-1}}-\tilde Z_{s_{j-1}}|
		+
		\|B-\tilde B\|_\alpha
		\Big).
	\end{align}
	Moreover,
	\[
	|Z_{s_j}-\tilde Z_{s_j}|
	\leq
	|Z_{s_{j-1}}-\tilde Z_{s_{j-1}}|
	+
	\Lambda^{\alpha^\prime}
	\|Z-\tilde Z\|_{\alpha^\prime;[s_{j-1},s_j]}.
	\]
	It follows from \eqref{eq-62} that
	\[
	|Z_{s_j}-\tilde Z_{s_j}|
	\leq
	\big(
	1+C_M\Lambda^{\alpha^\prime}
	\big)
	|Z_{s_{j-1}}-\tilde Z_{s_{j-1}}|
	+
	C_M\Lambda^{\alpha^\prime}
	\|B-\tilde B\|_\alpha.
	\]
	Iterating this inequality and recalling that $N$ and $\Lambda$ depend
	only on $M$, $T$ and the coefficients, we obtain
	\[
	\max_{0\leq j\leq N}
	|Z_{s_j}-\tilde Z_{s_j}|
	\leq
	C_M
	\big(
	|z_0-\tilde z_0|
	+
	\|B-\tilde B\|_\alpha
	\big).
	\]
	Consequently, by \eqref{eq-62},
	\[
	\max_{1\leq j\leq N}
	\|Z-\tilde Z\|_{\alpha^\prime;[s_{j-1},s_j]}
	\leq
	C_M
	\big(
	|z_0-\tilde z_0|
	+
	\|B-\tilde B\|_\alpha
	\big).
	\]
	Therefore, by the finite-sum H\"older inequality, we obtain
	\[
	\|Z-\tilde Z\|_{\alpha^\prime}
	\leq
	C_M
	\big(
	|z_0-\tilde z_0|
	+
	\|B-\tilde B\|_\alpha
	\big).
	\]
	This completes the proof.
\end{proof}

\section{Averaging Principle}\label{sec-4}
Throughout this section, let $\alpha\in(1/2,1)$ and $\beta\in(1/3,1/2)$ satisfy $\alpha+\beta>1$.
We consider the following slow--fast Young-rough system:
\begin{align}\label{eq-29}
\begin{cases}
dX^\delta_t=b_1(t,X^\delta_t,Y^\delta_t)dt+\sigma_1(t,X^\delta_t)dB_t,\\
dY^\delta_t=\frac1\delta b_2(X^\delta_t,Y^\delta_t)dt+
\frac{1}{\sqrt{\delta}}\sigma_2(X^\delta_t,Y^\delta_t)d\mathbb{W}_t,\\
(X^\delta_0,Y^\delta_0)=(x_0,y_0)\in\mathbb{R}^{d_1+d_2},~ t\in[0,T],~ \delta\in(0,1].
\end{cases}
\end{align}
Here, $(x_0,y_0)$ is an arbitrary deterministic initial condition, $B=(B_t)_{t\in[0,T]}$ is an $\alpha$-H\"older continuous stochastic process, and $\mathbb{W}=(W^1,W^2)$ denotes the canonical It\^o-type Brownian rough path, associated with an $r$-dimensional
standard BM $w=(w_t)_{t\geq0}$, independent of $B$, whose precise definition will be given below. We impose the following regularity assumptions on the coefficients.
\begin{itemize}
\item For each fixed $t\in[0,T]$, $\sigma_1(t,\cdot)\in C^2(\mathbb{R}^{d_1},\mathcal{L}(\mathbb{R}^m,\mathbb{R}^{d_1}))$ and $\sigma_2\in C^3(\mathbb{R}^{d_1}\times\mathbb{R}^{d_2},\mathcal{L}(\mathbb{R}^r,\mathbb{R}^{d_2}))$. 
Moreover, 
for every $N>0$, there exists a constant $C_N>0$ such that
\begin{align*}
	\sup_{|x|\leq N}|\nabla^i\sigma_1(t,x)|&\leq C_N,\\
	|\sigma_1(t,x)-\sigma_1(s,x)|+|\nabla\sigma_1(t,x)-\nabla\sigma_1(s,x)|&\leq C_N|t-s|^\beta, 
\end{align*}
for all $i=0,1,2$, $s,t\in[0,T]$ and $x\in\mathbb{R}^{d_1}$ with $|x|\leq N$.

\item For each fixed $t\in[0,T]$, the functions $b_1(t,\cdot,\cdot)\in\mathcal{C}(\mathbb{R}^{d_1}\times\mathbb{R}^{d_2},\mathbb{R}^{d_1})$ and  $b_2\in\mathcal{C}(\mathbb{R}^{d_1}\times\mathbb{R}^{d_2},\mathbb{R}^{d_2})$ are locally Lipschitz continuous. 
More precisely, for every $N>0$, there exists a constant $L_N>0$ such that
$$
|b_1(t,x_1,y_1)-b_1(t,x_2,y_2)|+|b_2(x_1,y_1)-b_2(x_2,y_2)|
\leq L_N(|x_1-x_2|+|y_1-y_2|),
$$
for all $t\in[0,T] $ and $(x_i,y_i)\in\mathbb{R}^{d_1}\times\mathbb{R}^{d_2}$ with $|x_i|+|y_i|\leq N$, $i=1,2$.
\end{itemize}

By Remark \ref{rem-1}, the above assumptions ensure that the slow--fast system \eqref{eq-29} admits a unique local solution.

To establish the AP, we impose the following additional assumptions.
\paragraph{\textbf{(H1)}} 
The function $\sigma_1(t,\cdot)\in C^2_b(\mathbb{R}^{d_1}, \mathcal{L}(\mathbb{R}^m,\mathbb{R}^{d_1}))$  in space, for any fixed $t\in[0,T]$, 
and is $\beta$-H\"{o}lder continuous in time. More precisely,
\begin{align*}
	S_k:=\sup_{t\in[0,T]}\Big\{\sum_{i=0}^{k}\|\nabla^i\sigma_1(t,\cdot)\|_\infty\Big\}<\infty,\quad k=0,1,2.
\end{align*} 
Moreover, there exist constants 
$L_1> 0$ such that for all $(x,y)\in \mathbb{R}^{d_1}\times\mathbb{R}^{d_2}$ and $s,t\in[0,T]$,
\begin{align*}
	|\nabla\sigma_1(t,x) - \nabla\sigma_1(s,x)| 
	+ |\sigma_1(t,x) - \sigma_1(s,x)| 
	&\le L_1 |t-s|^{\beta}.
\end{align*}
No differentiability with respect to the time variable is assumed.

\paragraph{\textbf{(H2)}} 
The drift coefficient $b_1$ is bounded and 
Lipschitz continuous in the spatial variables, uniformly in time. 
More precisely,
$$
\|b_1\|_\infty
:=\sup_{t\in[0,T],\,(x,y)\in\mathbb{R}^{d_1+d_2}}
|b_1(t,x,y)|<\infty,
$$
and there exists a constant $L_2>0$ such that
\begin{align*}
	|b_1(t,x_1,y_1)-b_1(t,x_2,y_2)|
	\leq L_2\big(|x_1-x_2|+|y_1-y_2|\big),
\end{align*}
for all $t\in[0,T]$ and
$(x_i,y_i)\in\mathbb{R}^{d_1}\times\mathbb{R}^{d_2}$,
$i=1,2$.

\paragraph{\textbf{(H3)}} 
The function $\sigma_2$ is globally Lipschitz continuous.\\

We set the following conditions for $2\leq q<\infty$. If $q^\prime>q$, then \textbf{(H4)$_{q^\prime}$} obviously implies \textbf{(H4)$_q$} without changing $\gamma_1$, $C$ and $\eta_1$. 
\paragraph{\textbf{(H4)}$_q$} 
There exist constants $\gamma_1>0$, $C>0$ and $\eta_1\geq 0$ such that, for all $(x,y)\in\mathbb{R}^{d_1+d_2}$,
$$
2\langle y,b_2(x,y)\rangle+(q-1)|\sigma_2(x,y)|^2\leq -\gamma_1|y|^2+C(|x|^{\eta_1}+1).
$$

\paragraph{\textbf{(H5)}} 
There exist constants $\eta_2\geq 0$ and $C>0$ such that, for all $(x,y)\in\mathbb{R}^{d_1+d_2}$,
$$
|b_2(x,y)|\leq C(|x|^{\eta_2}+|y|^{\eta_2}+1).
$$

\paragraph{\textbf{(H6)}} 
There exists a constant $\gamma_2>0$ such that, for all $x\in\mathbb{R}^{d_1}$ and $y_1,y_2\in\mathbb{R}^{d_2}$,
$$
2\langle b_2(x,y_1)-b_2(x,y_2),y_1-y_2\rangle +|\sigma_2(x,y_1)-\sigma_2(x,y_2)|^2\leq -\gamma_2|y_1-y_2|^2
.$$

\vspace{2mm}
We set the following condition for $r\geq 0$.
\paragraph{\textbf{(H7)}$_{r}$} 
There exist constants $\eta_3\geq 0$ and $C>0$ such that, for all $x_1,x_2\in\mathbb{R}^{d_1}$ and $y\in\mathbb{R}^{d_2}$,
$$
|b_2(x_1,y)-b_2(x_2,y)|\leq C|x_1-x_2|(1+|x_1|^{\eta_3}+|x_2|^{\eta_3}+|y|^{r}).
$$

\vspace{2mm}
Finally, we impose the following assumption on the driving noise $B$.
\paragraph{\textbf{(HN)}} 
For every $p\in[1,\infty)$, $\mathbb{E}\big[\|B\|^p_\alpha\big]<\infty$, where $B$ is a $\mathcal{C}^\alpha([0,T],\mathbb{R}^m)$-valued random variable.

\vspace{2mm}
Let $(\Omega,\mathscr{F},\mathbb{P})$ be a probability space. 
Let $w=(w_t)_{t\in[0,T]}$ be a standard $r$-dimensional BM and $B=(B_t)_{t\in[0,T]}$ be a $\mathcal{C}^\alpha([0,T],\mathbb{R}^m)$-valued random variable defined on $(\Omega,\mathscr{F},\mathbb{P})$.
Assume that $B$ and $w$ are independent.
Regarding the integrability of $B$, we impose Assumption \textbf{(HN)} throughout this section. Let $\{\mathscr{F}_t\}_{t\in[0,T]}$ be a filtration satisfying the usual conditions as well as the following two conditions:
(i) $w$ is an $\{\mathscr{F}_t\}$-BM and (ii) $B$ is  $\{\mathscr{F}_t\}$-adapted.

\begin{definition}\label{de-1}
	Let $\mathbb{W}=(W^1,W^2)$ denote the canonical It\^o-type Brownian rough path, that is 
	$$
	W^1_{s,t}=w_t-w_s,\quad
	W^2_{s,t}=\int_s^t(w_u-w_s)\otimes d^I w_u,
	$$
where $d^I w$ stands for the standard It\^o integration. It is well-known that $\mathbb{W}\in\Omega^\beta([0,T],\mathbb{R}^r)$ for every $\beta\in(1/3,1/2)$. 
\end{definition}

Throughout this section, we take $\alpha^\prime$ and $\beta^\prime$ sufficiently close to $\alpha$ and $\beta$, respectively, so that $\frac13<\beta^\prime<\beta\leq \frac12<\alpha^\prime<\alpha<1$
and $\alpha^\prime+\beta^\prime>1$.

We extend the filtration $\{\mathscr{F}_t\}$ to $[0,\infty)$ by setting $\mathscr{F}_t=\mathscr{F}_{t\wedge T}$ for $t\geq 0$. Let $\hat{\mathbb{R}}^{d_1+d_2}:=\mathbb{R}^{d_1+d_2}\cup\{\infty\}$ denote the one-point compactification of $\mathbb{R}^{d_1+d_2}$. If a global solution $(Z^\delta_t)_{t\in[0,T]}:=(X^\delta_t,Y^\delta_t)_{t\in[0,T]}$ exists, then we extend it to $[0,\infty)$ by setting $Z^\delta_t=Z^\delta_T$ for $t\geq T$. Otherwise, let $(Z^\delta_t)_{t\in[0,u^\delta)}$, where $0<u^\delta\leq T$, be a maximal local solution with lifetime $u^\delta$, and set $Z^\delta_t=\infty$ for $t\geq u^\delta$.

For each $N\in\mathbb{N}$, define
$
\tau^\delta_N:=\inf\{t\geq 0:|Z^\delta_t|\geq N\},
$
and let $\tau^\delta_\infty:=\lim_{N\to\infty}\tau^\delta_N$, where we adopt the conventions $\inf\emptyset:=\infty$. Since $Z^\delta$ is adapted and has continuous sample paths up to its explosion time, both $\tau^\delta_N$ and $\tau^\delta_\infty$ are $\{\mathscr{F}_t\}$-stopping times.
Then, the following conditions are equivalent:
\begin{itemize}
\item 
A global solution $(Z^\delta_t)_{t\in[0,T]}$ of Young-rough system \eqref{eq-29} exists.
\item 
$(Z^\delta_t)_{t\in[0,T]}$ defined as above is bounded in  $\mathbb{R}^{d_1+d_2}$.
\item 
$\tau^\delta_N=\infty$ for some $N$.
\item 
$\tau^\delta_\infty>T$.
\end{itemize} 

It should be noted that while a solution of \eqref{eq-29} moves in a bounded set, its trajectory is uniformly continuous in time, since its H\"older norm is bounded by Theorem \ref{th-1} and Remark \ref{rem-1}. consequently, if $(Z^\delta_t)_{t\in[0,s)}$ is bounded for some $0 < s \leq T$, then the solution admits a continuous extension to $t=s$, and the extended process solves \eqref{eq-29} on $[0,s]$.

Conversely, if no global solution exists, then we have $u^\delta=\tau^\delta_\infty\in[0,T]$ and $\limsup_{t\uparrow\tau^\delta_\infty}|Z^\delta_t|=\infty$. Moreover, $\lim_{t\uparrow\tau^\delta_\infty}|Z^\delta_t|=\infty$ because of the H\"older estimates and uniformly continuity mentioned above. Therefore, $(Z^\delta_t)_{t\geq 0}$ has continuous sample paths in $\hat{\mathbb{R}}^{d_1+d_2}$.

The following proposition follows directly from \cite[Proposition 4.7]{inahama2025averaging}.
\begin{proposition}\label{pro-2}
Under Assumption \textbf{(HN)}, for every $\delta\in(0,1]$, $Y^\delta$ satisfies the following It\^o SDE up to the explosion time of $Z^\delta=(X^\delta,Y^\delta)$:
\begin{align*}
dY^\delta_t=\frac1\delta b_2(X^\delta_t,Y^\delta_t)dt+\frac{1}{\sqrt\delta}\sigma_2(X^\delta_t,Y^\delta_t)d^Iw_t,
\quad
Y^\delta_0=y_0\in\mathbb{R}^{d_2},\quad t\in[0,\tau^\delta_\infty).
\end{align*}
\end{proposition}

\subsection{Preliminary estimates}
In this subsection, we establish some preliminary estimates that are essential for proving the global well-posedness of both the original slow--fast system and the averaged equation.

We begin by establishing uniform estimates for the solutions.
\begin{proposition}\label{pro-3}
Assume \textbf{(HN)}, \textbf{(H1)}--\textbf{(H3)},  \textbf{(H4)$_q$} and \textbf{(H5)} for some $2\leq q<\infty$. Then, the solution $Z^\delta=(X^\delta,Y^\delta)$ of \eqref{eq-29} does not explode on $[0,T]$ almost surely. Moreover, the following uniform estimates hold:
\begin{align}
&\sup_{0<\delta\leq 1}\mathbb{E}\big[\|X^\delta\|^p_{\alpha^\prime;[0,T]}\big]<\infty,\quad1\leq p<\infty,\label{eq-31}\\
&\sup_{0<\delta\leq 1}\sup_{t\in[0,T]}\mathbb{E}\big[|Y^\delta_t|^q\big]<\infty.\label{eq-32}
\end{align}
\end{proposition}
\begin{proof}
For each $N\in\mathbb{N}$, define
$\lambda^\delta_N
:=\inf\big\{t\geq 0:|Y^\delta_t|\geq N\big\},$
and let $\lambda^\delta_\infty:=\lim_{N\to\infty}\lambda^\delta_N$.
Clearly, $\lambda^\delta_\infty\geq\tau^\delta_\infty$ almost surely.
We first show that $\lambda^\delta_\infty= \tau^\delta_\infty$ almost surely. In the other words, $X^\delta$ does not explode before $Y^\delta$ does.
Suppose, to the contrary, that $\lambda^\delta_\infty> \tau^\delta_\infty$ for some $(B,\mathbb{W})$.
This implies that $\tau^\delta_\infty\leq T$ and that $Y^\delta$ remains bounded on $[0,\tau^\delta_\infty)$.
By Lemma \ref{lem-3}, $X^\delta$ also remains bounded on $[0,\tau^\delta_\infty)$, and hence so does $Z^\delta$. This contradicts the definition of $\tau^\delta_\infty$.

Lemma \ref{lem-3} also implies that if $\tau^\delta_\infty\leq T$, that is, if $Z^\delta$ explodes, then $X^\delta$
remains bounded on $[0,\tau^\delta_\infty)$, and consequently,
$\lim_{t\uparrow\tau^\delta_\infty}|Y^\delta_t|=\infty.$
Moreover, there exist constants $C,l>0$ such that
$$
\mathbb{E}\big[
\|X^\delta\|_{\alpha^\prime;
	[0,\tau^\delta_\infty)}
\big]
\leq C\mathbb{E}\big[\|B\|_\alpha^l+1\big].
$$

We first establish the non-explosion of $Z^\delta$ for each fixed $\delta\in(0,1]$ by using Proposition \ref{pro-2}.
By It\^o's formula and Assumption \textbf{(H4)$_2$}, we obtain
\begin{align*}
\mathbb{E}\big[|Y^\delta_{t\wedge\tau^\delta_N}|^2\big]
=&|y_0|^2+\delta^{-1}\mathbb{E}\Big[\int^{t\wedge\tau^\delta_N}_0\big(2\langle Y^\delta_s,b_2(X^\delta_s,Y^\delta_s)\rangle 
+|\sigma_2(X^\delta_s,Y^\delta_s)|^2\big)ds\Big]\\
\leq &
|y_0|^2+\delta^{-1}\mathbb{E}\Big[\int^t_0\big(-\gamma_1|Y^\delta_s|^2+C|X^\delta_s|^{\eta_1}+C\big)\mathbf{1}_{\{s\leq \tau^\delta_N\}}ds\Big]\\
\leq &
|y_0|^2+\delta^{-1}CT\big\{\mathbb{E}\big[(\|B\|^l_\alpha+1)^{\eta_1}\big]+1\big\}\\
\leq &
|y_0|^2+\delta^{-1}C,
\end{align*}  
where $C>0$ is independent of $N$ and $\delta$.
Since we stopped the time by $\tau^\delta_N$, the martingale part in the It\^o's formula is a true martingale and its expectation vanishes.
Recall that, as $N\to \infty$, $|Y^\delta_{T\wedge\tau^\delta_N}|\to \infty$ if $\tau^\delta_\infty\leq T$, and $|Y^\delta_{T\wedge\tau^\delta_N}|\to |Y^\delta_{T}|$ if $\tau^\delta_\infty=\infty$.
Applying Fatou's lemma to above inequality with $t=T$, we have 
$$
\mathbb{E}\big[\liminf_{N\to \infty}|Y^\delta_{T\wedge\tau^\delta_N}|^2\big]
\leq \liminf_{N\to \infty}\big(|y_0|^2+\delta^{-1}C\big)<\infty.
$$
Then, we see that $\mathbb{P}(\tau^\delta_\infty\leq T)=0$, and $\sup_{t\in[0,T]}\mathbb{E}\big[|Y^\delta_t|^2\big]\leq |y_0|^2+\delta^{-1}C$, that is, $Z^\delta$ does not explode on $[0,T]$ almost surely.
Furthermore, \eqref{eq-31} follows directly from Lemma \ref{lem-3} and \textbf{(HN)}.

We next establish moment estimates of arbitrary order $2\leq p<\infty$, with constants possibly depending on \(\delta\). These estimates will be used to justify the removal of the stopping time. Finally, taking $p=q$ and
exploiting the dissipativity condition \textbf{(H4)$_q$}, we derive the desired moment estimate \eqref{eq-32}, uniformly in $\delta$.

For \(2\leq p<\infty\), It\^{o}'s formula yields
\begin{align*}
\mathbb{E}\big[|Y^\delta_{t\wedge\tau^\delta_N}|^p\big]
\leq &
|y_0|^p
+\delta^{-1}\mathbb{E}\Big[\int_0^{t\wedge\tau^\delta_N}p\Big\{|Y^\delta_s|^{p-2}\langle Y^\delta_s,b_2(X^\delta_s,Y^\delta_s)\rangle+\frac12 |Y^\delta_s|^{p-2}|\sigma_2(X^\delta_s,Y^\delta_s)|^2\\
&\qquad+\frac{p-2}{2}|Y^\delta_s|^{p-4}(Y^\delta_s)^\mathrm{T}\sigma_2(X^\delta_s,Y^\delta_s)\sigma_2(X^\delta_s,Y^\delta_s)^\mathrm{T}Y^\delta_s
\Big\}ds\Big] \\
\leq &
|y_0|^p+\delta^{-1}\mathbb{E}\Big[\int_0^{t\wedge\tau^\delta_N}\frac p2|Y^\delta_s|^{p-2}\big\{2\langle Y^\delta_s,b_2(X^\delta_s,Y^\delta_s)\rangle+(p-1)|\sigma_2(X^\delta_s,Y^\delta_s)|^2\big\}ds\Big]\\
\leq &
|y_0|^p+\delta^{-1}C\mathbb{E}\Big[\int_0^{t\wedge\tau^\delta_N}|Y^\delta_s|^{p-2}\big(|Y^\delta_s|^2+|X^\delta_s|^{\eta_1\vee 2}+1\big)ds\Big]\\
\leq &
|y_0|^p+\delta^{-1}C\mathbb{E}\Big[\int_0^{t\wedge\tau^\delta_N}\big(|Y^\delta_s|^p+|X^\delta_s|^{\frac{(\eta_1\vee 2)p}2}+1\big)ds\Big]\\
\leq &
|y_0|^p+\delta^{-1}C\int_0^t\mathbb{E}\big[|Y^\delta_s|^p \mathbf{1}_{\{s\leq \tau^\delta_N\}}\big]ds+\delta^{-1}C.
\end{align*}
Here, the third-to-last inequality follows from Assumptions \textbf{(H3)} and \textbf{(H4)$_2$}, the second-to-last inequality follows from Young's inequality, and the last inequality follows from \eqref{eq-31}. 
By Gronwall's inequality, we obtain
\begin{align*}
\mathbb{E}\big[|Y^\delta_{t\wedge\tau^\delta_N}|^p\big]\leq (|y_0|^p+\delta^{-1}C)e^{\frac{CT}{\delta}},\quad t\in[0,T].
\end{align*}
Consequently, there exists a constant \(C_\delta>0\) independent of \(N\) and \(t\), such that
\begin{align*}
\mathbb{E}\big[|Y^\delta_{t\wedge\tau^\delta_N}|^p\big]\leq C_\delta,\quad t\in[0,T],\quad 2\leq p<\infty.
\end{align*}

Choosing \(r>1\) and applying the above estimate with exponent \(pr\), we obtain
$$ \sup_{N\in\mathbb{N}}\mathbb E\big[ |Y^\delta_{t\wedge\tau_N^\delta}|^{pr} \big]<\infty. $$
Hence, the family $\{|Y^\delta_{t\wedge\tau^\delta_N}|^p\}_{N\in\mathbb{N}}$ is uniformly integrable.
Since $|Y^\delta_{t\wedge\tau^\delta_N}|\to|Y^\delta_t|$ as $N\to \infty$ almost surely, 
it follows from Vitali's convergence theorem that
\[
\lim_{N\to\infty}
\mathbb{E}\big[
|Y^\delta_{t\wedge\tau^\delta_N}|^p
\big]
=
\mathbb{E}\big[|Y^\delta_t|^p\big].
\]
Consequently, 
\begin{align}\label{eq-33}
\sup_{t\in[0,T]}\mathbb{E}\big[|Y^\delta_t|^p\big]=\sup_{t\in[0,T]}\lim_{N\to\infty}\mathbb{E}\big[|Y^\delta_{t\wedge\tau^\delta_N}|^p\big]\leq C_\delta<\infty
\end{align}
for every $\delta\in(0,1]$ and $2\leq p<\infty$.

Define
\begin{align*}
F_s:=&|Y^\delta_s|^{p-2}\langle
Y^\delta_s,b_2(X^\delta_s,Y^\delta_s)\rangle +\frac12|Y^\delta_s|^{p-2}|\sigma_2(X^\delta_s,Y^\delta_s)|^2\\
&+\frac{p-2}{2}|Y^\delta_s|^{p-4}(Y^\delta_s)^\mathrm{T}\sigma_2(X^\delta_s,Y^\delta_s)\sigma_2(X^\delta_s,Y^\delta_s)^\mathrm{T}Y^\delta_s.
\end{align*}
Then, by Young's inequality and Assumptions \textbf{(H3)} and \textbf{(H5)}, we obtain
\begin{align*}
|F_s\mathbf{1}_{\{s\leq \tau^\delta_N\}}|
\leq &
|Y^\delta_s|^{p-1}C\big(|X^\delta_s|^{\eta_2}+|Y^\delta_s|^{\eta_2}+1\big)
+\frac{p-1}2 |Y^\delta_s|^{p-2}\big(|X^\delta_s|+|Y^\delta_s|+|\sigma_2(x_0,y_0)|\big)^2\\
\leq &
C\big(|X^\delta_s|^v+|Y^\delta_s|^v+1\big),
\end{align*}
where $v:=(p+\eta_2-1)\vee p\geq2$.
Due to \eqref{eq-31} and \eqref{eq-33}, the right hand side is independent of $N$ and integrable over $\Omega\times[0,t]$ with respect to the measure $\mathbb{P}\otimes ds$. More precisely,
\begin{align*}
	|F_s\mathbf{1}_{\{s\leq \tau^\delta_N\}}|\leq 
	C\big(|X^\delta_s|^v+|Y^\delta_s|^v+1\big)=:G_s,
\end{align*}
where $G\in L^1(\Omega\times[0,t],\mathbb{P}\otimes ds)$.
Since $F_s\mathbf{1}_{\{s\leq \tau^\delta_N\}}\to F_s$ as $N\to \infty$ almost surely,  the dominated convergence theorem yields $F\in L^1$ and
$$
\mathbb{E}\Big[\int_0^tF_s\mathbf{1}_{\{s\leq \tau^\delta_N\}}ds\Big]\to \mathbb{E}\Big[\int_0^tF_sds\Big] \quad\text{as }N\to \infty.
$$

Hence, combining the identity obtained from It\^o's formula with \eqref{eq-33} and the preceding convergence results, we obtain
\begin{align*}
\mathbb{E}\big[|Y^\delta_t|^p\big]
=
\lim_{N\to \infty}\mathbb{E}\big[|Y^\delta_{t\wedge\tau^\delta_N}|^p\big]
=|y_0|^p+\frac p\delta\lim_{N\to \infty}\mathbb{E}\Big[\int_0^tF_s\mathbf{1}_{\{s\leq \tau^\delta_N\}}ds\Big]
=|y_0|^p+\frac p\delta\mathbb{E}\Big[\int_0^tF_sds\Big].
\end{align*}
Differentiating both sides with respect to $t$, taking $p=q$, and applying Assumption \textbf{(H4)$_q$}, we obtain
\begin{align*}
\frac{d}{dt}\mathbb{E}\big[|Y^\delta_t|^q\big]
=&
\frac q\delta\mathbb{E}\Big[|Y^\delta_t|^{q-2}\langle
Y^\delta_t,b_2(X^\delta_t,Y^\delta_t)\rangle +\frac12|Y^\delta_t|^{q-2}|\sigma_2(X^\delta_t,Y^\delta_t)|^2\\
&\qquad+\frac{q-2}{2}|Y^\delta_t|^{q-4}(Y^\delta_t)^\mathrm{T}\sigma_2(X^\delta_t,Y^\delta_t)\sigma_2(X^\delta_t,Y^\delta_t)^\mathrm{T}Y^\delta_t\Big]\\
\leq &
\frac{q}{2\delta}\mathbb{E}\big[|Y^\delta_t|^{q-2}\big\{2\langle
Y^\delta_t,b_2(X^\delta_t,Y^\delta_t)\rangle+(q-1)|\sigma_2(X^\delta_t,Y^\delta_t)|^2\big\}\big]\\
\leq &
\frac{q}{2\delta}\mathbb{E}\big[|Y^\delta_t|^{q-2}\big(-\gamma_1|Y^\delta_t|^2+C|X^\delta_t|^{\eta_1}+C\big)\big]\\
\leq &
\frac{-q\gamma_1}{4\delta}\mathbb{E}\big[|Y^\delta_t|^q\big]+\frac C\delta
\end{align*} 
for almost every $t\in[0,T]$, 
where the last inequality follows from Young's inequality and \eqref{eq-32}.

By the comparison theorem for ordinary differential equations, we obtain
\begin{align*}
\mathbb{E}\big[|Y^\delta_t|^q\big]\leq &
|y_0|^qe^{\frac{-q\gamma_1}{4\delta}t}+\frac C\delta \frac{4\delta}{q\gamma_1}(1-e^{\frac{-q\gamma_1}{4\delta}t})
\leq 
|y_0|^q+\frac{4C}{q\gamma_1}
\end{align*}
for every $t\in[0,T]$,
Since the right-hand side is independent of $\delta$, the desired uniform moment estimate follows.
This completes the proof.
\end{proof}

For any fixed $x\in\mathbb{R}^{d_1}$, we first consider the following frozen equation:
\begin{align}\label{eq-37}
dY^{x,y}_t=b_2(x,Y^{x,y}_t)dt+\sigma_2(x,Y^{x,y}_t)d^Iw_t,\quad
Y^{x,y}_0=y\in\mathbb{R}^{d_2},\quad t\geq0. 
\end{align}
Under the basic assumptions together with \textbf{(H3)} and \textbf{(H4)$_2$}, the frozen equation \eqref{eq-37} admits a unique global solution $(Y^{x,y}_t)_{t\geq0}$ for every $x\in\mathbb{R}^{d_1}$ and $y\in\mathbb{R}^{d_2}$. The proof is analogous to that of Proposition \ref{pro-3}, with $q=2$.

For each fixed $x\in\mathbb{R}^{d_1}$,
the family $\{Y^{x,y}_\cdot\}_{y\in\mathbb{R}^{d_2}}$
defines a Markov process on $\mathbb{R}^{d_2}$.
Its associated Markov semigroup
$(P_t^x)_{t\geq0}$ is defined by
\[
P_t^x\varphi(y)
:=\mathbb{E}\big[\varphi(Y^{x,y}_t)\big],
\quad t\geq0,\quad y\in\mathbb{R}^{d_2},
\]
for every bounded Borel-measurable function
$\varphi:\mathbb{R}^{d_2}\to\mathbb{R}$.

We next establish the existence and uniqueness of an invariant probability measure for $(P_t^x)_{t\geq0}$.
To this end, we recall Lemmas \ref{lem-c1}--\ref{lem-c4} from \cite{inahama2025averaging,liu2020averaging}.

The first lemma provides a uniform-in-time moment estimate for the frozen dynamics.
In particular, it ensures the tightness of the family of time-averaged distributions.
\begin{lemma}\label{lem-c1}
Assume \textbf{(H3)} and \textbf{(H4)}$_q$ for some $2\leq q<\infty$. Then, there exists a positive constant $C_q$, independent of $t$, $x$ and $y$,such that
$$\mathbb{E}[|Y^{x,y}_t|^q] \le e^{\frac{-q \gamma_1}4 t } |y|^q + C_q (1 + |x|^{\frac{q \eta_1}2}), \quad t \ge 0, \; x \in \mathbb{R}^{d_1}, \; y \in \mathbb{R}^{d_2}.$$
Here, $\gamma_1>0$ and $\eta_1\geq 0$ are the constants in \textbf{(H4)$_q$}.
\end{lemma}
The next lemma establishes exponential contractivity of the frozen dynamics with respect to the initial condition.
\begin{lemma}\label{lem-c2}
Assume \textbf{(H3)}, \textbf{(H4)$_2$} and \textbf{(H6)}. Then, it holds that
$$\mathbb{E}[|Y_t^{x,y_1}-Y_t^{x,y_2}|^2]\leq e^{-\gamma_2t}|y_1-y_2|^2,\quad t\geq0,\:x\in\mathbb{R}^{d_1},\:y_1,y_2\in\mathbb{R}^{d_2}.$$
Here, $\gamma_2>0$ is the constant in \textbf{(H6)}.
\end{lemma}
Combining the preceding moment estimate in Lemma \ref{lem-c1} with the Krylov--Bogoliubov argument yields the existence of an invariant probability measure.
Moreover, the exponential contractivity established in Lemma \ref{lem-c2} ensures its uniqueness and provides a quantitative convergence estimate, as stated in the following lemma.
\begin{lemma}\label{lem-c3}
Assume \textbf{(H3)} and \textbf{(H4)}$_q$ for some $2\leq q<\infty$ and \textbf{(H6)}. Then, for every $x\in\mathbb{R}^{d_1}$, the semigroup $(P^x_t)_{t\geq0}$ has a unique invariant probability measure $\mu^{x}$. Furthermore, the following two estimates hold:
\begin{itemize}
\item [(i)] There exists a constant $C_q>0$, independent of $x$, such that, for every $x\in\mathbb{R}^{d_1}$,
$$\int_{\mathbb{R}^{d_2}}|z|^q\mu^x(dz)\leq C_q(1+|x|^{\frac{q\eta_1}2}).$$
Here, $\eta _{1}\geq 0$ $is$ the constant in \textbf{(H4)}$_q$.
\item [(ii)] There exists a constant $C>0$, independent of $t$, $x$, $y$ and $\varphi$, such that, for every $t\geq0$,  $x\in \mathbb{R}^{d_1}, y\in \mathbb{R}^{d_2}$ and
Lipschitz function $\varphi : \mathbb{R}^{d_2}\to \mathbb{R}$,
$$\Big|P_t^x\varphi(y)-\int_{\mathbb{R}^{d_2}}\varphi(z)\mu^x(dz)\Big|\leq C\|\varphi\|_{\mathrm{Lip}}e^{\frac{-\gamma_2}2 t}(1+|x|^{\frac{\eta_1}2}+|y|).$$
Here, $\|\varphi\|_{\mathrm{Lip}}$ is the Lipschitz constant of $\varphi$ and $\gamma_2>0$ is the constant in \textbf{(H6)}.
\end{itemize}  
\end{lemma}
Under the assumptions of Lemma \ref{lem-c3}, with the invariant probability measure
$\mu^x$ established, we define the
averaged drift coefficient by
\[
\bar{b}_1(t,x)
:=
\int_{\mathbb{R}^{d_2}}
b_1(t,x,y)\mu^x(dy),
\quad
(t,x)\in[0,T]\times\mathbb{R}^{d_1},
\]
provided that the integral is well-defined.
The corresponding averaged equation is given by
\begin{align}\label{eq-38}
	d\bar{X}_t
	=\bar{b}_1(t,\bar{X}_t)dt
	+\sigma_1(t,\bar{X}_t)dB_t,
	\quad
	\bar{X}_0=x_0\in\mathbb{R}^{d_1},
	\quad t\in[0,T],
\end{align}
where $x_0$ is the initial condition
of the slow component in the original
system \eqref{eq-29}.

To establish the regularity of the averaged drift coefficient, we first recall the following estimate on the dependence of the frozen dynamics on the parameter \(x\).
\begin{lemma}\label{lem-c4}
Assume \textbf{(H3)}, \textbf{(H4)$_{2(r\vee1)}$} and \textbf{(H7)}$_{r}$ for some $r\geq0$. Then, there exist positive
constants $C$ and $\xi$ independent of $t$, $x_1$, $x_2$ and $y$, such that
$$\mathbb{E}\big[|Y^{x_1,y}_t-Y^{x_2,y}_t|^2\big]\leq C|x_1-x_2|^2(1+|x_1|^{2\xi}+|x_2|^{2\xi}+|y|^{2(r\vee1)})$$
for all $t\geq 0$, $x_1, x_2\in \mathbb{R}^{d_1}$ and $y\in \mathbb{R}^{d_2}.$
\end{lemma}


\begin{proposition}\label{pro-5}
Assume \textbf{(H2)}, \textbf{(H3)}, \textbf{(H4)$_{2(r\vee 1)}$},  \textbf{(H6)} and \textbf{(H7)}$_{r}$ for some $r\geq 0$. Then, the map $\bar{b}_1:[0,T]\times\mathbb{R}^{d_1}\to \mathbb{R}^{d_1}$defined above is bounded and locally Lipschitz continuous with respect to the spatial variable. More precisely, we have $\|\bar{b}_1\|_\infty\leq \|b_1\|_\infty$
and there exists a constant $C>0$, independent of $x_1$, $x_2$ and $t$, such that
$$
|\bar{b}_1(t,x_1)-\bar{b}_1(t,x_2)|\leq 
C|x_1-x_2|(1+|x_1|^\xi+|x_2|^\xi),\quad x_1,x_2\in\mathbb{R}^{d_1},~t\in[0,T].
$$
Here, $\xi>0$ is the constant appearing in Lemma \ref{lem-c4}
\end{proposition}
\begin{proof}
The boundedness of $\bar{b}_1$ follows immediately from its definition. Indeed, $\|\bar{b}_1\|_\infty\leq \|b_1\|_\infty$.
Moreover, we obtain
\begin{align*}
|\bar{b}_1(t,x_1)-\bar{b}_1(t,x_2)
|
\leq &
\Big|
\int_{\mathbb{R}^{d_2}}b_1(t,x_1,z)\mu^{x_1}(dz)-\mathbb{E}\big[b_1(t,x_1,Y^{x_1,0}_s)\big]
\Big|\\
&+
\Big|
\int_{\mathbb{R}^{d_2}}b_1(t,x_2,z)\mu^{x_2}(dz)-\mathbb{E}\big[b_1(t,x_2,Y^{x_2,0}_s)\big]
\Big|\\
&+
\Big|
\mathbb{E}\big[b_1(t,x_1,Y^{x_1,0}_s)-b_1(t,x_2,Y^{x_2,0}_s)\big]
\Big|\\
\leq &
CL_2e^{\frac{-\gamma_2}2 s}(1+|x_1|^{\frac{\eta_1}2}+|x_2|^{\frac{\eta_1}2})
+L_2\big(|x_1-x_2|+\mathbb{E}\big[|Y^{x_1,0}_s-Y^{x_2,0}_s|\big]\big)
\\\leq &CL_2e^{\frac{-\gamma_2 }2 s}(1+|x_1|^{\frac{\eta_1}2}+|x_2|^{\frac{\eta_1}2})+CL_2|x_1-x_2|(1+|x_1|^\xi+|x_2|^\xi)
.
\end{align*}
Here, 
the second inequality follows from Lemma \ref{lem-c3} and \textbf{(H2)}, and the third from Lemma \ref{lem-c4} and H\"older's inequality.
Letting $s\to\infty$, we obtain the desired estimate.
\end{proof}

\begin{remark}\label{rem-3}
	Assume \textbf{(HN)}, \textbf{(H1)}--\textbf{(H3)}, \textbf{(H4)$_{2(r\vee 1)}$}, \textbf{(H6)} and \textbf{(H7)}$_{r}$ for some $r\geq 0$.
	By Propositions \ref{pro-4} and \ref{pro-5}, the averaged equation \eqref{eq-38} admits a unique global solution. Moreover, there exist constants $C,l>0$ such that
	$$\mathbb{E}[\|\bar{X}\|^p_{\alpha^\prime}]\leq C\mathbb{E}[\|B\|^{lp}_\alpha+1]<\infty,$$
	for any $1\leq p <\infty$.
\end{remark}

\subsection{Statement of AP}
We are now in a position to establish the AP. The following theorem shows that, as $\delta\to0$, the slow component $X^\delta$ converges to the solution of the averaged equation in $L^p$ with respect to the H\"older topology.
\begin{theorem}\label{th-2}
Assume \textbf{(HN)}, \textbf{(H1)}--\textbf{(H3)}, \textbf{(H4)}$_q$, \textbf{(H5)}, \textbf{(H6)} and \textbf{(H7)}$_{r}$ for some $2\leq q<\infty$ and $0\leq r<\infty$ satisfying $q>2r$.
Then, for every $1\leq p<\infty$ and $\alpha^\prime\in(1/2, \alpha)$, we have 
$$
\lim_{\delta\to 0}\mathbb{E}[\|X^\delta-\bar{X}\|^p_{\alpha^\prime}]=0.
$$
\end{theorem}

\begin{remark}\label{rem-4}
The prominent example of stochastic processes satisfying \textbf{(HN)} is FBM with Hurst parameter $H \in (1/2, 1)$.
The study of such slow-fast systems driven by 
mixed FBM with $H \in (1/2, 1)$ was initiated by \cite{hairer2020averaging}, followed by \cite{pei2023averaging} and has since attracted increasing attention \cite{pei2020stochastic, shen2022averaging,liu2020averaging,li2025averaging,gailus2025large, inahama2023large, xu2025large,yang2026fluctuation,bourguin2021typical}. 
	
The majority of the papers along this research line use fractional calculus (see e.g. \cite{pei2023averaging,pei2020stochastic, shen2022averaging,liu2020averaging}), while a few use the stochastic sewing lemma (see e.g. \cite{hairer2020averaging,li2025averaging}).
Our approach is based on neither of these methods, but instead on the controlled rough path framework.
We believe that our method is more accessible to readers who are not familiar with fractional calculus.
Compared with the stochastic sewing approach, our method is based on Young--rough integration and allows for a more direct pathwise analysis without relying on the stochastic sewing lemma.

Compared with preceding works on AP for this kind of slow-fast system, our result (Theorem \ref{th-2}) has the following features:
\begin{itemize}
\item 
The convergence takes place at the path space level and the topology of the path space is the H\"older topology with suitable H\"older exponent. Moreover, our strong AP holds in $L^p$-norm with any $p \in [1,\infty)$. 
\item 
Our driving process of the slow equation is quite general.
\item
Our work basically follows \cite{pei2023averaging}, but our assumptions on the coefficients are weaker than those in \cite{pei2023averaging}. (It should be noted that those in \cite{hairer2020averaging} are very different from those in other papers on this research topic. For example, the assumptions in  \cite{hairer2020averaging} and \cite{pei2023averaging} are mutually exclusive.)
\end{itemize}

It is worth noting that strong convergence 
($L^1$ or $L^p$ convergence) may fail in general for fully coupled slow--fast systems when the diffusion coefficient of the slow equation depends on the fast variable. A counterexample in the Brownian case illustrating this phenomenon is provided in \cite{liu2010strong}. In view of this difficulty, we focus on systems in which the diffusion coefficient of the slow equation is independent of the fast variable.
 
On the other hand,
under a substantially different set of assumptions, Hairer and Li \cite{hairer2020averaging} established an AP for fully coupled slow--fast systems driven by mixed FBM,  proving convergence in probability using the stochastic sewing lemma. Subsequently, Li, Gao and Qu \cite{li2025averaging} investigated the AP in an infinite-dimensional setting with additive noise in the fast equation.
\end{remark}

\subsection{Proof of AP}
In this subsection, we prove Theorem \ref{th-2} by employing Khasminskii's time discretization method.

Introduce an auxiliary parameter $\Delta$ satisfying $0<\delta<\Delta\leq 1\wedge T$.
We will eventually let $\Delta\to0$ and $\delta/\Delta\to0$. Partition the time interval $[0,T]$ into subintervals of  length $\Delta$, except possibly for the last subinterval.
For each $s\geq 0$, set $s(\Delta):= \lfloor\frac s\Delta\rfloor \Delta$. Define the auxiliary process $\hat{Y}^\delta$
by
\begin{align}\label{eq-41}
\hat{Y}^\delta_t=y_0+\frac1\delta\int_0^tb_2(X^\delta_{s(\Delta)},\hat{Y}_s)ds
+\frac{1}{\sqrt\delta}\int_0^t\sigma_2(X^\delta_{s(\Delta)},\hat{Y}_s)d^Iw_s,\quad t\in[0,T].
\end{align}
The following lemma provides a uniform moment estimate for the auxiliary process $\hat{Y}^\delta$.

\begin{lemma}\label{lem-4}
Assume \textbf{(HN)}, \textbf{(H1)}--\textbf{(H3)}, \textbf{(H4)}$_q$ and \textbf{(H5)} for some $2\leq q<\infty$.
Then, for every $\delta$ and $\Delta$ satisfying $0<\delta<\Delta\leq 1$, equation \eqref{eq-41} admits a unique global solution $\hat{Y}^\delta$. Moreover, 
$$\sup_{0<\delta<\Delta\leq1}\sup_{t\in[0,T]}\mathbb{E}[|\hat{Y}^\delta_t|^q]<\infty.$$
\end{lemma}
\begin{proof}
The proof follows essentially the same argument as that of Proposition \ref{pro-3}.
\end{proof}

\begin{lemma}\label{lem-5}
Assume \textbf{(HN)}, \textbf{(H1)}--\textbf{(H3)}, \textbf{(H4)}$_q$, \textbf{(H5)}, \textbf{(H6)} and \textbf{(H7)}$_{r}$ for some $2\leq q<\infty$ and $0\leq r<\infty$ satisfying $q>2r$.
Then, there exists a constant $C>0$ independent of $\delta$ and $\Delta$, such that 
$$\sup_{0<\delta<\Delta\leq 1}\sup_{t\in[0,T]}\mathbb{E}[|Y^\delta_t-\hat{Y}^\delta_t|^2]\leq C\Delta^{2\alpha^\prime}.$$
\end{lemma}
\begin{proof}
By It\^o's formula,
\begin{align*}
\mathbb{E}\big[|Y^\delta_t-\hat{Y}^\delta_t|^2\big]
=&
\frac2\delta\mathbb{E}\Big[
\int_0^t \big\langle b_2(X^\delta_s,Y^\delta_s)-b_2(X^\delta_{s(\Delta)},\hat{Y}^\delta_s),Y^\delta_s-\hat{Y}^\delta_s\big\rangle ds
\Big]\\
&+
\frac1\delta\mathbb{E}\Big[
\int_0^t\big|\sigma_2(X^\delta_s,Y^\delta_s)-\sigma_2(X^\delta_{s(\Delta)},\hat{Y}^\delta_s)\big|^2 ds
\Big].
\end{align*}
Consequently, for almost every $t\in[0,T]$, 
\begin{align*}
\frac{d}{dt}\mathbb{E}[|Y^\delta_t-\hat{Y}^\delta_t|^2]
\leq &
\frac1\delta
\mathbb{E}\big[
2\langle b_2(X^\delta_t,Y^\delta_t)-b_2(X^\delta_t,\hat{Y}^\delta_t),Y^\delta_t-\hat{Y}^\delta_t\rangle
+|\sigma_2(X^\delta_t,Y^\delta_t)-\sigma_2(X^\delta_t,\hat{Y}^\delta_t) |^2
\big]\\
&+
\frac2\delta 
\mathbb{E}\big[
\langle b_2(X^\delta_t,\hat{Y}^\delta_t)-b_2(X^\delta_{t(\Delta)},\hat{Y}^\delta_t),Y^\delta_t-\hat{Y}^\delta_t\rangle
\big]\\
&+
\frac2\delta 
\mathbb{E}\big[
\langle \sigma_2(X^\delta_t,Y^\delta_t)-\sigma_2(X^\delta_t,\hat{Y}^\delta_t), \sigma_2(X^\delta_t,\hat{Y}^\delta_t)-\sigma_2(X^\delta_{t(\Delta)},\hat{Y}^\delta_t)\rangle
\big]\\
&+
\frac1\delta\mathbb{E}\Big[
\big|\sigma_2(X^\delta_t,\hat{Y}^\delta_t)-\sigma_2(X^\delta_{t(\Delta)},\hat{Y}^\delta_t)\big|^2 
\Big]\\
=:&
I_1+I_2+I_3+I_4.
\end{align*}
Throughout the following, $C$ denotes a positive constant independent of $t$, $\delta$ and $\Delta$.
By the dissipativity condition in \textbf{(H6)}, we have 
\begin{align*}
I_1
\leq &
\frac{-\gamma_2}{\delta}\mathbb{E}[|Y^\delta_t-\hat{Y}^\delta_t|^2].
\end{align*}
By \textbf{(H7)}$_{r}$, Proposition \ref{pro-3} and Lemma \ref{lem-4}, we have 
\begin{align*}
I_2
\leq &
\frac{C}{\delta}\mathbb{E}\big[
\big(1+|X^\delta_t|^{\eta_3}+|\hat{Y}^\delta_t|^{r}\big)|X^\delta_t-X^\delta_{t(\Delta)}| |Y^\delta_t-\hat{Y}^\delta_t|
\big]\\
\leq &
\frac{\gamma_2}{4\delta} \mathbb{E}[|Y^\delta_t-\hat{Y}^\delta_t|^2]
+\frac C\delta
\mathbb{E}\big[
\big\{1+|X^\delta_t|^{2\eta_3}+|\hat{Y}^\delta_t|^{2r}\big\}|X^\delta_t-X^\delta_{t(\Delta)}|^2
\big]\\
\leq &
\frac{\gamma_2}{4\delta} \mathbb{E}[|Y^\delta_t-\hat{Y}^\delta_t|^2]
+\frac C\delta
\mathbb{E}\big[
1+|X^\delta_t|^{4\eta_3}\big]^{\frac12}
\mathbb{E}\big[|X^\delta_t-X^\delta_{t(\Delta)}|^4\big]^{\frac12}
\\
&+
\frac C\delta
\mathbb{E}\big[
|\hat{Y}^\delta_t|^q\big]^{\frac{2r}q}
\mathbb{E}\big[|X^\delta_t-X^\delta_{t(\Delta)}|^{\frac{2q}{q-2r}}\big]^{\frac{q-2r}q}\\
\leq &
\frac{\gamma_2}{4\delta} \mathbb{E}\big[|Y^\delta_t-\hat{Y}^\delta_t|^2\big]
+\frac C\delta \Delta^{2\alpha^\prime},
\end{align*}
where the second and third inequalities follow from Young's and H\"older's inequalities, respectively. 
Using the Lipschitz continuity of $\sigma_2$ and Young's inequality, we have 
\begin{align*}
I_3+I_4
\leq &
\frac2\delta \|\sigma_2\|^2_{C^1_b}\mathbb{E}\big[
|Y^\delta_t-\hat{Y}^\delta_t| |X^\delta_t-X^\delta_{t(\Delta)}|+|X^\delta_t-X^\delta_{t(\Delta)}|^2
\big]\\
\leq &
\frac{\gamma_2}{4\delta} \mathbb{E}\big[|Y^\delta_t-\hat{Y}^\delta_t|^2\big]
+\frac C\delta \Delta^{2\alpha^\prime}.
\end{align*}
Combining the above estimates, we obtain
\begin{align*}
\frac{d}{dt}\mathbb{E}[|Y^\delta_t-\hat{Y}^\delta_t|^2]
\leq &
\frac{-\gamma_2}{2\delta}
\mathbb{E}[|Y^\delta_t-\hat{Y}^\delta_t|^2]
+\frac C\delta \Delta^{2\alpha^\prime} 
\end{align*}
for almost every $t\in[0,T]$.
Since $Y^\delta_0=\hat{Y}^\delta_0=y_0$, an application of the differential form of Gronwall's inequality yields
\begin{align*}
\mathbb{E}[|Y^\delta_t-\hat{Y}^\delta_t|^2]\leq 
\frac{C\Delta^{2\alpha^\prime}}{\gamma_2}(1-e^{\frac{-\gamma_2}{2\delta}t})
\leq C\Delta^{2\alpha^\prime}
\end{align*}
for all $t\in[0,T]$ and $\delta\in(0,\Delta)$.
This completes the proof.
\end{proof}

Define 
\begin{align*}
M_t:=&\int_0^t\big\{
b_1(s,X^\delta_s,Y^\delta_s)-b_1(s,X^\delta_s,\hat{Y}^\delta_s)
\big\}ds\\
&+\int_0^t\big\{
b_1(s,X^\delta_s,\hat{Y}^\delta_s)-b_1(s,X^\delta_{s(\Delta)},\hat{Y}^\delta_s)
\big\}ds\\
&+\int_0^t\big\{
b_1(s,X^\delta_{s(\Delta)},\hat{Y}^\delta_s)-\bar{b}_1(s,X^\delta_{s(\Delta)})
\big\}ds\\
&+\int_0^t\big\{
\bar{b}_1(s,X^\delta_{s(\Delta)})-\bar{b}_1(s,X^\delta_s)
\big\}ds
.
\end{align*}
Then, by \eqref{eq-29} and \eqref{eq-38}, we have
\begin{align*}
X^\delta_t-\bar{X}_t
=&
\int_0^t\big\{
b_1(s,X^\delta_s,Y^\delta_s)-\bar{b}_1(s,\bar{X}_s)
\big\}ds
+
\int_0^t\big\{
\sigma_1(s,X^\delta_s)-\sigma_1(s,\bar{X}_s)
\big\}dB_s\\
=&
M_t+\int_0^t\big\{
\bar{b}_1(s,X^\delta_s)-\bar{b}_1(s,\bar{X}_s)
\big\}ds+\int_0^t\big\{
\sigma_1(s,X^\delta_s)-\sigma_1(s,\bar{X}_s)
\big\}dB_s.
\end{align*}

In the following lemma, we estimate the H\"{o}lder norms of the first three terms in the decomposition of $M_t$.
\begin{lemma}\label{lem-14}
Assume \textbf{(HN)}, \textbf{(H1)}--\textbf{(H3)}, \textbf{(H4)}$_q$, \textbf{(H5)}, \textbf{(H6)} and  \textbf{(H7)}$_{r}$,  for some $2\leq q<\infty$ and $0\leq r<\infty$ satisfying $q>2r$. Let $\gamma\in(1/2,1)$. 
Then, there exists a constant  $C>0$ independent of $\delta$ and $\Delta$, such that 
\begin{align*}
	\mathbb{E}\Big[
	\Big\|
	\int_0^\cdot\big\{
	b_1(s,X^\delta_s,Y^\delta_s)-b_1(s,X^\delta_s,\hat{Y}^\delta_s)
	\big\} ds 
	\Big\|^2_\gamma\Big]
	&\leq 
	C\Delta^{2\alpha^\prime(1-\gamma)},\\
	\mathbb{E}\Big[
	\Big\|
	\int_0^\cdot\big\{
	b_1(s,X^\delta_s,\hat{Y}^\delta_s)-b_1(s,X^\delta_{s(\Delta)},\hat{Y}^\delta_s)
	\big\} ds 
	\Big\|^2_1
	\Big]
	&\leq 
	C\Delta^{2\alpha^\prime},\\
	\mathbb{E}\Big[
	\Big\|
	\int_0^\cdot\big\{
	\bar{b}_1(s,X^\delta_{s(\Delta)})-\bar{b}_1(s,X^\delta_s)
	\big\} ds 
	\Big\|^2_1
	\Big]
	&\leq 
	C\Delta^{2\alpha^\prime}.
\end{align*}
\end{lemma}

\begin{proof}
Since $b_1$ is bounded and globally Lipschitz continuous from \textbf{(H2)}, 
\begin{align*}
&\mathbb{E}\Big[
\Big\|
\int_0^\cdot\big\{
b_1(s,X^\delta_s,Y^\delta_s)-b_1(s,X^\delta_s,\hat{Y}^\delta_s)
\big\} ds 
\Big\|^2_\gamma\Big]\\
&\leq 
\mathbb{E}\Big[\Big\{
\sup_{0\leq s<t\leq T}
\frac{\int_s^t \big(2\|b_1\|_\infty\big)^\gamma|b_1(u,X^\delta_u,Y^\delta_u)-b_1(u,X^\delta_u,\hat{Y}^\delta_u)|^{1-\gamma}du}{|t-s|^\gamma}
\Big\}^2\Big]\\
&\leq 
C
\mathbb{E}\Big[
\sup_{0\leq s<t\leq T}
\Big(\int_s^t |Y^\delta_u-\hat{Y}^\delta_u|du\Big)^{2(1-\gamma)}
\Big]\\
&\leq 
C
\mathbb{E}\Big[
\Big(\int_0^T |Y^\delta_t-\hat{Y}^\delta_t|^2dt\Big)^{1-\gamma}
\Big]\\
&\leq 
C
\mathbb{E}\Big[
\int_0^T |Y^\delta_t-\hat{Y}^\delta_t|^2dt
\Big]^{1-\gamma}\\
&\leq 
C
\sup_{t\in[0,T]}\mathbb{E}\Big[
 |Y^\delta_t-\hat{Y}^\delta_t|^2
\Big]^{1-\gamma}\\
&\leq C\Delta^{2\alpha^\prime(1-\gamma)},
\end{align*}
where we have used H\"older's inequality in the second and third inequalities, Jensen's inequality in the fourth inequality, Fubini's theorem in the fifth inequality, and Lemma \ref{lem-5} in the last inequality.

By the global Lipschitz continuity of $b_1$, the second term can be estimated as follows
\begin{align*}
\mathbb{E}\Big[
\Big\|
\int_0^\cdot\big\{
b_1(s,X^\delta_s,\hat{Y}^\delta_s)-b_1(s,X^\delta_{s(\Delta)},\hat{Y}^\delta_s)
\big\} ds 
\Big\|^2_1
\Big]
\leq &
C
\mathbb{E}\Big[\sup_{t\in[0,T]}
\big|
b_1(t,X^\delta_t,\hat{Y}^\delta_t)-b_1(t,X^\delta_{t(\Delta)},\hat{Y}^\delta_t)
\big|^2
\Big]\\
\leq &
C \mathbb{E}\Big[\sup_{t\in[0,T]}|X^\delta_t-X^\delta_{t(\Delta)}|^2\Big]\\
\leq &
C\Delta^{2\alpha^\prime},
\end{align*}
where the last inequality follows from the uniform $\alpha^\prime$-H\"older estimate for $X^\delta$.
The third term requires a slightly more delicate argument since the averaged coefficient $\bar{b}_1$ may not be globally Lipschitz continuous. Since $q>2r$, Proposition \ref{pro-5} yields
\begin{align*}
\mathbb{E}\Big[
\Big\|
\int_0^\cdot\big\{
\bar{b}_1(s,X^\delta_{s(\Delta)})-\bar{b}_1(s,X^\delta_s)
\big\} ds 
\Big\|^2_1
\Big]
\leq &
C
\mathbb{E}\Big[\sup_{t\in[0,T]}
\big|
\bar{b}_1(t,X^\delta_{t(\Delta)})-\bar{b}_1(t,X^\delta_t)
\big|^2
\Big]\\
\leq &
C 
\mathbb{E}\Big[\sup_{t\in[0,T]}
\big|X^\delta_{t(\Delta)}-X^\delta_t\big|^2\big(1+|X^\delta_{t(\Delta)}|^{2\xi}+|X^\delta_t|^{2\xi}\big)
\Big]\\
\leq &
C 
\mathbb{E}\Big[\sup_{t\in[0,T]}\big|X^\delta_{t(\Delta)}-X^\delta_t\big|^4\Big]^{\frac12}
 \mathbb{E}\Big[1+\sup_{t\in[0,T]}|X^\delta_t|^{4\xi}\Big]^{\frac12}\\
\leq& C\Delta^{2\alpha^\prime},
\end{align*} 
where $\xi>0$ is the constant appearing in Proposition \ref{pro-5}. Here, the third inequality
follows from H\"older's inequality, while the last inequality follows from the uniform $\alpha^\prime$-H\"older moment estimate for $X^\delta$. Combining the estimates for the three terms completes the proof.
\end{proof}

\begin{lemma}\label{lem-6}
Assume \textbf{(H2)}, \textbf{(H3)}, \textbf{(H4)}$_q$ and \textbf{(H6)} for some $2\leq q<\infty$.
Then, there exists a constant $C>0$, independent of $x$, $y$ and $t$, such that, for all $(x,y)\in\mathbb{R}^{d_1+d_2}$ and $t\geq 0$,
\begin{align*}
\mathbb{E}\Big[
\Big|
\int_0^t\big\{
b_1(s,x,Y^{x,y}_s)-\bar{b}_1(s,x)
\big\} ds \Big|^2
\Big]
\leq &
C\big(1+|x|^{\frac{\eta_1}2}+|y|\big)t.
\end{align*}
Here, $\eta_1>0$ is the constant appearing in \textbf{(H4)}$_q$.
\end{lemma}
\begin{proof}
Let $\hat{\mathbb{E}}^{x,y}$ denote the expectation with respect to the law of the process $Y^{x,y}$ on $\mathcal{C}([0,\infty),\mathbb{R}^{d_2})$. Under this law, $(Y_t)_{t\geq0}$ denotes the
canonical coordinate process. For $0\leq s\leq u$, by the Markov property, we have
\begin{align*}
&\big|
\mathbb{E}\big[
\langle
b_1(u,x,Y^{x,y}_u)-\bar{b}_1(u,x), b_1(s,x,Y^{x,y}_s)-\bar{b}_1(s,x)
\rangle
\big]
\big|\\
&=
\big|
\hat{\mathbb{E}}^{x,y}\big[
\langle
b_1(u,x,Y_u)-\bar{b}_1(u,x), b_1(s,x,Y_s)-\bar{b}_1(s,x)
\rangle
\big]
\big|\\
&=
\big|
\hat{\mathbb{E}}^{x,y}\big[\hat{\mathbb{E}}^{x,y}\big[
\langle
b_1(u,x,Y_u)-\bar{b}_1(u,x), b_1(s,x,Y_s)-\bar{b}_1(s,x)
\rangle
|\mathcal{F}_s\big]\big]
\big|\\
&=
\big|
\hat{\mathbb{E}}^{x,y}\big[
\langle
\hat{\mathbb{E}}^{x,y}\big[
b_1(u,x,Y_u)-\bar{b}_1(u,x)|\mathcal{F}_s\big]
, b_1(s,x,Y_s)-\bar{b}_1(s,x)
\rangle
\big]
\big|\\
&=
\big|
\hat{\mathbb{E}}^{x,y}\big[
\langle
\hat{\mathbb{E}}^{x,Y_s}\big[
b_1(u,x,Y_{u-s})-\bar{b}_1(u,x)\big]
, b_1(s,x,Y_s)-\bar{b}_1(s,x)
\rangle
\big]
\big|\\
&\leq 
2\|b_1\|_\infty \hat{\mathbb{E}}^{x,y} \big[
\big|
\hat{\mathbb{E}}^{x,Y_s}[
b_1(u,x,Y_{u-s})]-\bar{b}_1(u,x)
\big|
\big]\\
&=
2\|b_1\|_\infty \hat{\mathbb{E}}^{x,y} \big[
\big|
P^x_{u-s}b_1(u,x,Y_s)-\int_{\mathbb{R}^{d_2}}b_1(u,x,z)\mu^x(dz)
\big|
\big]\\
&\leq 
C e^{\frac{-\gamma_2(u-s)}2}
\mathbb{E} \big[
1+|x|^{\frac{\eta_1}2}+|Y^{x,y}_s|
\big]\\
&\leq 
C e^{\frac{-\gamma_2(u-s)}2}\big(1+|x|^{\frac{\eta_1}2}+|y|
\big),
\end{align*}
where $C>0$ is a constant independent of $x$, $y$, $s$ and $u$. The penultimate inequality follows from Lemma \ref{lem-c3},
while the last follows from Lemma \ref{lem-c1} and H\"older's inequality.
Consequently, by Fubini's theorem, we obtain
\begin{align*}
&\mathbb{E}\Big[
\Big|
\int_0^t \big\{
b_1(s,x,Y^{x,y}_s)-\bar{b}_1(s,x)
\big\}ds
\Big|^2
\Big]\\
&=
2\int_0^t ds \int_s^t du 
\mathbb{E}\big[
\langle
b_1(s,x,Y^{x,y}_s)-\bar{b}_1(s,x), b_1(u,x,Y^{x,y}_u)-\bar{b}_1(u,x)
\rangle
\big]\\
&\leq 
C(1+|x|^{\frac{\eta_1}2}+|y|)\int_0^t\int_s^te^{\frac{-\gamma_2(u-s)}2} duds  \\
&\leq 
C(1+|x|^{\frac{\eta_1}2}+|y|)t.
\end{align*}
This completes the proof.
\end{proof}

\begin{lemma}\label{lem-7}
Assume \textbf{(HN)}, \textbf{(H1)}--\textbf{(H3)}, \textbf{(H4)}$_q$, \textbf{(H5)}, \textbf{(H6)} and \textbf{(H7)}$_{r}$ for some $2\leq q<\infty$ and $0\leq r<\infty$ satisfying $q>2r$.
Let $\gamma\in(1/2,1)$. 
Then, there exists a constant  $C>0$ independent of $\delta$ and $\Delta$, such that 
\begin{align*}
	\mathbb{E}\Big[
	\Big\|
	\int_0^\cdot\big\{
	b_1(s,X^\delta_{s(\Delta)},\hat{Y}^\delta_s)-\bar{b}_1(s,X^\delta_{s(\Delta)})
	\big\} ds 
	\Big\|^2_\gamma
	\Big]
	\leq 
	C\big(\Delta^{2(1-\gamma)}+ \Delta^{-2\gamma}\delta\big).
\end{align*}
\end{lemma}
\begin{proof}
Let 
$$G_t:=\int_0^t\big\{
b_1(s,X^\delta_{s(\Delta)}, \hat{Y}^\delta_s)-\bar{b}_1(s,X^\delta_{s(\Delta)})
\big\}ds.$$
For $0\leq s<t\leq T$, denote
$G^1_{s,t}:=G_t-G_s$. Throughout the proof, $C$ denotes a positive
constant independent of $s,t,\delta,\Delta$ and $k$, which may change
from line to line.

We first consider the case $0<t-s<2\Delta$. Since both $b_1$ and
$\bar b_1$ are bounded, we have
\begin{align}\label{eq-42}
|G^1_{s,t}|
\leq 
\big(\|b_1\|_\infty+\|\bar{b}_1\|_\infty\big)(t-s)
\leq 
2\|b_1\|_\infty(2\Delta)^{1-\gamma}(t-s)^\gamma.
\end{align}

Next, suppose that $t-s\geq 2\Delta$. Then, there exists at least one $k\in\mathbb{N}$ such that $[k\Delta,(k+1)\Delta]\subset[s,t]$.
For simplicity, we assume that neither $s/\Delta$ nor $t/\Delta$ is an integer. The cases where either endpoint lies on the grid can be treated similarly and are slightly simpler.
Decomposing $[s,t]$ into two boundary intervals and a collection of all intervals of length $\Delta$, and applying Cauchy-Schwarz's inequality, we obtain
\begin{align*}
|G^1_{s,t}|^2
= &
\Big|
G^1_{s,(\lfloor\frac s\Delta\rfloor+1)\Delta}+\sum_{j=\lfloor\frac s\Delta\rfloor+1}^{\lfloor\frac t\Delta\rfloor-1}G^1_{j\Delta,(j+1)\Delta}+G^1_{\lfloor\frac t\Delta\rfloor\Delta,t}
\Big|^2\\
\leq &
C\Delta^2 +2\Big|
\sum_{j=\lfloor\frac s\Delta\rfloor+1}^{\lfloor\frac t\Delta\rfloor-1}G^1_{j\Delta,(j+1)\Delta}
\Big|^2\\
\leq &
C\Delta^2+2\sum_{j=\lfloor\frac s\Delta\rfloor+1}^{\lfloor\frac t\Delta\rfloor-1} 1^2\times \sum_{j=\lfloor\frac s\Delta\rfloor+1}^{\lfloor\frac t\Delta\rfloor-1}|G^1_{j\Delta,(j+1)\Delta}|^2\\
\leq &
C\Delta^2+\frac{2(t-s)}\Delta\sum_{j=\lfloor\frac s\Delta\rfloor+1}^{\lfloor\frac t\Delta\rfloor-1}|G^1_{j\Delta,(j+1)\Delta}|^2,
\end{align*}
where the second inequality follows from \eqref{eq-42}, which bounds the sum of the two boundary terms by $C\Delta$. 
Dividing both sides by $|t-s|^{2\gamma}$ gives
\begin{align}\label{eq-43}
\frac{|G^1_{s,t}|^2}{|t-s|^{2\gamma}}
\leq &
C\Delta^{2(1-\gamma)}+\frac{2(t-s)^{1-2\gamma}}\Delta\sum_{j=\lfloor\frac s\Delta\rfloor+1}^{\lfloor\frac t\Delta\rfloor-1}|G^1_{j\Delta,(j+1)\Delta}|^2.
\end{align}
Therefore, combining \eqref{eq-42} and \eqref{eq-43} and noting that $1-2\gamma<0$, we obtain
\begin{align*}
\mathbb{E}\big[\|G\|^2_\gamma\big]
\leq &
C\Delta^{2(1-\gamma)}
+C   \Delta^{-(1+2\gamma)}\max_{0\leq k\leq \lfloor\frac T\Delta\rfloor-1}\mathbb{E}\Big[\Big|
\int_{k\Delta}^{(k+1)\Delta} \big\{
b_1(s,X^\delta_{s(\Delta)},\hat{Y}^\delta_s)-\bar{b}_1(s,X^\delta_{s(\Delta)})
\big\}
\Big|^2\Big].
\end{align*}

Let $\bar{w}=(\bar{w}_t)_{t\geq 0}$ be another standard $r$-dimensional BM which is independent of $(B,w)$ and let $(\bar{Y}^{x,y}_t)_{t\geq 0}$ denote the unique solution of the frozen SDE driven by $\bar{w}$.
By Lemma \ref{lem-6}, for every$(x,y)\in\mathbb R^{d_1+d_2}$ and $t\geq0$,
\begin{align*}
\mathbb{E}\Big[
\Big|
\int_0^t \big\{
b_1(s,x,\bar{Y}^{x,y}_s)-\bar{b}_1(s,x)
\big\}ds
\Big|^2
\Big]
\leq &
C\big(1+|x|^{\frac{\eta_1}2}+|y|\big)t,
\end{align*}
where the positive constant $C$ is independent of $x$, $y$ and $t$.

For a fixed $k$, define $\bar{Y}_t:=\bar{Y}^{x,y}_t|_{(x,y)=(X^\delta_{k\Delta},\hat{Y}^\delta_{k\Delta})}$.
Since $\bar{w}$ is independent of $(B,w)$, it is also independent of $(X^\delta_{k\Delta},\hat{Y}^\delta_{k\Delta})$.
Then, $(\bar{Y}_{t/\delta})_{t\geq 0}$ satisfies 
\begin{align}\label{eq-44}
\bar{Y}_{t/\delta}
=&\hat{Y}^\delta_{k\Delta}+\int_0^{\frac t\delta} b_2(X^\delta_{k\Delta},\bar{Y}_s)ds +\int_0^{\frac t\delta} \sigma_2(X^\delta_{k\Delta},\bar{Y}_s)d^I\bar{w}_s\notag\\
=&\hat{Y}^\delta_{k\Delta}+\frac1\delta\int_0^t b_2(X^\delta_{k\Delta},\bar{Y}_{s/\delta})ds +\frac 1{\sqrt\delta}\int_0^t \sigma_2(X^\delta_{k\Delta},\bar{Y}_{s/\delta})d^I(\sqrt\delta\bar{w}_{s/\delta}),
\end{align}
where $(\sqrt\delta\bar{w}_{t/\delta})_{t\geq 0}$ is again a standard $r$-dimensional BM independent of $(X^\delta_{k\Delta},\hat{Y}^\delta_{k\Delta})$.
Comparing \eqref{eq-44} with \eqref{eq-41},
we see that, on $[0,\Delta)$,
$(\hat{Y}^\delta_{k\Delta+t})_{t\in[0,\Delta)}$
and $(\bar{Y}_{t/\delta})_{t\in[0,\Delta)}$
satisfy the same frozen SDE in distribution,
with the same initial condition
$\hat{Y}^\delta_{k\Delta}$ and frozen parameter
$X^\delta_{k\Delta}$.
Hence, by uniqueness in law of the frozen SDE,
\[
\big(
X^\delta_{k\Delta},
(\hat{Y}^\delta_{k\Delta+t})_{t\in[0,\Delta)}
\big)
\stackrel{\mathrm{law}}{=}
\big(
X^\delta_{k\Delta},
(\bar{Y}_{t/\delta})_{t\in[0,\Delta)}
\big).
\]
By the same argument as in the proof of Lemma \ref{lem-6}, after
rescaling time, we have
\begin{align*}
	\mathbb{E}\Big[
	\Big|
	\int_0^t\big\{
	b_1(\delta s,x,Y^{x,y}_s)-\bar{b}_1(\delta s,x)
	\big\} ds \Big|^2
	\Big]
	\leq &
	C(1+|x|^{\frac{\eta_1}2}+|y|)t,
\end{align*}
for some constant $C>0$ independent of $x$, $y$ and $t$.
Hence, 
\begin{align*}
&\mathbb{E}\Big[
\Big|
\int_{k\Delta}^{(k+1)\Delta}\big\{
b_1(s,X^\delta_{k\Delta},\hat{Y}^\delta_s)-\bar{b}_1(s,X^\delta_{k\Delta})
\big\}ds
\Big|^2
\Big]\\
&=
\mathbb{E}\Big[
\Big|
\int_0^\Delta\big\{
b_1(s,X^\delta_{k\Delta},\bar{Y}_{s/\delta})-\bar{b}_1(s,X^\delta_{k\Delta})
\big\}ds
\Big|^2
\Big]\\
&=
\delta^2 \mathbb{E}\Big[
\Big|
\int_0^{\frac\Delta\delta}\big\{
b_1(\delta s,X^\delta_{k\Delta},\bar{Y}_s)-\bar{b}_1(\delta s,X^\delta_{k\Delta})
\big\}ds
\Big|^2
\Big]\\
&\leq 
\delta^2 C \mathbb{E}\Big[
\big(1+|X^\delta_{k\Delta}|^{\frac{\eta_1}2}+|\hat{Y}^\delta_{k\Delta}|\big)\frac\Delta\delta
\Big]\\
&\leq 
C\delta\Delta,
\end{align*}
where the last inequality follows from Proposition \ref{pro-3} and
Lemma \ref{lem-4}, together with H\"older's inequality.
Substituting this estimate into the preceding bound yields
\begin{align*}
\mathbb{E}\big[\|G\|^2_\gamma\big]\leq C\big(\Delta^{2(1-\gamma)}+\Delta^{-2\gamma}\delta\big).
\end{align*}
This completes the proof.
\end{proof}

\begin{proof}[Proof of Theorem \ref{th-2}]

Combining Lemmas \ref{lem-14}--\ref{lem-7}, we obtain
\begin{align*}
\mathbb{E}\big[\|M\|^2_{\alpha^\prime}\big]\leq C\big(\Delta^{2\alpha^\prime(1-\alpha^\prime)}+\Delta^{-2\alpha^\prime}\delta\big).
\end{align*}
If we choose $\Delta:=\delta^{\frac12}$ for example, then $\|M\|_{\alpha^\prime}$ converges to 0 in $L^2$-sense as $\delta\to 0$. Consequently, $\|M\|_{\alpha^\prime}$ converges to 0 in probability as $\delta\to 0$.

On the other hand, Assumption \textbf{(HN)} implies that $\exp(C\|B\|^l_\alpha)<\infty$ almost surely for any $l\leq  1$.
Combining this fact with the stability estimate in Proposition \ref{pro-10}, we conclude that $\|X^\delta-\bar{X}\|_{\alpha^\prime}\to 0$ in probability as $\delta\to 0$.
By continuous mapping theorem, for every $1\leq p<\infty$, $\|X^\delta-\bar{X}\|^p_{\alpha^\prime}\to 0$ in probability as $\delta\to 0$. 
Moreover, the uniform moment estimates yield
\begin{align*}
\sup_{0<\delta\leq 1}\mathbb{E}\big[\|X^\delta-\bar{X}\|^{2p}_{\alpha^\prime}\big]
\leq 2^{2p-1}\sup_{0<\delta\leq 1}\mathbb{E}\big[\|X^\delta\|^{2p}_{\alpha^\prime}\big]+2^{2p-1} \mathbb{E}\big[\|\bar{X}\|^{2p}_{\alpha^\prime}\big]
\leq C\big(\mathbb{E}\big[\|B\|^{2lp}_\alpha\big]+1\big)<\infty.
\end{align*}
Thus, the family
$\{
\|X^\delta-\bar{X}\|^{p}_{\alpha^\prime}
\}_{0<\delta\leq 1}$ is uniformly integrable for every fixed $1\leq p<\infty$.
Finally, Vitali's convergence theorem yields
$$\mathbb{E}\big[\|X^\delta-\bar{X}\|^{p}_{\alpha^\prime}\big]\to 0 \quad \mathrm{as}~ \delta\to 0.$$
This completes the proof of Theorem \ref{th-2}.
\end{proof}

\section{Uniform Large Deviation Principle}\label{sec-5}
Throughout this section, we assume $H\in(1/2,1)$. Moreover, let $\alpha\in(1/2,H)$ and $\beta\in(1/3,1/2)$ satisfy $\alpha+\beta>1$.
We now turn to the uniform LDP for the slow--fast Young-rough system with a small noise parameter $\varepsilon$. More precisely, consider the following system:
\begin{align}\label{eq-45}
	\begin{cases}
		dX^{\varepsilon,\delta;x_0,y_0}_t=b_1(t,X^{\varepsilon,\delta;x_0,y_0}_t,Y^{\varepsilon,\delta;x_0,y_0}_t)dt+\sqrt\varepsilon\sigma_1(t,X^{\varepsilon,\delta;x_0,y_0}_t)dB^H_t,\\
		dY^{\varepsilon,\delta;x_0,y_0}_t=\frac1\delta b_2(X^{\varepsilon,\delta;x_0,y_0}_t,Y^{\varepsilon,\delta;x_0,y_0}_t)dt+
		\frac{1}{\sqrt{\delta}}\sigma_2(X^{\varepsilon,\delta;x_0,y_0}_t,Y^{\varepsilon,\delta;x_0,y_0}_t)d\mathbb{W}_t,\\
		(X^{\varepsilon,\delta;x_0,y_0}_0,Y^{\varepsilon,\delta;x_0,y_0}_0)=(x_0,y_0)\in\mathbb{R}^{d_1+d_2},~ t\in[0,T],~ 0<\delta<\varepsilon\leq 1,
	\end{cases}
\end{align}
Here, $(x_0,y_0)$ is a deterministic initial condition,
$B^H=(B^H_t)_{t\geq0}$ is an $m$-dimensional FBM with Hurst parameter $H$, and
$\mathbb{W}=(W^1,W^2)$ is the canonical It\^o-type Brownian rough path associated with an $r$-dimensional
standard BM $w=(w_t)_{t\geq0}$, independent of $B^H$, in the sense of Definition \ref{de-1}.
By Proposition \ref{pro-2}, the fast component $Y^{\varepsilon,\delta;x_0,y_0}$ satisfies the following It\^o SDE up to the explosion time of $(X^{\varepsilon,\delta;x_0,y_0},Y^{\varepsilon,\delta;x_0,y_0})$:
\begin{align*}
dY^{\varepsilon,\delta;x_0,y_0}_t=\frac1\delta b_2(X^{\varepsilon,\delta;x_0,y_0}_t,Y^{\varepsilon,\delta;x_0,y_0}_t)dt+
\frac{1}{\sqrt{\delta}}\sigma_2(X^{\varepsilon,\delta;x_0,y_0}_t,Y^{\varepsilon,\delta;x_0,y_0}_t)d^Iw_t,\quad Y^{\varepsilon,\delta;x_0,y_0}_0=y_0.
\end{align*}
For notational simplicity, we suppress the dependence on the initial condition $(x_0,y_0)$ whenever no confusion arises.

In this section, we establish an LDP for the slow component $X^{\varepsilon,\delta}$, uniformly over initial conditions $(x_0,y_0)$ belonging to compact subsets of $\mathbb{R}^{d_1+d_2}$.
More precisely, let $\mathscr{E}:= \big\{
E\subset\mathbb{R}^{d_1+d_2}:
E\text{ is compact}
\big\}$.
We show that the large deviation estimates hold uniformly over $(x_0,y_0)\in E$ for every $E\in\mathscr{E}$, as $\varepsilon\to0$ and $\delta/\varepsilon\to0$.

Throughout this section, we impose the following local regularity conditions on the slow--fast system \eqref{eq-45}:
\begin{itemize}
	\item For each fixed $t\in[0,T]$, $\sigma_1(t,\cdot)\in C^2(\mathbb{R}^{d_1},\mathcal{L}(\mathbb{R}^m,\mathbb{R}^{d_1}))$ and $\sigma_2\in C^3(\mathbb{R}^{d_1}\times\mathbb{R}^{d_2},\mathcal{L}(\mathbb{R}^r,\mathbb{R}^{d_2}))$. 
	Moreover, 
	for every $N>0$, there exists a constant $C_N>0$ such that
	\begin{align*}
		\sup_{|x|\leq N}|\nabla^i\sigma_1(t,x)|&\leq C_N,\\
		|\sigma_1(t,x)-\sigma_1(s,x)|+|\nabla\sigma_1(t,x)-\nabla\sigma_1(s,x)|&\leq C_N|t-s|^\beta, 
	\end{align*}
	for all $i=0,1,2$, $s,t\in[0,T]$ and $x\in\mathbb{R}^{d_1}$ with $|x|\leq N$.

	\item For each fixed $t\in[0,T]$, the functions $b_1(t,\cdot,\cdot)\in\mathcal{C}(\mathbb{R}^{d_1}\times\mathbb{R}^{d_2},\mathbb{R}^{d_1})$ and  $b_2\in\mathcal{C}(\mathbb{R}^{d_1}\times\mathbb{R}^{d_2},\mathbb{R}^{d_2})$ are locally Lipschitz continuous. 
	More precisely, for every $N>0$, there exists a constant $L_N>0$ such that
	$$
	|b_1(t,x_1,y_1)-b_1(t,x_2,y_2)|+|b_2(x_1,y_1)-b_2(x_2,y_2)|
	\leq L_N(|x_1-x_2|+|y_1-y_2|),
	$$
	for all $t\in[0,T] $ and $(x_i,y_i)\in\mathbb{R}^{d_1}\times\mathbb{R}^{d_2}$ with $|x_i|+|y_i|\leq N$, $i=1,2$.
\end{itemize}

By Remark \ref{rem-1}, these conditions guarantee that the slow–fast system \eqref{eq-45} always admits a unique time-local solution.

To formulate our main result, we introduce additional assumptions on the coefficients. Assumptions \textbf{(H1)}--\textbf{(H7)$_r$} remain the same as those used to establish the AP in Section \ref{sec-4}, except that \textbf{(H5)} is replaced by the stronger condition \textbf{(H5)$^\prime$} stated below.
In addition, we impose a further assumption \textbf{(H8)} to establish the uniform LDP.

\paragraph{\textbf{(H5)}$^\prime$} There exist constants $\eta_2\geq 0$ and $C>0$  such that, for all $(x,y)\in\mathbb{R}^{d_1+d_2}$,$$|b_2(x,y)|\leq C(|x|^{\eta_2}+|y|+1).$$

\paragraph{\textbf{(H8)}} 
There exists a constant $L_4>0$ such that, for all $x\in\mathbb{R}^{d_1}$,
$$\sup_{y\in\mathbb{R}^{d_2}}|\sigma_2(x,y)|\leq L_4(1+|x|).$$

Throughout this section, we take $\alpha^\prime$ and $\beta^\prime$ sufficiently close to $\alpha$ and $\beta$, respectively, so that $\frac13<\beta^\prime<\beta\leq \frac12<\alpha^\prime<\alpha<H$
and $\alpha^\prime+\beta^\prime>1$.

By Proposition \ref{pro-3}, under Assumptions \textbf{(H1)}--\textbf{(H3)}, \textbf{(H4)$_q$} and \textbf{(H5)$^\prime$} for some $2\leq q<\infty$, slow--fast system \eqref{eq-45} admits a unique global solution.
 \begin{proposition}\label{pro-9}
 Assume \textbf{(H1)}--\textbf{(H3)},  \textbf{(H4)$_q$} and \textbf{(H5)$^\prime$} for some $2\leq q<\infty$. Then, the probability that $Z^{\varepsilon,\delta}:=(X^{\varepsilon,\delta},Y^{\varepsilon,\delta})$ explodes on $[0,T]$ is zero. Moreover, for any compact subset $E$ of $\mathbb{R}^{d_1+d_2}$, we have 
 \begin{align}
 	&\sup_{0<\delta<\varepsilon\leq 1}\sup_{(x_0,y_0)\in E}\mathbb{E}\big[\|X^{\varepsilon,\delta}\|^p_{\alpha^\prime;[0,T]}\big]<\infty,\quad1\leq p<\infty,\label{eq-46}\\
 	&\sup_{0<\delta<\varepsilon\leq 1}\sup_{(x_0,y_0)\in E}\sup_{t\in[0,T]}\mathbb{E}\big[|Y^{\varepsilon,\delta}_t|^q\big]<\infty.\label{eq-47}
 \end{align}
 \end{proposition}
 \begin{proof}
 	Since $E$ is compact, the initial value $(x_0,y_0)\in E$ is bounded. By making minor modifications to the proof of Proposition \ref{pro-3}, and using the boundedness of $E$, we obtain the uniform estimates \eqref{eq-46} and \eqref{eq-47}. The same argument also shows that the corresponding solution cannot explode on $[0,T]$. Therefore, the probability of explosion is zero.
 \end{proof}
 
For any fixed $(x_0,y_0)\in \mathbb{R}^{d_1+d_2}$,
let
$$\mathcal{G}_{x_0,y_0}^{\varepsilon,\delta}(\sqrt\varepsilon\bullet,\sqrt\varepsilon\star):\mathcal{C}_0([0,T],\mathbb{R}^{m+r})\to \mathcal{C}^\alpha([0,T],\mathbb{R}^{d_1})$$
be a Borel-measurable version of $X^{\varepsilon,\delta;x_0,y_0}$.
In particular, $X^{\varepsilon,\delta;x_0,y_0}=\mathcal{G}_{x_0,y_0}^{\varepsilon,\delta}(\sqrt\varepsilon B^H,\sqrt\varepsilon w)$ almost surely. Here, $\mathcal{C}_0([0,T],\mathbb{R}^{m+r})$ denotes the space of $\mathbb{R}^{m+r}$-valued continuous paths starting from the origin $0$.

\subsection{Cameron-Martin space}

\setlength{\parindent}{2em} 
In this subsection, we introduce an $\mathbb{R}^{m+r}$-valued mixed FBM with Hurst parameter $H$ and recall some basic properties that will be used throughout the sequel. Let $(\Omega,\mathcal{F},\mathbb{P})$ be a suitable probability space equipped with a filtration satisfying the usual conditions.

Let $\beta^H=(\beta_t^H)_{t\in[0,T]}$ denote a one-dimensional FBM with Hurst parameter $H$. That is, $\beta^H$ is a centered Gaussian process with covariance function (cf. \cite{Mandelbrot1968}) given by
$$\mathbb{E}\left[\beta^H_t\beta^H_s\right]=\frac12\left(t^{2H}+s^{2H}-|t-s|^{2H}\right).$$
Let $\beta=(\beta_t)_{t\in[0,T]}$ be a standard scalar BM. 
Define $$c_H:=\Big[\frac{2H\Gamma\left(\frac32-H\right)\Gamma\left(H+\frac12\right)}{\Gamma(2-2H)}\Big]^{\frac12},$$
and denote by $F$ the Gauss hypergeometric function.
For all $0\leq s\leq t\leq T$, define
$$k_H(t,s):=\frac{c_H}{\Gamma\left(H+\frac12\right)}(t-s)^{H-\frac12}F\Big(H-\frac12,\frac12-H,H+\frac12;1-\frac ts\Big).$$
We can construct an FBM $\beta^H$ through the representation
\begin{align}\label{eq-48}
\beta^H_t:=\int_{0}^{t}k_H(t,s)d\beta_s.
\end{align}
Furthermore, we introduce the integral operator $\mathbb{K}_H$ associated with the kernel $k_H$ defined by
$$\mathbb{K}_Hf(t):=\int_{0}^{t}k_H(t,s)f(s)ds,$$
for any $f\in L^2([0,T],\mathbb{R})$.

The multidimensional counterparts are defined componentwise.
In particular, $B^H=(B_t^H)_{t\in[0,T]}$ is an $\mathbb{R}^m$-valued FBM with independent components, each having the same Hurst parameter $H$. It is well known that $B^H$ satisfies \textbf{(HN)} for every $\alpha\in(1/2,H)$.
Similarly, $w=(w_t)_{t\in[0,T]}$ is an $\mathbb{R}^r$-valued standard BM with independent components.

The Cameron-Martin Hilbert space $\mathcal{H}^H=\mathcal{H}^H([0,T],\mathbb{R}^m)$ associated with $B^H$ is defined by 
$$\mathcal{H}^H:=\Big\{
h=\mathbb{K}_H \dot{h}:\dot{h}\in L^2([0,T],\mathbb{R}^m)
\Big\}.
$$
Here, $\dot{h}$ denotes the corresponding element of $L^2([0,T],\mathbb{R}^m)$, and the dot does not indicate differentiation.
It is well known that 
$$\mathcal{H}^H\subset \mathcal{C}^{H'}([0,T],\mathbb{R}^m),
\quad H'\in(0,H).$$
The inner product on $\mathcal{H}^H$ is defined by
$\langle h,g\rangle_{\mathcal{H}^H}
:=\langle\dot{h},\dot{g}\rangle_{L^2([0,T],\mathbb{R}^m)}.$
Similarly, the Cameron--Martin Hilbert space $\mathcal{H}^\frac12=\mathcal{H}^\frac12([0,T],\mathbb{R}^r)$ associated with $w$  is given by 
$$\mathcal{H}^\frac12:=\Big\{v\in \mathcal{C}([0,T],\mathbb{R}^r):v=\int_0^\cdot v_s^\prime ds,~v^\prime \in L^2([0,T],\mathbb{R}^r)\Big\},$$
with inner product $\langle v,k\rangle_{\mathcal{H}^\frac12}:=\int_0^T\langle v_t',k_t'\rangle dt.$ 
Here, the prime denotes the usual time derivative.
Since $B^H$ and $w$ are independent, the $\mathbb{R}^{m+r}$-valued process $(B_t^H,w_t)_{t\in[0,T]}$ is referred to a mixed FBM with Hurst parameter $H$. Its Cameron-Martin space is given by $\mathcal{H}:=\mathcal{H}^H\oplus\mathcal{H}^\frac12$.

\subsection{Variational representation for mixed FBM}

For $N>0$, we set
$$\mathcal{S}_N=\Big\{(u,v)\in\mathcal{H}:\frac12\|(u,v)\|^2_\mathcal{H}:=\frac12(\|u\|^2_{\mathcal{H}^H}+\|v\|^2_{\mathcal{H}^\frac12})\leq N\Big\}.$$  
Equipped with the weak topology, the ball $\mathcal{S}_N$ is a compact Polish space (a complete separable metric space), since it is a bounded closed subset of the separable Hilbert space $\mathcal{H}$ and the weak topology is metrizable on bounded subsets.


Let $(\Omega,\mathcal{F},\mathbb{P})$ be the classical Wiener space, where
$
\Omega=\mathcal{C}_0([0,T],\mathbb{R}^{m+r})
$
is equipped with the uniform topology, $\mathbb{P}$ is the Wiener measure on $\Omega$, and $\mathcal{F}$ is the $\mathbb{P}$-completion of the Borel $\sigma$-field on $\Omega$. The coordinate process
$
(B_t,w_t)_{t\in[0,T]}
$
is an $\mathbb{R}^{m+r}$-valued standard BM under $\mathbb{P}$.
Let $\{\mathcal{F}_t\}_{t\in[0,T]}$ denote the natural augmented filtration generated by $(B,w)$, that is,
$
\mathcal{F}_t
:=
\sigma\{(B_s,w_s):0\leq s\leq t\}\vee\mathcal{N},
$
where $\mathcal{N}$ denotes the collection of all $\mathbb{P}$-null sets in $\mathcal{F}$.
Using the representation $\eqref{eq-48}$, we obtain a FBM $(B^H_t)_{t\in[0,T]}$ with the Hurst parameter $H$ defined on the same probability space $(\Omega,\mathcal{F},\mathbb{P})$.  Moreover, the augmented natural filtrations generated by $(B^H,w)$ and $(B,w)$ coincide.
Throughout this section, all stochastic processes and random variables are defined on the probability space $(\Omega,\mathcal{F},\mathbb{P})$, and we work with the filtration $\{\mathcal{F}_t\}_{t\in[0,T]}$.

For $N>0$, we denote by $\mathcal{B}_b^N$ the set of all $\mathbb{R}^{m+r}$-valued $\{\mathcal{F}_t\}$-progressively measurable processes
$(\varphi_t,\psi_t)_{t\in[0,T]}$ satisfying 
$\frac{1}{2}\int_{0}^{T}(|\varphi_s|^2+|\psi_s|^2)ds\leq N$ almost surly.
Set $\mathcal{B}_b:=\bigcup_{N>0}\mathcal{B}_b^N$.

For $N>0$, let $\mathcal{A}_b^N$ denote the set of all $\mathbb{R}^{m+r}$-valued $\{\mathcal{F}_t\}$-progressively measurable processes $(u_t,v_t)_{t\in[0,T]}$ such that
$(u,v)=(\mathbb{K}_H\dot{u},\int_0^\cdot v_s^\prime ds)$ almost surly, for some  $(\dot{u},v^\prime )\in\mathcal{B}_b^N$.
Set $\mathcal{A}_b:=\bigcup_{N>0}\mathcal{A}_b^N$. Every $(u,v)\in\mathcal{A}_b^N$ can be regarded as an $\mathcal{S}_N$-valued random variable. Since $\mathcal{S}_N$ is compact, the family $\{\mathbb{P}\circ(u,v)^{-1}:(u,v)\in\mathcal{A}_b^N\}$ is automatically tight. Moreover, by Girsanov's formula \cite[Theorem 10.18]{Prato2014}, the law of $(B^H_t+u_t,w_t+v_t)_{t\in[0,T]}$ is mutually absolutely continuous with respect to that of $(B^H_t,w_t)_{t\in[0,T]}$ for every $(u,v)\in\mathcal{A}_b^N$.

The following variational representation formula for $\mathbb{R}^{m+r}$-valued mixed FBM plays a key role in our analysis (see \cite{inahama2023large}). 
\begin{lemma}\label{lem2.3}
	Let $F:\mathcal{C}([0, T],\mathbb{R}^{m+r})\to \mathbb{R}$ be a bounded Borel-measurable function. Then
	$$-\log \mathbb{E}\left[\exp\{-F(B^H,w)\}\right]\\
	=\inf_{(u,v)\in\mathcal{A}_b}\mathbb{E}\left[F\left(B^H+u,w+v\right)+\frac12\|(u,v)\|_\mathcal{H}^2\right].$$
\end{lemma}

\subsection{Statement of uniform LDP}

For any $(u,v)\in\mathcal{H}$, the skeleton equation is defined by
\begin{align}\label{eq-49}
d\bar{X}^{u;x_0}_t=\bar{b}_1(t,\bar{X}^{u;x_0}_t)dt+\sigma_1(t,\bar{X}^{u;x_0}_t)du_t,\quad \bar{X}^{u;x_0}_0=x_0,\quad t\in[0,T],
\end{align}
where $\bar{b}_1(t,x)=\int_{\mathbb{R}^{d_2}}b_1(t,x,z)\mu^x(dz)$ with the invariant measure $\mu^x$ of the frozen equation \eqref{eq-37} as defined in Section \ref{sec-4}. Here, \(u\) denotes the first component of \((u,v)\).
By Proposition \ref{pro-5} and Remark \ref{rem-3}, under Assumptions \textbf{(H1)}, \textbf{(H2)}, \textbf{(H4)$_{2(r\vee 1)}$}, \textbf{(H6)} and \textbf{(H7)}$_{r}$ for some $r\geq 0$, skeleton equation \eqref{eq-49} admits a unique global solution. Moreover, for any $1\leq p <\infty$ and \(N>0\), there exist constants $C, l>0$ such that 
$$\sup_{(x_0,y_0)\in E}\sup_{(u,v)\in\mathcal{S}_N}\|\bar{X}^{u;x_0}\|^p_{\alpha^\prime}\leq \sup_{(u,v)\in\mathcal{S}_N}C(\|u\|^{lp}_\alpha+1)<\infty,$$
where $l$ and $C$ are independent of $(x_0,y_0)\in E$ and $(u,v)\in\mathcal{S}_N$.
For each $x_0\in\mathbb{R}^{d_1}$, define the
solution map
$
\mathcal{G}_{x_0}^{0}:
\mathcal{H}\to
\mathcal{C}^{\alpha}([0,T],\mathbb{R}^{d_1})
$
by
$
\mathcal{G}_{x_0}^{0}(u,v)
:=\bar{X}^{u;x_0}.
$
Note that $\mathcal{G}_{x_0}^{0}(u,v)$ depends only on $u$ and $x_0$, and is independent of $v$ and the initial condition $y_0$.

We now present the main result of this section, which establishes a Freidlin--Wentzell uniform LDP for \(X^{\varepsilon,\delta;x_0,y_0}\), uniformly with respect to the initial conditions \((x_0,y_0)\in E\), as \(\varepsilon\to0\) and \(\delta/\varepsilon\to0\).

Define the rate function by
\begin{align}\label{eq-65}
	I_{x_0,y_0}(\varphi):=\inf_{\{(u,v)\in\mathcal{H}:~\varphi=\mathcal{G}_{x_0}^0(u,v)\}}{\frac12\|(u,v)\|^2_{\mathcal{H}}}
	=\inf_{\{u\in\mathcal{H}^H:~\varphi=\mathcal{G}_{x_0}^0(u,0)\}}{\frac12\|u\|^2_{\mathcal{H}^H}},
\end{align}
where the infimum over the empty set is taken as $+\infty$.
Since the skeleton equation does not depend on the initial condition $y_0$ of the fast component, neither does the rate function. Hence, we simply write $I_{x_0}$ in place of $I_{x_0,y_0}$.
For each $M\geq 0$, define the level set of $I_{x_0}$ by
$$
\Phi_{x_0}(M):=\{\varphi\in\mathcal{C}^\alpha([0,T],\mathbb{R}^{d_1}):I_{x_0}(\varphi)\leq M\}.
$$
\begin{theorem}\label{th-3}
Let $H\in(1/2,1)$.
Assume \textbf{(H1)}--\textbf{(H3)}, \textbf{(H4)$_q$}, \textbf{(H5)$^\prime$}, \textbf{(H6)}, \textbf{(H7)$_{r}$} and \textbf{(H8)} for some $2\leq q<\infty$ and $0\leq r<\infty$ satisfying $q>2r$.
Let $\varepsilon\to 0$ and $\delta/\varepsilon\to 0$. Then, the slow components $X^{\varepsilon, \delta;x_0,y_0}$ of the slow--fast system \eqref{eq-45} satisfy the Freidlin–Wentzell uniform LDP in $\mathcal{C}^{\alpha^\prime}([0,T],\mathbb{R}^{d_1})$, for every $\alpha^\prime\in(1/2,H)$, uniformly over the class $\mathscr{E}$ of all compact subsets of $\mathbb{R}^{d_1+d_2}$,
with good rate functions $I_{x_0}$.
More precisely, 
	\begin{itemize}
	\item For any $E\in\mathscr{E}$, $M_0>0$, and $\zeta>0$,
	$$\liminf_{\substack{\varepsilon\to0\\\delta/\varepsilon\to 0}}\inf_{(x_0,y_0)\in E}
	\inf_{\varphi\in\Phi_{x_0}(M_0)}
	\big(\varepsilon\log\mathbb{P}\big(\|X^{\varepsilon,\delta;x_0,y_0}-\varphi\|_{\alpha^\prime}<\zeta\big)+I_{x_0}(\varphi)\big)\geq 0.$$
	\item For any $E\in\mathscr{E}$, $M_0>0$, and $\zeta>0$,
	$$\limsup_{\substack{\varepsilon\to0\\\delta/\varepsilon\to 0}}\sup_{(x_0,y_0)\in E}
	\sup_{M\in[0,M_0]}
	\big(\varepsilon\log\mathbb{P}\big(
	\inf_{\psi\in \Phi_{x_0}(M)} \|X^{\varepsilon,\delta;x_0,y_0}-\psi\|_{\alpha^\prime}\geq \zeta\big)
	+M\big)\leq 0.$$
\end{itemize}
\end{theorem}

\begin{remark}\label{rem-5}

To the author's knowledge, 
there are two preceding LDP
results \cite{gailus2025large,inahama2023large} on this kind of slow-fast systems driven by mixed FBM with $H >1/2$ with rather general coefficients.
Compared with them, our LDP is uniform with respect to the initial value and the topology of the path space is satisfactory, namely, the H\"older exponent is $H-\kappa$ for any sufficiently small $\kappa >0$.
The former \cite{gailus2025large} employs the viable pair approach within the weak convergence framework, while the latter \cite{inahama2023large} relies directly on the weak convergence method, leading to a simpler proof.
The present paper basically follows \cite{inahama2023large}, but more general.
For instance, the conditions on the coefficients are much weaker than those in \cite{inahama2023large}, which benefits from the Young-rough framework.

It is worth noting that, while the LP and the Freidlin--Wentzell LDP are equivalent under standard conditions, the uniform setting requires more careful treatment. In this context, Salins \cite{salins2019equivalences} introduced the equicontinuous uniform Laplace principle, a stronger version of the uniform Laplace principle, and established its equivalence with the Freidlin--Wentzell uniform LDP without additional assumptions. This equivalence plays an important role in our analysis (see Theorem \ref{th-6}), as the equicontinuous uniform Laplace principle can be conveniently verified through the weak convergence approach (see Theorem \ref{th-4}).

Unlike the case of single SDEs,
the uniform LDPs remain relatively unexplored  for slow--fast systems.
To the best of our knowledge, the only preceding work addressing a uniform LDP for slow–fast systems is \cite{li2026uniform},
in which a uniform LDP for 
 McKean-Vlasov type slow-fast system driven by usual BMs is proved. 
Compared with \cite{li2026uniform}, 
our result establishes a uniform LDP in a H\"older topology and extends the theory beyond the classical Brownian setting. It also provides a systematic approach to establishing uniform LDPs for slow–fast systems driven by FBM.



\end{remark}


\subsection{Preliminary estimates}

Throughout this subsection, unless otherwise stated, let $N>0$ and $E\in\mathcal{E}$ be arbitrary.
For notational simplicity, we suppress the dependence of the constants on $N$ and $E$ throughout this subsection.

To prove Theorem \ref{th-3}, we first derive several crucial estimates. Let $(u,v)\in\mathcal{A}_b^N$. We consider the following controlled system associated with \eqref{eq-45}.
\begin{align}\label{eq-50}
	\begin{cases}
		d\tilde{X}^{\varepsilon,\delta}_t=b_1(t,\tilde{X}^{\varepsilon,\delta}_t,\tilde{Y}^{\varepsilon,\delta}_t)dt+\sqrt\varepsilon\sigma_1(t,\tilde{X}^{\varepsilon,\delta}_t)dB^H_t+\sigma_1(t,\tilde{X}^{\varepsilon,\delta}_t)du_t,\\
		d\tilde{Y}^{\varepsilon,\delta}_t=\frac1\delta b_2(\tilde{X}^{\varepsilon,\delta}_t,\tilde{Y}^{\varepsilon,\delta}_t)dt+
		\frac{1}{\sqrt{\delta}}\sigma_2(\tilde{X}^{\varepsilon,\delta}_t,\tilde{Y}^{\varepsilon,\delta}_t)d^Iw_t+\frac{1}{\sqrt{\delta\varepsilon}}\sigma_2(\tilde{X}^{\varepsilon,\delta}_t,\tilde{Y}^{\varepsilon,\delta}_t)dv_t,\\
		(\tilde{X}^{\varepsilon,\delta}_0,\tilde{Y}^{\varepsilon,\delta}_0)=(x_0,y_0)\in\mathbb{R}^{d_1+d_2},~ t\in[0,T],~0<\delta<\varepsilon\leq 1.
	\end{cases}
\end{align}
By an argument analogous to that in the proof of Proposition \ref{pro-9}, under Assumptions \textbf{(H1)}--\textbf{(H3)}, \textbf{(H4)$_q$} and \textbf{(H5)$^\prime$} for some $2\leq q<\infty$, the controlled system \eqref{eq-50} admits a unique global solution $(\tilde{X}^{\varepsilon,\delta},\tilde{Y}^{\varepsilon,\delta})$.
Moreover, by Girsanov's theorem, for any $(u,v)\in\mathcal{A}^N_b$ and $(x_0,y_0)\in\mathbb{R}^{d_1+d_2}$, we have
\[
\tilde{X}^{\varepsilon,\delta;x_0,y_0}
=\mathcal{G}_{x_0,y_0}^{\varepsilon,\delta}
(\sqrt{\varepsilon}B^H+u,\sqrt{\varepsilon}w+v)
\quad\text{a.s.},
\]
where $\mathcal{G}_{x_0,y_0}^{\varepsilon,\delta}$ denotes the solution map associated with the slow component of \eqref{eq-45}.
\begin{proposition}\label{pro-7}
Assume \textbf{(H1)}--\textbf{(H3)}, \textbf{(H4)$_q$} and \textbf{(H5)$^\prime$} for some $2\leq q<\infty$. 
Then, for any $1\leq p<\infty$, we have 
$$\sup_{0<\delta<\varepsilon\leq 1}\sup_{(x_0,y_0)\in E}\sup_{(u,v)\in\mathcal{A}^N_b} \mathbb{E}[\|\tilde{X}^{\varepsilon,\delta}\|^p_{\alpha^\prime}]<\infty.$$
\end{proposition}
\begin{proof}
	
By an argument analogous to that in the proof of Lemma \ref{lem-3}, there exist constants $C,l>0$, independent of $\varepsilon$, $\delta$ and $(x_0,y_0)\in E$, such that
\begin{align*}
\|\tilde{X}^{\varepsilon,\delta}\|_{\alpha^\prime}\leq 
C\big(\|B^H\|_\alpha+\|u\|_\alpha+1\big)^l.
\end{align*}
By the fact that $(u,v)\in\mathcal{A}^N_b$ and $\|u\|_\alpha\leq C\|u\|_{\mathcal{H}^H}$, for any $1\leq p<\infty$,
$$\sup_{0<\delta<\varepsilon\leq 1}\sup_{(x_0,y_0)\in E}\sup_{(u,v)\in\mathcal{A}^N_b} \mathbb{E}[\|\tilde{X}^{\varepsilon,\delta}\|^p_{\alpha^\prime}]<\infty.$$
This completes the proof.
\end{proof}

\begin{proposition}\label{pro-8}
Assume \textbf{(H1)}--\textbf{(H3)}, \textbf{(H4)$_q$} and \textbf{(H5)$^\prime$} for some $1\leq q<\infty$.
Then, for any $2\leq p<\infty$ and any $0<\delta<\varepsilon\leq 1$, 
$$\sup_{(x_0,y_0)\in E}\sup_{(u,v)\in\mathcal{A}^N_b}\mathbb{E}\Big[\sup_{0\leq t\leq T}|\tilde{Y}^{\varepsilon,\delta}_t|^p
\Big]<\infty.$$
\end{proposition}
\begin{proof}
Throughout this proof, $0<\delta<\varepsilon\leq 1$ are fixed. The constant $C=C_{\varepsilon,\delta}>0$ may depend on $\varepsilon$ and $\delta$ and may change from line to line.
Define $(X,Y):=(\tilde{X}^{\varepsilon,\delta},\tilde{Y}^{\varepsilon,\delta})$, $Y^*_t:=\sup_{0\leq s\leq t}|\tilde{Y}^{\varepsilon,\delta}_s|$ to simplify the symbols.

For any $2\leq p<\infty$, by H\"older's inequality, Burkholder--Davis--Gundy's inequality, and the fact that $(u,v)\in\mathcal{A}^N_b$, we obtain
\begin{align*}
\mathbb{E}[|Y^*_t|^p]
\leq &
C|y_0|^p
+C\mathbb{E}\Big[
\Big(
\int_0^t |\sigma_2(X_s,Y_s)|^2ds
\Big)^{\frac p2}
\Big(
\int_0^t |\frac{dv_s}{ds}|^2ds
\Big)^{\frac p2}
\Big] \\
&+ C\mathbb{E}\Big[
\Big(
\int_0^t |b_2(X_s,Y_s)|ds
\Big)^p
\Big] 
+C\mathbb{E}\Big[
\sup_{0\leq s\leq t}\Big|
\int_0^s \sigma_2(X_r,Y_r)d^Iw_r
\Big|^p
\Big] 
\\
\leq &
C|y_0|^p
+C\mathbb{E}\Big[
\Big(
\int_0^t |\sigma_2(X_s,Y_s)|^2ds
\Big)^{\frac p2}
\Big] 
+ C\mathbb{E}\Big[
\Big(
\int_0^t |b_2(X_s,Y_s)|ds
\Big)^p
\Big] 
\\
\leq &
C|y_0|^p+C\mathbb{E}\Big[\int_0^t\big(1+|X_s|^{(\eta_2\vee 1) p}+|Y_s|^p\big)ds\Big]\\
\leq &
C|y_0|^p+C\int_0^t\mathbb{E}\Big[|Y^*_s|^p\Big]ds+C,
\end{align*}
where $\eta_2\geq0$ is the constant appearing in \textbf{(H5)$^\prime$}, and the last inequality follows from Tonelli's theorem and  Propositon \ref{pro-7}.
Since $E\subset\mathbb{R}^{d_1+d_2}$ is compact, $|y_0|$ is uniformly bounded for $(x_0,y_0)\in E$. Hence, by Gronwall's inequality,
$
\mathbb{E}[|Y^*_t|^p]
\leq C,
$
where constant $C>0$ is independent of $t\in[0,T]$, $(u,v)\in\mathcal{A}^N_b$ and $(x_0,y_0)\in E$.

For \(1\leq p<2\), Jensen's inequality and the estimate for \(p=2\) yield
\begin{align*}
\mathbb E[|Y_T^*|^p]
&=\mathbb E\big[(|Y_T^*|^2)^{p/2}\big]
\leq\left(\mathbb E[|Y_T^*|^2]\right)^{p/2}
\leq C.
\end{align*}
This completes the proof.
\end{proof}

\begin{lemma}\label{lem-12}
Assume \textbf{(H1)}--\textbf{(H3)}, \textbf{(H4)$_q$}, \textbf{(H5)$^\prime$} and \textbf{(H8)} for some $2\leq q<\infty$. Suppose that $ \delta/\varepsilon$ satisfies
\begin{align}\label{eq-66}
	\frac{\delta}{\varepsilon}
	\leq
	\frac{\gamma_1^2}{4},
	\end{align}
	where the constants \(\gamma_1>0\) is the constant appearing in  \textbf{(H4)$_q$}.
	Then, there exists a constant \(C>0\) such that, for any \((u,v)\in\mathcal{A}^N_b\), \((x_0,y_0)\in E\) and  \(0<\delta<\varepsilon\leq 1\) satisfying \eqref{eq-66},
	$$
	\int_0^T
	\mathbb{E}\big[|\tilde{Y}^{\varepsilon,\delta}_t|^2\big]\,dt
	\leq C.
	$$
\end{lemma}
\begin{proof}
For notational simplicity, we write $(X,Y):=(\tilde{X}^{\varepsilon,\delta},\tilde{Y}^{\varepsilon,\delta})$.
Applying It\^o's formula to $|Y_t|^2$, we obtain
\begin{align*}
|Y_t|^2
=&
|y_0|^2+\frac2\delta\int_0^t\big\langle Y_s,b_2(X_s,Y_s)\big\rangle ds+
\frac2{\sqrt\delta}\int_0^t\big\langle Y_s,\sigma_2(X_s,Y_s)d^Iw_s\big\rangle +
\frac2{\sqrt{\delta\varepsilon}}\int_0^t\Big\langle Y_s,\sigma_2(X_s,Y_s)\frac{dv_s}{ds}\Big\rangle ds\\
&+
\frac1\delta\int_0^t\big|\sigma_2(X_s,Y_s)\big|^2 ds.
\end{align*}
From Propositions \ref{pro-7} and \ref{pro-8}, we have $\mathbb{E}\big[\int_0^t\langle Y_s,\sigma_2(X_s,Y_s)d^Iw_s\rangle \big]=0$.
Then, from \textbf{(H4)}$_2$ and \textbf{(H8)}, by taking expectations, we have
\begin{align}\label{eq-63}
\frac{d\mathbb{E}[|Y_t|^2]}{dt}
\leq &
\frac2\delta\mathbb{E}\big[\langle Y_t,b_2(X_t,Y_t)\rangle\big]
+\frac1\delta\mathbb{E}\big[|\sigma_2(X_t,Y_t)|^2\big]
+\frac{2}{\sqrt{\delta\varepsilon}}\mathbb{E}\Big[ \Big\langle Y_t, \sigma_2(X_t,Y_t)\frac{dv_t}{dt}\Big\rangle\Big]\notag\\
\leq &
\frac1\delta\mathbb{E}\big[-\gamma_1|Y_t|^2+C|X_t|^{\eta_1}+C\big]
+\frac{C}{\sqrt{\delta\varepsilon}}\mathbb{E}\Big[(1+|X_t|^2)\Big|\frac{dv_t}{dt}\Big|^2\Big]+\frac{1}{\sqrt{\delta\varepsilon}}\mathbb{E}\big[|Y_t|^2\big]\notag\\
\leq &
\frac{-\gamma_1}{2\delta}\mathbb{E}\big[|Y_t|^2\big]+\frac C\delta
+\frac{C}{\sqrt{\delta\varepsilon}}\mathbb{E}\Big[\|X\|^2_{\alpha^\prime}\Big|\frac{dv_t}{dt}\Big|^2\Big]
+\frac{C}{\sqrt{\delta\varepsilon}}\mathbb{E}\Big[\Big|\frac{dv_t}{dt}\Big|^2\Big]
\end{align}
for almost every $t\in[0,T]$, where the second inequality follows from Young's inequality, while the last inequality follows from condition \eqref{eq-66} and Proposition \ref{pro-7}.
Consider the following ordinary differential equation 
\begin{align*}
\frac{dA_t}{dt}=\frac{-\gamma_1}{2\delta}A_t+\frac C{\sqrt{\delta\varepsilon}}\mathbb{E}\Big[\|X\|^2_{\alpha^\prime}\Big|\frac{dv_t}{dt}\Big|^2\Big]
+\frac{C}{\sqrt{\delta\varepsilon}}\mathbb{E}\Big[\Big|\frac{dv_t}{dt}\Big|^2\Big]+\frac C\delta,\quad
A_0=|y_0|^2.
\end{align*}
Moreover, we have
\begin{align*}
A_t=&
|y_0|^2e^{-\frac{\gamma_1}{2\delta}t}+\frac C{\sqrt{\delta\varepsilon}}\int_0^t e^{-\frac{ \gamma_1}{2\delta}(t-s)}\mathbb{E}\Big[\|X\|^2_{\alpha^\prime}\Big|\frac{dv_s}{ds}\Big|^2\Big]ds
\\
&+\frac C{\sqrt{\delta\varepsilon}}\int_0^t e^{-\frac{\gamma_1}{2\delta}(t-s)}\mathbb{E}\Big[\Big|\frac{dv_s}{ds}\Big|^2\Big]ds+\frac C\delta\int_0^t e^{-\frac{ \gamma_1}{2\delta}(t-s)}ds.
\end{align*}
Since $t\mapsto\mathbb{E}[|Y_t|^2]$ is absolutely continuous, the comparison principle applied to \eqref{eq-63} yields
\begin{align*}
\mathbb{E}[|Y_t|^2]\leq A_t
\end{align*}
for all $t\in[0,T]$.
Therefore, by Tonelli's theorem, we obtain
\begin{align*}
\int_0^T \mathbb{E}[|Y_t|^2]dt\leq &
|y_0|^2\int_0^T e^{-\frac{\gamma_1}{2\delta}t}dt
+\frac C{\sqrt{\delta\varepsilon}}\int_0^T\int_0^t e^{-\frac{\gamma_1}{2\delta}(t-s)}\mathbb{E}\Big[\|X\|^2_{\alpha^\prime}\Big|\frac{dv_s}{ds}\Big|^2\Big]dsdt\\&+
\frac C{\sqrt{\delta\varepsilon}}\int_0^T\int_0^t e^{-\frac{\gamma_1}{2\delta}(t-s)}\mathbb{E}\Big[\Big|\frac{dv_s}{ds}\Big|^2\Big]dsdt+\frac C\delta\int_0^T\int_0^t e^{-\frac{\gamma_1}{2\delta}(t-s)}dsdt\\
\leq &
|y_0|^2\frac{2\delta}{\gamma_1}+\frac C{\sqrt{\delta\varepsilon}}\mathbb{E}\Big[\|X\|^2_{\alpha^\prime}\int_0^T\int_s^Te^{-\frac{\gamma_1}{2\delta}(t-s)}dt\Big|\frac{dv_s}{ds}\Big|^2ds\Big]\\
&+\frac C{\sqrt{\delta\varepsilon}}\int_0^T\int_s^Te^{-\frac{\gamma_1}{2\delta}(t-s)}dt\mathbb{E}\Big[\Big|\frac{dv_s}{ds}\Big|^2\Big]ds+C\\
\leq &
C|y_0|^2+C\sqrt{\frac\delta\varepsilon}\mathbb{E}\Big[\|X\|^2_{\alpha^\prime}\int_0^T\Big|\frac{dv_s}{ds}\Big|^2ds\Big]+C\sqrt{\frac\delta\varepsilon}\int_0^T\mathbb{E}\Big[\Big|\frac{dv_s}{ds}\Big|^2\Big]ds+C\\
\leq &
C|y_0|^2+C \sqrt{\frac\delta\varepsilon}\mathbb{E}\big[\|X\|^2_{\alpha^\prime} \big]+C \sqrt{\frac\delta\varepsilon}+C\\
\leq &
C|y_0|^2+C\sqrt{\frac\delta\varepsilon}+C,
\end{align*}
where the last two inequalities follow from the fact that $(u,v)\in\mathcal{A}^N_b$ and Proposition~\ref{pro-7}, respectively.

Since $E$ is compact, we obtain
$$
\int_0^T
\mathbb{E}\big[|Y_t|^2\big]dt
\leq C,
$$
where $C>0$ is independent of $\varepsilon$, $\delta$ satisfying \eqref{eq-66}, $(u,v)\in\mathcal{A}_b^N$, and $(x_0,y_0)\in E$.
This completes the proof.
\end{proof}

We consider the following time-discretized auxiliary system associated with  $\eqref{eq-50}$:
\begin{align}\label{eq-51}
	\begin{cases}
		d\hat{X}^{\varepsilon,\delta}_t=\bar{b}_1(t,\hat{X}^{\varepsilon,\delta}_t)dt+\sqrt\varepsilon\sigma_1(t,\hat{X}^{\varepsilon,\delta}_t)dB^H_t+\sigma_1(t,\hat{X}^{\varepsilon,\delta}_t)du_t,\\
		d\hat{Y}^{\varepsilon,\delta}_t=\frac1\delta b_2(\tilde{X}^{\varepsilon,\delta}_{t(\Delta)},\hat{Y}^{\varepsilon,\delta}_t)dt+
		\frac{1}{\sqrt{\delta}}\sigma_2(\tilde{X}^{\varepsilon,\delta}_{t(\Delta)},\tilde{Y}^{\varepsilon,\delta}_t)d^Iw_t,\\
		(\hat{X}^{\varepsilon,\delta}_0,\hat{Y}^{\varepsilon,\delta}_0)=(x_0,y_0)\in\mathbb{R}^{d_1+d_2},~ t\in[0,T],~0<\delta<\varepsilon\leq 1,
	\end{cases}
\end{align}
where $t(\Delta):=\lfloor\frac t\Delta\rfloor\Delta$ is the nearest breakpoint preceding $t$ with $t\in[0,T]$. Without loss of generality, we assume that $0<\Delta<1\wedge T$.

Since $\hat{X}^{\varepsilon,\delta}$ does not depend on the fast component and $\hat{Y}^{\varepsilon,\delta}$ contains no control term, from Propositions \ref{pro-4} and \ref{pro-9}, we can obtain the boundedness of both $\hat{X}^{\varepsilon,\delta}$ and $\hat{Y}^{\varepsilon,\delta}$.
\begin{lemma}\label{lem-8}
 Assume \textbf{(H1)}--\textbf{(H3)}, \textbf{(H4)$_q$} and \textbf{(H5)$^\prime$} for some $2\leq q<\infty$. Then, we have 
\begin{align}
	&\sup_{\substack{0<\delta<\varepsilon\leq 1\\0<\Delta<1\wedge T}}\sup_{(x_0,y_0)\in E}
	\sup_{(u,v)\in\mathcal{A}^N_b}
	\mathbb{E}\big[\|\hat{X}^{\varepsilon,\delta}\|^p_{\alpha^\prime;[0,T]}\big]<\infty,\quad1\leq p<\infty,\label{eq-52}\\
	&\sup_{\substack{0<\delta<\varepsilon\leq 1\\0<\Delta<1\wedge T}}\sup_{0<\delta<\varepsilon\leq 1}\sup_{(x_0,y_0)\in E}
	\sup_{(u,v)\in\mathcal{A}^N_b}
	\sup_{t\in[0,T]}\mathbb{E}\big[|\hat{Y}^{\varepsilon,\delta}_t|^q\big]<\infty.\label{eq-53}
\end{align}
\end{lemma}

\begin{lemma}\label{lem-13}
Assume \textbf{(H1)}--\textbf{(H3)}, \textbf{(H4)}$_q$, \textbf{(H5)$^\prime$}, \textbf{(H6)}, \textbf{(H7)}$_{r}$ and \textbf{(H8)} for some $2\leq q<\infty$ and $0\leq r<\infty$ satisfying $q>2r$.
Then, there exists a constant $C>0$ such that
$$\sup_{(x_0,y_0)\in E}
\sup_{(u,v)\in\mathcal{A}^N_b}
\int_0^T \mathbb{E}\big[|\tilde{Y}^{\varepsilon,\delta}_t-\hat{Y}^{\varepsilon,\delta}_t|^2\big]dt\leq C\Big(\Delta+\frac\delta\varepsilon\Big),$$
where the constant $C$ is independent of $\varepsilon$, $\delta$ and $\Delta$.
\end{lemma}
\begin{proof}
For notational simplicity, we write $J:=\tilde{Y}^{\varepsilon,\delta}-\hat{Y}^{\varepsilon,\delta}$, $(X,Y):=(\tilde{X}^{\varepsilon,\delta}, \tilde{Y}^{\varepsilon,\delta})$.
By It\^o's formula, 
\begin{align*}
\mathbb{E}[|J_t|^2]
=&
\frac2\delta\mathbb{E}\Big[
\int_0^t\big\langle J_s,b_2(X_s,Y_s)-b_2(X_{s(\Delta)},\hat{Y}^{\varepsilon,\delta}_s)\big\rangle ds
\Big]
+\frac1\delta\mathbb{E}\Big[
\int_0^t\big|\sigma_2(X_s,Y_s)-\sigma_2(X_{s(\Delta)},\hat{Y}^{\varepsilon,\delta}_s)\big|^2 ds
\Big]\\
&+
\frac2{\sqrt{\delta}}\mathbb{E}\Big[
\int_0^t\big\langle J_s,\big\{\sigma_2(X_s,Y_s)-\sigma_2(X_{s(\Delta)},\hat{Y}^{\varepsilon,\delta}_s)\big\}d^Iw_s\big\rangle 
\Big]\\
&+
\frac2{\sqrt{\delta\varepsilon}}\mathbb{E}\Big[
\int_0^t\Big\langle J_s,\sigma_2(X_s,Y_s)\frac{dv_s}{ds}\Big\rangle ds
\Big].
\end{align*}
By the boundedness estimates for $X$, $Y$ and $\hat{Y}^{\varepsilon,\delta}$, the It\^o integral term is a true martingal in $L^2$. Hence, its expectation vanishes.
Taking the derivative with respect to $t$, we have 
\begin{align*}
\frac{d\mathbb{E}[|J_t|^2]}{dt}
=&
\frac1\delta \mathbb{E}\big[
2\big\langle J_t,b_2(X_t,Y_t)-b_2(X_t,\hat{Y}^{\varepsilon,\delta}_t)\big\rangle 
+\big|\sigma_2(X_t,Y_t)-\sigma_2(X_t,\hat{Y}^{\varepsilon,\delta}_t)\big|^2
\big]\\
&+
\frac2\delta\mathbb{E}\big[
\big\langle J_t,b_2(X_t,\hat{Y}^{\varepsilon,\delta}_t)-b_2(X_{t(\Delta)},\hat{Y}^{\varepsilon,\delta}_t)\big\rangle 
\big]\\
&+
\frac2\delta\mathbb{E}\big[
\big\langle \sigma_2(X_t,Y_t)-\sigma_2(X_t,\hat{Y}^{\varepsilon,\delta}_t),\sigma_2(X_t,\hat{Y}^{\varepsilon,\delta}_t)-\sigma_2(X_{t(\Delta)},\hat{Y}^{\varepsilon,\delta}_t)\big\rangle 
\big]\\
&+
\frac1\delta \mathbb{E}\big[
\big|\sigma_2(X_t,\hat{Y}^{\varepsilon,\delta}_t)-\sigma_2(X_{t(\Delta)},\hat{Y}^{\varepsilon,\delta}_t)\big|^2
\big]\\
&+
\frac2{\sqrt{\delta\varepsilon}}\mathbb{E}\Big[
\Big\langle J_t,\sigma_2(X_t,Y_t)\frac{dv_t}{dt}\Big\rangle 
\Big]\\
=: &I_1+I_2+I_3+I_4+I_5.
\end{align*}
By \textbf{(H6)}, we have 
\begin{align*}
I_1\leq \frac{-\gamma_2}\delta\mathbb{E}[|J_t|^2].
\end{align*}
Next, by \textbf{(H7)$_r$}, Young's inequality and H\"older's inequality, for any $\xi>0$, we have
\begin{align*}
	I_2\leq & \frac C\delta\mathbb{E}\big[|J_t||X_t-X_{t(\Delta)}|\big(1+|X_t|^{\eta_3}+|X_{t(\Delta)}|^{\eta_3}+|\hat{Y}^{\varepsilon,\delta}_t|^r\big)\big]\\
	\leq &
	\frac{\gamma_2}{4\delta}\mathbb{E}[|J_t|^2]+\frac C\delta\mathbb{E}\big[|X_t-X_{t(\Delta)}|^2\big(1+|X_t|^{2\eta_3}+|X_{t(\Delta)}|^{2\eta_3}+|\hat{Y}^{\varepsilon,\delta}_t|^{2r}\big)\big]\\
	\leq &
	\frac{\gamma_2}{4\delta}\mathbb{E}[|J_t|^2]+\frac C\delta\mathbb{E}\big[|X_t-X_{t(\Delta)}|^{\frac{2(1+\xi)}{\xi}}\big]^{\frac{\xi}{1+\xi}}
	\mathbb{E}\big[1+|X_t|^{2\eta_3(1+\xi)}+|X_{t(\Delta)}|^{2\eta_3(1+\xi)}+|\hat{Y}^{\varepsilon,\delta}_t|^{2r(1+\xi)}
	\big]^{\frac1{1+\xi}}.
\end{align*}
Since $q>2r$, we may choose $\xi>0$ sufficiently small such that $2r(1+\xi)\leq q$.
Hence, by H\"older's inequality, Proposition \ref{pro-7} and Lemma \ref{lem-8}, we obtain
\begin{align*}
I_2
\leq{}&
\frac{\gamma_2}{4\delta}\mathbb{E}\big[|J_t|^2\big]
+
C \frac{\Delta^{2\alpha'}}{\delta}.
\end{align*}
Since $\sigma_2$ is globally Lipschitz continuous, Young's inequality yields
	\begin{align*}
		I_3+I_4\leq &
		\frac C\delta\mathbb{E}\big[|J_t||X_t-X_{t(\Delta)}|+|X_t-X_{t(\Delta)}|^2\big]\\
		\leq &
		\frac{\gamma_2}{4\delta}\mathbb{E}[|J_t|^2]+\frac C\delta\mathbb{E}\big[|X_t-X_{t(\Delta)}|^2\big]\\
		\leq &
		\frac{\gamma_2}{4\delta}\mathbb{E}[|J_t|^2]+C\frac{\Delta^{2\alpha^\prime}}\delta.
	\end{align*}
Finally, by Young's inequality and \textbf{(H8)},
\begin{align*}
I_5\leq \frac{\gamma_2}{4\delta}\mathbb{E}[|J_t|^2]+\frac C\varepsilon\mathbb{E}\Big[\big(1+|X_t|^2\big)\Big|\frac{dv_t}{dt}\Big|^2\Big].
\end{align*}
Therefore,
\begin{align*}
\frac{d\mathbb{E}[|J_t|^2]}{dt}
\leq &
\frac{-\gamma_2}{4\delta}\mathbb{E}[|J_t|^2]+C\frac{\Delta^{2\alpha^\prime}}\delta+\frac C\varepsilon\mathbb{E}\Big[\big(1+|X_t|^2\big)\Big|\frac{dv_t}{dt}\Big|^2\Big].
\end{align*}
Since $J_0=0$, Gronwall's inequality yields
\begin{align*}
	\mathbb{E}[|J_t|^2]
	\leq&
	C\frac{\Delta^{2\alpha^\prime}}{\delta}\int_0^t
	e^{-\frac{\gamma_2}{4\delta}(t-s)}
	ds
	+\frac{C}{\varepsilon}\int_0^t
	e^{-\frac{\gamma_2}{4\delta}(t-s)}
	\mathbb{E}\Big[
	\big(1+|X_s|^2\big)
	\Big|\frac{dv_s}{ds}\Big|^2
	\Big]
	ds.
\end{align*}
Integrating over $[0,T]$ and recalling that $\alpha^\prime>1/2$, so that
$\Delta^{2\alpha^\prime}<\Delta$ for $0<\Delta<1$, an argument analogous to that in the proof of Lemma \ref{lem-12} yields
$$
\sup_{(x_0,y_0)\in E}
\sup_{(u,v)\in\mathcal{A}^N_b}
\int_0^T
\mathbb{E}\big[
|J_t|^2
\big]dt
\leq 
C\Big(\Delta^{2\alpha^\prime}+\frac\delta\varepsilon\Big)
\leq
C\Big(
\Delta+\frac{\delta}{\varepsilon}
\Big)
$$
where we have used Proposition \ref{pro-7} and the fact that $(u,v)\in\mathcal{A}^N_b$ to obtain the first inequality.
This completes the proof.
\end{proof}

\subsection{Proof of uniform LDP}

This section is devoted to proving Theorem \ref{th-3}, which establishes the uniform LDP for $X^{\varepsilon,\delta;x_0,y_0}$ uniformly over $\mathscr{E}$, with the good rate function $I_{x_0}$. Throughout this section, we work under the assumptions of Theorem \ref{th-3}.

We first establish the continuity of the solution map
\[
\mathcal{G}^0:\mathbb{R}^{d_1}\times\mathcal{H}
\to\mathcal{C}^\alpha([0,T],\mathbb{R}^{d_1}),
\]
which will be used to show that $I_{x_0}$ is a good rate function. The following lemma is purely deterministic.

\begin{lemma}\label{lem-9}
Let $N>0$ and $E\in\mathscr{E}$. Suppose that
$\{(u^n,v^n)\}_{n\geq 1}\subset\mathcal{S}_N$ converges to $(u,v)\in\mathcal{S}_N$ with respect to the weak topology of $\mathcal{H}$, and that
$\{(x_0^n,y_0^n)\}_{n\geq 1}\subset E$ converges to $(x_0,y_0)\in E$ in $\mathbb{R}^{d_1+d_2}$. Then,
$$
\mathcal{G}_{x_0^n}^0(u^n,v^n)
\longrightarrow
\mathcal{G}_{x_0}^0(u,v)
\quad\text{in }\mathcal{C}^{\alpha^\prime}
([0,T],\mathbb{R}^{d_1})
\quad\text{as }n\to\infty.
$$
\end{lemma}
\begin{proof}

For each $n\geq1$, let $\bar{X}^{u^n;x^n_0}$ denote the solution to the skeleton equation
\begin{align*}
	d\bar{X}^{u^n;x^n_0}_t=\bar{b}_1(t,\bar{X}^{u^n;x^n_0}_t)dt+\sigma_1(t,\bar{X}^{u^n;x^n_0}_t)du^n_t,\quad \bar{X}^{u^n;x^n_0}_0=x^n_0,\quad t\in[0,T].
\end{align*}
Throughout the proof, for simplicity, we write $\bar{X}^{u^n;x^n_0}=:\bar{X}^{(n)}$.
Since $\mathcal{H}^H$ is continuously embedded into $\mathcal{C}^\alpha([0,T],\mathbb{R}^m)$ and $(u^n,v^n)\in\mathcal{S}_N$, the sequence $\{u^n\}_{n\geq1}$ is uniformly bounded in $\mathcal{C}^\alpha$.
By Propositions \ref{pro-4} and \ref{pro-5}, 
for any $\alpha\in(\frac12,H)$ and $1\leq p<\infty$, the skeleton equation admits a
unique global solution $\bar{X}^{(n)}\in\mathcal{C}^\alpha([0,T],\mathbb{R}^{d_1})$. Moreover, choosing $\iota>0$ sufficiently small such that $\frac12<\alpha^\prime<\alpha-\iota<H$, there exists a constant $C>0$, independent of $n$, such that $\|\bar{X}^{(n)}\|_{\alpha-\iota}\leq C$.
Here we have used the fact that $\{(x_0^n,y_0^n)\}_{n\geq1}\subset E$ is bounded.
By the compact embedding $\mathcal{C}^{\alpha-\iota}([0,T],\mathbb{R}^{d_1})\Subset \mathcal{C}^{\alpha^\prime}([0,T],\mathbb{R}^{d_1})$, the sequence $\{\bar{X}^{(n)}\}_{n\geq 1}$ is precompact in $\mathcal{C}^{\alpha^\prime}([0,T],\mathbb{R}^{d_1})$.
Hence, there exist a subsequence, still denoted by $\{\bar{X}^{(n)}\}_{n\geq1}$, and a function $\bar{X}\in\mathcal{C}^{\alpha^\prime}([0,T],\mathbb{R}^{d_1})$ such that
$$
\bar{X}^{(n)}\longrightarrow\bar{X}
\quad\text{in }
C^{\alpha^\prime}([0,T],\mathbb{R}^{d_1})\quad\text{as }n\to\infty.
$$

We next show that $\bar{X}$ satisfies the YDE
\begin{align*}
	d\bar{X}_t=\bar{b}_1(t,\bar{X}_t)dt+\sigma_1(t,\bar{X}_t)du_t,\quad \bar{X}_0=x_0,\quad t\in[0,T].
\end{align*}
In the following, $t\in[0,T]$ is arbitary but fixed, unless otherwise stated.
First, by Propositions \ref{pro-4} and \ref{pro-5}, we have
\begin{align*}
	\Big|
	\int_0^t\big\{
	\bar{b}_1(s,\bar{X}^{(n)}_s)-\bar{b}_1(s,\bar{X}_s)
	\big\}ds
	\Big|
	\leq &
	C\int_0^t |\bar{X}^{(n)}_s-\bar{X}_s|\big(1+|\bar{X}^{(n)}_s|^\xi+|\bar{X}_s|^\xi\big)ds\\
	\leq &
	C\sup_{t\in[0,T]}|\bar{X}^{(n)}_t-\bar{X}_t|.
\end{align*}
From the estimate of Young integral \eqref{eq-57} and the fact that $(u^n,v^n)\in\mathcal{S}_N$, we have 
\begin{align*}
	\Big|
	\int_0^t\big\{
	\sigma_1(s,\bar{X}^{(n)}_s)-\sigma_1(s,\bar{X}_s)
	\big\}du^n_s
	\Big|
	\leq &
	\Big\|\int_0^\cdot\big\{
	\sigma_1(s,\bar{X}^{(n)}_s)-\sigma_1(s,\bar{X}_s)
	\big\}du^n_s\Big\|_{\alpha-\iota}\\
	\leq &
	C\big\|	\sigma_1(\cdot,\bar{X}^{(n)}_\cdot)-\sigma_1(\cdot,\bar{X}_\cdot)\big\|_\beta\|u^n\|_\alpha\\
	\leq &
	C\big\|	\sigma_1(\cdot,\bar{X}^{(n)}_\cdot)-\sigma_1(\cdot,\bar{X}_\cdot)\big\|_\beta.
\end{align*}
By 
\textbf{(H1)}, 
for any $0\leq s<t\leq T$, we obtain 
\begin{align*}
&\big|\sigma_1(t,\bar{X}^{(n)}_t)-\sigma_1(t,\bar{X}_t)-\sigma_1(s,\bar{X}^{(n)}_s)+\sigma_1(s,\bar{X}_s)\big|\\
&=
\Big|\int_0^1 \nabla\sigma_1\big(t,\theta\bar{X}^{(n)}_t+(1-\theta)\bar{X}_t\big)\big\langle\bar{X}^{(n)}_t-\bar{X}_t\big\rangle
d\theta
-\int_0^1 \nabla\sigma_1\big(s,\theta\bar{X}^{(n)}_s+(1-\theta)\bar{X}_s\big)\big\langle\bar{X}^{(n)}_s-\bar{X}_s\big\rangle
d\theta\Big|\\
&=
\Big|\int_0^1  \big\{\nabla\sigma_1\big(t,\theta\bar{X}^{(n)}_t+(1-\theta)\bar{X}_t\big)-\nabla\sigma_1\big(s,\theta\bar{X}^{(n)}_s+(1-\theta)\bar{X}_s\big)\big\}\big\langle\bar{X}^{(n)}_t-\bar{X}_t\big\rangle d\theta\\
&\qquad+
\int_0^1 \nabla\sigma_1\big(s,\theta\bar{X}^{(n)}_s+(1-\theta)\bar{X}_s\big)\big\langle\bar{X}^{(n)}_t-\bar{X}_t-\bar{X}^{(n)}_s+\bar{X}_s\big\rangle
d\theta\Big|\\
&\leq C\big(
|t-s|^\beta+\|\bar{X}^{(n)}\|_{\alpha^\prime}|t-s|^{\alpha^\prime}+\|\bar{X}\|_{\alpha^\prime}|t-s|^{\alpha^\prime}
\big) \|\bar{X}^{(n)}-\bar{X}\|_\infty
+
C\|\bar{X}^{(n)}-\bar{X}\|_\beta|t-s|^\beta
\\
&\leq 
C\big(1+\|\bar{X}^{(n)}\|_{\alpha^\prime}+\|\bar{X}\|_{\alpha^\prime}\big)\|\bar{X}^{(n)}-\bar{X}\|_\beta|t-s|^\beta\\
&\leq C\|\bar{X}^{(n)}-\bar{X}\|_\beta|t-s|^\beta,
\end{align*}
where the last inequality follows from the boundedness of $\bar{X}^{(n)}$ and $\bar{X}$ in $\mathcal{C}^{\alpha^\prime}$.
Hence, we have
\begin{align*}
\Big|
\int_0^t\big\{
\sigma_1(s,\bar{X}^{(n)}_s)-\sigma_1(s,\bar{X}_s)
\big\}du^n_s
\Big|
\leq 
C\|\bar{X}^{(n)}-\bar{X}\|_{\alpha^\prime}.
\end{align*}
It remains to show that
$$
\Big|\int_0^t\sigma_1(s,\bar{X}_s)\,d(u^n_s-u_s)\Big|
\longrightarrow0.
$$
Define $\Lambda_t(h):=\int_0^t\sigma_1(s,\bar{X}_s)dh_s$ for $h\in\mathcal{H}^H$.
Then,
$\Lambda_t:\mathcal{H}^H\to \mathbb{R}^{d_1}$ is a continuous linear operator.
Writing
\begin{align*}
	\Lambda_t(h)=:\big(
	\Lambda^1_t(h),\cdots, \Lambda^{d_1}_t(h)
	\big),
\end{align*}
where $\Lambda^i_t(h)\in \mathbb{R}$, $i=1,\cdots,d_1$.
Then, each $\Lambda^i_t:\mathcal{H}^H\to \mathbb{R}$ is a  continuous linear functional.
By Riesz representation theorem, there exists a unique element $k^i_t\in\mathcal{H}^H$ such that $\Lambda^i_t(h)=\langle h,k^i_t\rangle_{\mathcal{H}^H}$ for any $h\in\mathcal{H}^H$.
Then, we have 
$$\Lambda_t(h)=\big(
\langle h,k^1_t\rangle_{\mathcal{H}^H},\cdots,\langle h,k^{d_1}_t\rangle_{\mathcal{H}^H}
\big)$$
Since $u^n\to u$ weakly in $\mathcal{H}^H$ as $n\to\infty$, it follows that
\begin{align*}
	\Big|
	\int_0^t
	\sigma_1(s,\bar{X}_s)d(u^n_s-u_s)
	\Big|
	=&
	|\Lambda_t(u^n)-\Lambda_t(u)|\\
	= &
	\big|
	\big(
	\langle u^n-u,k^1_t\rangle_{\mathcal{H}^H},\cdots,\langle u^n-u,k^{d_1}_t\rangle_{\mathcal{H}^H}
	\big)
	\big| \to 0\quad\text{as }n\to \infty.
\end{align*}

Combining the above convergence results with
$x_0^n\to x_0$, we conclude that $\bar X$ satisfies
\begin{align*}
	\bar{X}_t
	=x_0
	+
	\int_0^t\bar{b}_1(s,\bar{X}_s)ds
	+
	\int_0^t\sigma_1(s,\bar{X}_s)du_s.
\end{align*}
By the uniqueness of the solution to the skeleton equation, we obtain
$
\bar{X}
=
\mathcal{G}_{x_0}^{0}(u,v).
$
Since every convergent subsequence has the same limit, we obtain
$$
\mathcal{G}_{x_0^n}^{0}(u^n,v^n)
\longrightarrow
\mathcal{G}_{x_0}^{0}(u,v)
$$
in $\mathcal{C}^{\alpha^\prime}([0,T],\mathbb{R}^{d_1})$.
This completes the proof.
\end{proof}

\begin{corollary}
Let $N>0$ and $E\in\mathscr{E}$.
Let $\{(u^n,v^n)\}_{n\geq 1}\subset\mathcal{S}_N$.
Suppose that
$\{(x_0^n,y_0^n)\}_{n\geq 1}\subset E$ converges to	$(x_0,y_0)\in E$ as $n\to\infty$.
Then, there exist a subsequence
$\{(u^{n_k},v^{n_k})\}_{k\geq 1}$ and $(u,v)\in\mathcal{S}_N$ such that
$(u^{n_k},v^{n_k})$ converges  to $(u,v)$ with respect to the weak topology of
$\mathcal{H}$ as $k\to\infty$. Moreover,
\[
\mathcal{G}_{x_0^{n_k}}^0(u^{n_k},v^{n_k})
\longrightarrow
\mathcal{G}_{x_0}^0(u,v)
\]
in $\mathcal{C}^{\alpha^\prime}([0,T],\mathbb{R}^{d_1})$ as $k\to\infty$.
\end{corollary}
\begin{proof}
Since $\mathcal{S}_N$ is compact with respect to the
weak topology of $\mathcal{H}$, there exist a subsequence,
still denoted by $\{(u^n,v^n)\}_{n\geq1}$, and
$(u,v)\in\mathcal{S}_N$ such that
$(u^n,v^n)\to(u,v)$ with respect to the weak topology of $\mathcal{H}$ as $n\to\infty$.
Then, by the same argument as in the proof of Lemma \ref{lem-9}, the desired result follows.
\end{proof}


\begin{lemma}\label{lem-10}
Let $E\in\mathscr{E}$. Then, for every $(x_0,y_0)\in E$, the rate function $I_{x_0}$ defined in \eqref{eq-65} is a good rate function.
\end{lemma}

\begin{proof}
It suffices to show that, for any $0\leq M<\infty$, the level set $\Phi_{x_0}(M)$ is compact. By Lemma \ref{lem-9}, for each $n\geq1$, the set
$$
\Gamma_{M+\frac1n}
:=
\left\{
\mathcal{G}_{x_0}^0(u,v):
(u,v)\in\mathcal{S}_{M+\frac1n}
\right\}
$$
is compact. We will show that
$$
\Phi_{x_0}(M)
=
\bigcap_{n=1}^{\infty}\Gamma_{M+\frac1n}.
$$

First, let $\varphi\in\Phi_{x_0}(M)$, that is, $I_{x_0}(\varphi)\leq M$. By the definition of the rate function, for each $n\geq1$, there exists $(u^n,v^n)\in\mathcal{H}$ such that
$\varphi=\mathcal{G}_{x_0}^0(u^n,v^n)$
and
$\frac12\|(u^n,v^n)\|_{\mathcal{H}}^2
\leq M+\frac1n$.
Hence, $(u^n,v^n)\in\mathcal{S}_{M+\frac1n}$ and therefore
$\varphi\in\Gamma_{M+\frac1n}$
for every $n\geq1$. Thus,
$$
\varphi\in\bigcap_{n=1}^{\infty}\Gamma_{M+\frac1n}.
$$

Conversely, suppose that
$$
\psi\in\bigcap_{n=1}^{\infty}\Gamma_{M+\frac1n}.
$$
Then, for each $n\geq1$, there exists
$(u^n,v^n)\in\mathcal{S}_{M+\frac1n}$ such that
$
\psi=\mathcal{G}_{x_0}^0(u^n,v^n).
$
Consequently,
$
I_{x_0}(\psi)
\leq
\frac12\|(u^n,v^n)\|_{\mathcal{H}}^2
\leq M+\frac1n
$
for every $n\geq1$. Letting $n\to\infty$, we obtain
$
I_{x_0}(\psi)\leq M.
$
Hence,
$
\psi\in\Phi_{x_0}(M).
$

Therefore,
$$
\Phi_{x_0}(M)
=
\bigcap_{n=1}^{\infty}\Gamma_{M+\frac1n}.
$$

Since each $\Gamma_{M+\frac1n}$ is compact, and
$
\Phi_{x_0}(M)\subset\Gamma_{M+1},
$
while $\Phi_{x_0}(M)$ is closed as the intersection of the compact sets
$\Gamma_{M+\frac1n}$, it follows that $\Phi_{x_0}(M)$ is compact. Thus, $I_{x_0}$ is a good rate function.
\end{proof}

Next, we introduce a condition on $\mathcal{G}^{\varepsilon,\delta}_{x_0,y_0}$ guaranteeing that $X^{\varepsilon,\delta;x_0,y_0}=\mathcal{G}^{\varepsilon,\delta}_{x_0,y_0}(\sqrt\varepsilon B^H,\sqrt\varepsilon w)$ satisfies an equicontinuous uniform Laplace principle.
\begin{theorem}\label{th-4}
Assume that for every $E\in \mathscr{E}$, $N>0$ and $\zeta>0$,
\begin{align}\label{eq-67}
	\lim_{\substack{\varepsilon\to 0\\\delta/\varepsilon\to 0}} \sup_{(x_0,y_0)\in E}\sup_{(u,v)\in\mathcal{A}_b^N}
	\mathbb{P}\big(
	\big\|
	\mathcal{G}^{\varepsilon,\delta}_{x_0,y_0}(\sqrt\varepsilon B^H+u,\sqrt\varepsilon w+v)		-\mathcal{G}^0_{x_0}(u)
	\big\|_{\alpha^\prime}>\zeta
	\big)=0.
\end{align}
Then, the $\mathcal{C}^\alpha([0,T],\mathbb{R}^{d_1})$-valued random variables $X^{\varepsilon,\delta;x_0,y_0}$ satisfy an equicontinuous uniform Laplace principle with speed $ 1/\varepsilon$ and rate function $I_{x_0}$ uniformly over $\mathscr{E}$.
\end{theorem}
\begin{proof}
By Lemma \ref{lem2.3} and an argument analogous to that used in the proof of \cite[Theorem 2.13]{salins2019equivalences}, the result follows.
\end{proof}

\begin{remark}
	Since $\mathcal{C}^{\alpha^\prime}([0,T],\mathbb{R}^{d_1})$ is not separable, we apply the weak convergence approach in the separable little H\"older space $\tilde{\mathcal{C}}^{\alpha^\prime}$. Indeed, for any $\frac12<\alpha^\prime<\alpha<H$, we can choose $\iota>0$ such that $\alpha^\prime+\iota<\alpha$. 
	Then our process takes values in
	$\mathcal{C}^{\alpha^\prime+\iota}([0,T],\mathbb{R}^{d_1})$
	and, in particular, also belongs to the space $\tilde{C}^{\alpha^\prime}$, since  $\mathcal{C}^{\alpha^\prime+\iota}\subset\tilde{\mathcal{C}}^{\alpha^\prime}$.
	The variational formula is applied on the space $\tilde{C}^{\alpha^\prime}$, and we only need to prove the weak convergence method under the $\alpha^\prime$-H\"older norm.
	The resulting LDP extends to $\mathcal{C}^{\alpha^\prime}$ by the contraction principle \cite[Theorem 4.2.1]{dembo2009large}.
\end{remark}

Therefore, by the equivalence between the equicontinuous uniform Laplace principle and the Freidlin--Wentzell uniform LDP in Theorem \ref{th-6}, it suffices to prove the following result to establish the Freidlin--Wentzell uniform LDP.
\begin{lemma}\label{lem-11}
For every $E\in \mathscr{E}$, $N>0$ and $\zeta>0$, condition \eqref{eq-67} holds.
\end{lemma}

To establish convergence in $\mathcal{C}^{\alpha^\prime}$, we first recall the following interpolation inequality, which will be used in the proof of Lemma \ref{lem-11}.
\begin{proposition}\label{pro-12}
	Let $\eta>0$ be sufficiently small such that
	$
	0<\frac12-\eta<1-2\eta<1-\eta<1.
	$
	Then, for any continuous function $M:[0,T]\to\mathbb{R}^{d_1}$,
	$$
	\|M\|_{1-2\eta}
	\leq
	\|M\|_{1-\eta}^{1-2\eta}
	\|M\|_{\frac12-\eta}^{2\eta}.
	$$
	In particular, the inequality holds almost surely when $M$ is a continuous stochastic process.
\end{proposition}

\begin{proof}
	For any $0\leq s<t\leq T$, we have
	\begin{align*}
		\frac{|M_{s,t}|}{|t-s|^{1-2\eta}}
		= &
		\Big(\frac{|M_{s,t}|}{|t-s|^{1-\eta}}\Big)^a\Big(\frac{|M_{s,t}|}{|t-s|^{\frac12-\eta}}\Big)^b,
	\end{align*}
if
	\begin{align*}
		\begin{cases}
			a+b=1,\\
			(1-\eta)a+(\frac12-\eta)b=1-2\eta.
		\end{cases}
	\end{align*}
	It is easy to see that $a=1-2\eta$, and $b=2\eta$.
	Then, we have 
	\begin{align*}
		\frac{|M_{s,t}|}{|t-s|^{1-2\eta}}
		= &
		\Big(\frac{|M_{s,t}|}{|t-s|^{1-\eta}}\Big)^{1-2\eta}\Big(\frac{|M_{s,t}|}{|t-s|^{\frac12-\eta}}\Big)^{2\eta}.
	\end{align*}
	Taking the supremum over $0\leq s<t\leq T$ yields
	\begin{align*}
		\sup_{0\leq s<t\leq T}	\frac{|M_{s,t}|}{|t-s|^{1-2\eta}}
		\leq  &
		\sup_{0\leq s<t\leq T}	\Big(\frac{|M_{s,t}|}{|t-s|^{1-\eta}}\Big)^{1-2\eta}\sup_{0\leq s<t\leq T}\Big(\frac{|M_{s,t}|}{|t-s|^{\frac12-\eta}}\Big)^{2\eta}.
	\end{align*}
	This proves the desired inequality.
\end{proof}

\begin{proof}[Proof of Lemma \ref{lem-11}]
Throughout this proof, the constant $C>0$ is independent of $\varepsilon$, $\delta$, $\Delta$, $(u,v)\in\mathcal{A}^N_b$ and the initial conditions $(x_0,y_0)\in E$, which my change from line to line.

By direct calculation, we have
\begin{align*}
\tilde{X}^{\varepsilon,\delta}_t-\hat{X}^{\varepsilon,\delta}_t
=&
\int_0^t \big\{
b_1(s,\tilde{X}^{\varepsilon,\delta}_s,\tilde{Y}^{\varepsilon,\delta}_s)-b_1(s,\tilde{X}^{\varepsilon,\delta}_s,\hat{Y}^{\varepsilon,\delta}_s)
\big\}ds\\
&+
\int_0^t \big\{
b_1(s,\tilde{X}^{\varepsilon,\delta}_s,\hat{Y}^{\varepsilon,\delta}_s)-b_1(s,\tilde{X}^{\varepsilon,\delta}_{s(\Delta)},\hat{Y}^{\varepsilon,\delta}_s)
\big\}ds\\
&+
\int_0^t \big\{
b_1(s,\tilde{X}^{\varepsilon,\delta}_{s(\Delta)},\hat{Y}^{\varepsilon,\delta}_s)-\bar{b}_1(s,\tilde{X}^{\varepsilon,\delta}_{s(\Delta)})
\big\}ds\\
&+
\int_0^t \big\{
\bar{b}_1(s,\tilde{X}^{\varepsilon,\delta}_{s(\Delta)})-\bar{b}_1(s,\tilde{X}^{\varepsilon,\delta}_s)
\big\}ds\\
&+
\int_0^t \big\{
\bar{b}_1(s,\tilde{X}^{\varepsilon,\delta}_s)-\bar{b}_1(s,\hat{X}^{\varepsilon,\delta}_s)
\big\}ds\\
&+
\int_0^t \big\{
\sigma_1(s,\tilde{X}^{\varepsilon,\delta}_s)-\sigma_1(s,\hat{X}^{\varepsilon,\delta}_s)
\big\}d(\sqrt\varepsilon B^H_s+u_s)\\
=:&K^1_t+K^2_t+K^3_t+K^4_t+K^5_t+K^6_t.
\end{align*}

For the first term, since $b_1$ is Lipschitz continuous with respect to its spatial variables, for any $\theta\leq \frac12$, we have
\begin{align}\label{eq-64}
\mathbb{E}\big[\|K^1\|^2_\theta\big]
\leq &
C\mathbb{E}\Big[
\Big(
\sup_{0\leq  s<t\leq T} \frac{\int_s^t |\tilde{Y}^{\varepsilon,\delta}_r-\hat{Y}^{\varepsilon,\delta}_r|ds}{|t-s|^\theta}
\Big)^2
\Big]\notag\\
\leq &
C\mathbb{E}\Big[
\sup_{0\leq  s<t\leq T} \frac{\int_s^t |\tilde{Y}^{\varepsilon,\delta}_r-\hat{Y}^{\varepsilon,\delta}_r|^2ds}{|t-s|^{2\theta-1}}
\Big]\notag\\
\leq &
C\int_0^T \mathbb{E}\big[|\tilde{Y}^{\varepsilon,\delta}_s-\hat{Y}^{\varepsilon,\delta}_s|^2\big]ds\notag\\
\leq &
C\Big(\Delta+\frac\delta\varepsilon\Big).
\end{align}
Here, the second inequality follows from H\"older's inequality and the fact that $2\theta-1\leq0$, and the last inequality follows from Lemma \ref{lem-13}. 
Let $\eta>0$ be sufficiently small such that
$
0<\frac12-\eta<1-2\eta<1-\eta<1.
$
By Proposition \ref{pro-12}, we have
\begin{align*}
\mathbb{E}\big[\|K^1\|^2_{1-2\eta}\big]
\leq& \mathbb{E}\big[\|K^1\|^{2-4\eta}_{1-\eta}\|K^1\|^{4\eta}_{\frac12-\eta}\big]
\leq 
\mathbb{E}\big[\|K^1\|^2_{1-\eta}]^{1-2\eta}\mathbb{E}[\|K^1\|^2_{\frac12-\eta}\big]^{2\eta}.
\end{align*}
Here, the last inequality follows from H\"older's inequality. 
Since $b_1$ is bounded, we have 
$$\sup_{0<\delta<\varepsilon\leq 1}\sup_{(x_0,y_0)\in E}\sup_{(u,v)\in\mathcal{A}^N_b}\mathbb{E}\big[\|K^1\|^2_{1-\zeta}\big]<\infty.$$
Moreover, by \eqref{eq-64}, 
$$\sup_{(x_0,y_0)\in E}\sup_{(u,v)\in\mathcal{A}^N_b}\mathbb{E}\big[\|K^1\|^2_{\frac12-\eta}\big]\leq C\Big(\Delta+\frac\delta\varepsilon\Big).$$
Consequently,
\begin{align*}
\mathbb{E}\big[\|K^1\|^2_{1-2\eta}\big]\leq  C\Big(\Delta^{2\eta}+\Big(\frac\delta\varepsilon\Big)^{2\eta}\Big).
\end{align*}

Since $b_1$ is Lipschitz continuous, we have 
\begin{align*}
\mathbb{E}\big[\|K^2\|^2_1\big]
\leq &
C\mathbb{E}\Big[
\Big(
\sup_{0\leq s<t\leq T}\frac{\int_s^t |\tilde{X}^{\varepsilon,\delta}_r-\tilde{X}^{\varepsilon,\delta}_{r(\Delta)}|dr}{|t-s|}
\Big)^2
\Big]
\leq 
C\mathbb{E}\big[\|\tilde{X}^{\varepsilon,\delta}\|^2_{\alpha^\prime}\big]\Delta^{2\alpha^\prime}
\leq 
C\Delta^{2\alpha^\prime}.
\end{align*}
Here, the last inequality follows from Proposition \ref{pro-7}. 

By Propositions \ref{pro-5} and \ref{pro-7}, we have
\begin{align*}
\mathbb{E}\big[\|K^4\|^2_1\big]
\leq &
C\mathbb{E}\Big[
\Big(
\sup_{0\leq s<t\leq T}\frac{\int_s^t |\tilde{X}^{\varepsilon,\delta}_r-\tilde{X}^{\varepsilon,\delta}_{r(\Delta)}|(1+|\tilde{X}^{\varepsilon,\delta}_r|^\xi+|\tilde{X}^{\varepsilon,\delta}_{r(\Delta)}|^\xi)dr}{|t-s|}
\Big)^2
\Big]\\
\leq &
C\mathbb{E}\big[
\|\tilde{X}^{\varepsilon,\delta}\|^2_{\alpha^\prime}+\|\tilde{X}^{\varepsilon,\delta}\|^{2(1+\xi)}_{\alpha^\prime}
\big] \Delta^{2\alpha^\prime}\\
\leq &
C\Delta^{2\alpha^\prime}.
\end{align*}

For the term $K^3$, following the argument in the proof of Lemma \ref{lem-7}, we obtain, for any $\gamma\in(1/2,1)$,
\begin{align*}
\mathbb{E}\big[\|K^3\|^2_\gamma\big]
\leq &
C\Delta^{2(1-\gamma)}
+C \Delta^{-(1+2\gamma)}\max_{0\leq k \leq \lfloor\frac T\Delta\rfloor-1}\mathbb{E}\Big[
\Big|
\int_{k\Delta}^{(k+1)\Delta}\big\{
b_1(s,\tilde{X}^{\varepsilon,\delta}_{k\Delta},\hat{Y}^{\varepsilon,\delta}_s)-\bar{b}_1(s,\tilde{X}^{\varepsilon,\delta}_{k\Delta})
\big\}ds
\Big|^2
\Big].
\end{align*}
Let $\bar{w}=(\bar{w_t})_{t\geq 0}$ be another standard $r$-dimensional BM which is independent of $(B^H,w)$ and let $(\bar{Y}^{x,y}_t)_{t\geq 0}$ denote the unique solution of the frozen SDE driven by $(\bar{w_t})_{t\geq 0}$.
By Lemma \ref{lem-6}, for any $(x,y)\in\mathbb{R}^{d_1+d_2}$ and $t\geq 0$, 
\begin{align*}
\mathbb{E}\Big[
\Big|
\int_0^t \big\{
b_1(s,x,\bar{Y}^{x,y}_s)-\bar{b}_1(s,x)
\big\}ds
\Big|^2
\Big]
\leq C(1+|x|^{\frac{\eta_1}2}+|y|)t.
\end{align*}
Define 
$
\bar{Y}_t:=\bar{Y}^{x,y}_t|_{(x,y)=(\tilde{X}^{\varepsilon,\delta}_{k\Delta},\hat{Y}^{\varepsilon,\delta}_{k\Delta})}.
$
Since $\bar{w}$ is independent of $(B^H,w)$, it is also independent of $(X^\delta_{k\Delta},\hat{Y}^\delta_{k\Delta})$.
Then, $(\bar{Y}_{t/\delta})_{t\geq 0}$ satisfies  
\begin{align*}
\bar{Y}_{t/\delta}
=& \hat{Y}^{\varepsilon,\delta}_{k\Delta}+\int_0^{\frac t\delta}b_2(\tilde{X}^{\varepsilon,\delta}_{k\Delta},\bar{Y}_s)ds+\int_0^{\frac t\delta}\sigma_2(\tilde{X}^{\varepsilon,\delta}_{k\Delta},\bar{Y}_s)d^I\bar{w}_s\\
=&\hat{Y}^{\varepsilon,\delta}_{k\Delta}+\frac1\delta\int_0^tb_2(\tilde{X}^{\varepsilon,\delta}_{k\Delta},\bar{Y}_{s/\delta})ds+\frac1{\sqrt\delta}\int_0^t\sigma_2(\tilde{X}^{\varepsilon,\delta}_{k\Delta},\bar{Y}_{s/\delta})d^I(\sqrt\delta \bar{w}_{s/\delta}),
\end{align*}
where $(\sqrt\delta\bar{w}_{t/\delta})_{t\geq 0}$ is again a standard $r$-dimensional BM independent of $(\tilde{X}^{\varepsilon,\delta}_{k\Delta},\hat{Y}^{\varepsilon,\delta}_{k\Delta})$.
Comparing the equation above with \eqref{eq-51}, we see that, on $[0,\Delta)$,  $(\hat{Y}^{\varepsilon,\delta}_{k\Delta+t})$ and $(\bar{Y}_{t/\delta})$ satisfy the same frozen SDE with the same initial condition $\hat{Y}^{\varepsilon,\delta}_{k\Delta}$ and the same frozen vatiable $\tilde{X}^{\varepsilon,\delta}_{k\Delta}$.
Hence, by the uniqueness in law of the frozen SDE, 
\begin{align*}
\big(\tilde{X}^{\varepsilon,\delta}_{k\Delta},(\hat{Y}^{\varepsilon,\delta}_{k\Delta+t})_{t\in[0,\Delta)}\big)\stackrel{\mathrm{law}}{=}\big(\tilde{X}^{\varepsilon,\delta}_{k\Delta},(\bar{Y}_{t/\delta})_{t\in[0,\Delta)}\big).
\end{align*}
By the same argument as in the proof of Lemma \ref{lem-6}, we obain, for any
$(x,y)\in\mathbb{R}^{d_1+d_2}$ and $t\geq0$,
\begin{align*}
\mathbb{E}\Big[\Big|\int_0^t\big\{b_1(\delta s, x , Y^{x,y}_s)-\bar{b}_1(\delta s, x)\big\}ds\Big|^2\Big]\leq C(1+|x|^{\frac{\eta_1}2}+|y|)t.
\end{align*}
Here, $\eta_1\geq 0$ is the constant in \textbf{(H4)$_q$}.
Consequently, 
\begin{align*}
&\mathbb{E}\Big[
\Big|\int_{k\Delta}^{(k+1)\Delta}\big\{
b_1(s,\tilde{X}^{\varepsilon,\delta}_{k\Delta},\hat{Y}^{\varepsilon,\delta}_s)-\bar{b}_1(s,\tilde{X}^{\varepsilon,\delta}_{k\Delta})
\big\}ds\Big|^2
\Big]\\
&= 
\mathbb{E}\Big[
\Big|\int_0^\Delta\big\{
b_1(s,\tilde{X}^{\varepsilon,\delta}_{k\Delta},\bar{Y}_{s/\delta})-\bar{b}_1(s,\tilde{X}^{\varepsilon,\delta}_{k\Delta})
\big\}ds\Big|^2
\Big]\\
&= 
\delta^2\mathbb{E}\Big[
\Big|\int_0^{\frac\Delta\delta}\big\{
b_1(\delta s,\tilde{X}^{\varepsilon,\delta}_{k\Delta},\bar{Y}_s)-\bar{b}_1(\delta s,\tilde{X}^{\varepsilon,\delta}_{k\Delta})
\big\}ds\Big|^2
\Big]\\
&\leq 
C\delta^2\mathbb{E}\Big[\big(1+|\tilde{X}^{\varepsilon,\delta}_{k\Delta}|^{\frac{\eta_1}2}+|\hat{Y}^{\varepsilon,\delta}_{k\Delta}|\big)\frac\Delta\delta\Big]\\
&\leq C\Delta\delta,
\end{align*}
Where the last inequality follows from Proposition \ref{pro-7} and Lemma \ref{lem-8}.
Combining the above estimates, we obtain
\begin{align*}
\mathbb{E}\big[\|K^3\|^2_\gamma\big]\leq C\big(\Delta^{2(1-\gamma)}+\Delta^{-2\gamma}\delta\big).
\end{align*}

Set
\begin{align*}
V^{\varepsilon,\delta}_t:=&(\tilde{X}^{\varepsilon,\delta}_t-\hat{X}^{\varepsilon,\delta}_t)-\int_0^t \big\{\bar{b}_1(s,\tilde{X}^{\varepsilon,\delta}_s)-\bar{b}_1(s,\hat{X}^{\varepsilon,\delta}_s)\big\}ds\\
&-\int_0^t \big\{\sigma_1(s,\tilde{X}^{\varepsilon,\delta}_s)-\sigma_1(s,\hat{X}^{\varepsilon,\delta}_s)\big\}d(\sqrt\varepsilon B^H_s+u_s).
\end{align*}
Then, $V^{\varepsilon,\delta}_t=K^1_t+K^2_t+K^3_t+K^4_t$.
For the  term $K^1$, we take $1-2\eta=\alpha^\prime$, that is $\eta=\frac{1-\alpha^\prime}2$. For the term $K^3$, we take $\gamma=\alpha^\prime$.
Therefore,
\begin{align*}
\sup_{(x_0,y_0)\in E}
\sup_{(u,v)\in\mathcal{A}^N_b}
\mathbb{E}\big[\|V^{\varepsilon,\delta}\|^2_{\alpha^\prime}\big]
\leq &
C\Big(\Delta^{1-\alpha^\prime}+\Big(\frac\delta\varepsilon\Big)^{1-\alpha^\prime}+\Delta^{2\alpha^\prime}+\Delta^{2(1-\alpha^\prime)}+\Delta^{-2\alpha^\prime}\delta\Big)\\
\leq &
C\Big(\Delta^{1-\alpha^\prime}+\Big(\frac\delta\varepsilon\Big)^{1-\alpha^\prime}+\Delta^{-2\alpha^\prime}\delta\Big)
\end{align*}
Choosing a suitable \(\Delta\), for instance, $\Delta=\delta^\frac12$, then $\|V^{\varepsilon,\delta}\|_{\alpha^\prime}\to 0$ in $L^2$-sense uniformly with respect to $(x_0,y_0)\in E$ as $\varepsilon\to 0$ and $\delta/\varepsilon\to 0$.

From Proposition \ref{pro-10}, there exist constants $C,l>0$ such that 
\begin{align*}
\|\tilde{X}^{\varepsilon,\delta}-\hat{X}^{\varepsilon,\delta}\|_{\alpha^\prime}\leq C\exp(C\|\sqrt\varepsilon B^H+u\|^l_\alpha)\|V^{\varepsilon,\delta}\|_{\alpha^\prime}
\leq C_{l,N} \exp{\big(C\|B^H\|^l_\alpha\big)}\|V^{\varepsilon,\delta}\|_{\alpha^\prime}
\quad\text{a.s.}
\end{align*}

Since
$\mathbb{E}\big[\|B^H\|^l_\alpha\big]<\infty$,
we have 
$\exp(C\|B^H\|^l_\alpha)<\infty$ almost surely.
On the other hand, since $\|V^{\varepsilon,\delta}\|_{\alpha^\prime}\to 0$ in $L^2$-sense uniformly, by Markov's inequality, for any $R>0$, we have
\begin{align*}
\lim_{\substack{\varepsilon\to 0\\\delta/\varepsilon\to 0}}\sup_{(x_0,y_0)\in E}
\sup_{(u,v)\in\mathcal{A}^N_b}\mathbb{P}\big(\|V^{\varepsilon,\delta}\|_{\alpha^\prime}>R\big)=0.
\end{align*}


Then, for any $\zeta>0$ and $K>0$, we have
\begin{align*}
&\limsup_{\substack{\varepsilon\to 0\\\delta/\varepsilon\to 0}}
\sup_{(x_0,y_0)\in E}
\sup_{(u,v)\in\mathcal{A}^N_b}\mathbb{P}\big(\exp{\big(C\|B^H\|^l_\alpha\big)}\|V^{\varepsilon,\delta}\|_{\alpha^\prime}>\zeta\big)\\
&\leq 
\mathbb{P}\big(\exp{\big(C\|B^H\|^l_\alpha\big)}\geq K\big)+
\limsup_{\substack{\varepsilon\to 0\\\delta/\varepsilon\to 0}}
\sup_{(x_0,y_0)\in E}
\sup_{(u,v)\in\mathcal{A}^N_b}\mathbb{P}\Big(\|V^{\varepsilon,\delta}\|_{\alpha^\prime}>\frac\zeta K\Big)\\
&=\mathbb{P}\big(\exp{\big(C\|B^H\|^l_\alpha\big)}\geq K\big),
\end{align*}
where the last indentity follows from the uniform convergence in probability of $\|V^{\varepsilon,\delta}\|_{\alpha^\prime}$.
Since $\exp(C\|B^H\|^l_\alpha)<\infty$ almost surely, letting $K\to \infty$, we obtain
\begin{align}\label{eq-68}
	\lim_{\substack{\varepsilon\to 0\\\delta/\varepsilon\to 0}}
	\sup_{(x_0,y_0)\in E}
	\sup_{(u,v)\in\mathcal{A}^N_b}
	\mathbb{P}\big(\|\tilde{X}^{\varepsilon,\delta}-\hat{X}^{\varepsilon,\delta}\|_{\alpha^\prime}>\zeta\big)=0.
\end{align}


Define $\tilde{X}:=\mathcal{G}_{x_0}^0(u,0)\in\mathcal{C}^\alpha([0,T],\mathbb{R}^{d_1})$, which is the unique global solution to the following YDE:
\begin{align*}
d\tilde{X}_t=\bar{b}_1(t,\tilde{X}_t)dt+\sigma_1(t,\tilde{X}_t)du_t,\quad \tilde{X}_0=x_0,\quad t\in[0,T].
\end{align*}
Moreover, by Proposition \ref{pro-4}, $\sup_{(x_0,y_0)\in E}\sup_{(u,v)\in\mathcal{A}^N_b}\mathbb{E}[\|\tilde{X}\|^p_{\alpha^\prime}]<\infty$ for any $1\leq p<\infty$.
By applying Proposition \ref{pro-11} to the first term and H\"older's inequality to the second term on the right-hand side of the first inequality, we obtain, 
\begin{align*}
&\sup_{(x_0,y_0)\in E}
\sup_{(u,v)\in \mathcal{A}^N_b}
\mathbb{E}\big[\|\hat{X}^{\varepsilon,\delta}-\tilde{X}\|_{\alpha^\prime}\big]\\
&\leq\sup_{(x_0,y_0)\in E}
\sup_{(u,v)\in\mathcal{A}^N_b}
\mathbb{E}\Big[\|\hat{X}^{\varepsilon,\delta}-\tilde{X}\|_{\alpha^\prime}\textbf{1}_{\{\sqrt\varepsilon\|B^H\|_\alpha\leq 1\}} \Big] \\
&\quad+
\sup_{(x_0,y_0)\in E}
\sup_{(u,v)\in\mathcal{A}^N_b}
\mathbb{E}\Big[\|\hat{X}^{\varepsilon,\delta}-\tilde{X}\|_{\alpha^\prime}\textbf{1}_{\{\sqrt\varepsilon\|B^H\|_\alpha> 1\}} \Big] \\
&\leq 
C\sqrt\varepsilon\mathbb{E}\big[\|B^H\|_\alpha\big]
+\sup_{(x_0,y_0)\in E}
\sup_{(u,v)\in\mathcal{A}^N_b}
\mathbb{E}\big[\|\hat{X}^{\varepsilon,\delta}-\tilde{X}\|^2_{\alpha^\prime}\big]^{\frac12}\times\mathbb{P}\Big(\|B^H\|_\alpha>\frac 1{\sqrt\varepsilon}\Big)
^{\frac12}\\
&\leq C\sqrt\varepsilon+C\mathbb{P}\Big(\|B^H\|_\alpha>\frac 1{\sqrt\varepsilon}\Big)
^{\frac12}\to 0\quad\text{as }\varepsilon\to 0.
\end{align*}
Here, the last inequality follows from the uniform boundedness of the second moments of $\hat{X}^{\varepsilon,\delta}$ and $\tilde{X}$.
Therefore, by Markov's inequality, we have, for any $\zeta>0$, 
\begin{align}\label{eq-69}
\lim_{\varepsilon\to 0}
\sup_{(x_0,y_0)\in E}
\sup_{(u,v)\in\mathcal{A}^N_b}
\mathbb{P}\big(\|\hat{X}^{\varepsilon,\delta}-\tilde{X}\|_{\alpha^\prime}>\zeta\big)=0
\end{align}

Combining \eqref{eq-68} and \eqref{eq-69}, together with the triangle inequality, we obtain, for any $\zeta>0$,
\begin{align*}
	\lim_{\substack{\varepsilon\to 0\\\delta/\varepsilon\to 0}}
	\sup_{(x_0,y_0)\in E}
	\sup_{(u,v)\in\mathcal{A}^N_b}
	\mathbb{P}\big(\|\tilde{X}^{\varepsilon,\delta}-\tilde{X}\|_{\alpha^\prime}>\zeta\big)=0.
\end{align*}
This completes the proof.
\end{proof}

\begin{proof}[Proof of Theorem \ref{th-3}]
By Lemma \ref{lem-11} and Theorem \ref{th-4}, the family
\(\{X^{\varepsilon,\delta;x_0,y_0}\}\) satisfies the equicontinuous uniform Laplace principle.
Combining Lemma \ref{lem-10} with Theorem \ref{th-6}, we obtain
the Freidlin--Wentzell uniform LDP over \(\mathscr{E}\), with the good rate function
\(I_{x_0}\) defined in \eqref{eq-65}, for
\((x_0,y_0)\in E\), where \(E\in\mathscr{E}\).
\end{proof}

\appendix
\section{Definition of Uniform Large Deviation Principle}\label{sec-A}

We recall the definition of uniform LDP. Let $(\mathcal{V},\rho)$ be a Polish space and let $\mathcal{V}_0$ be an index set. 
Let $(\Omega,\mathcal{F}, \mathbb{P})$ be a probability space, and let $\{X^\varepsilon_x:\varepsilon>0,x\in\mathcal{V}_0\}$ be a family of $\mathcal{V}$-valued random variables indexed by $x\in\mathcal{V}_0$. 
Denote by $\mathbb{E}$ the expectation with respect to $\mathbb{P}$. Let $\{I_x:x\in\mathcal{V}_0\}$ be a family of lower semicontinuous rate functions $I_x:\mathcal{V}\to[0,+\infty]$, with level sets
$
\Phi_x(s):=\{\varphi\in\mathcal{V}:I_x(\varphi)\leq s\},~s\geq0.
$
The rate function $I_x$ is called good if $\Phi_x(s)$ is compact in $\mathcal{V}$ for every $s\geq0$.
In fact, in our setting, $\mathcal{V}=\tilde{\mathcal{C}}^{\alpha^\prime}([0,T],\mathbb{R}^{d_1})$ and $\mathcal{V}_0=\mathbb{R}^{d_1+d_2}$.
\begin{definition}[Freidlin-Wentzell uniform large deviation principle]\label{FWULDP}
Let $\mathscr{V}$ be a collection of subsets of $\mathcal{V}_0$ and $a(\varepsilon)$ be a function converging to zero as $\varepsilon$ converges to zero. The random variables $\{X^\varepsilon_x\}$ are said to satisfy a Freidlin-Wentzell uniform large deviations principle with respect to the rate functions $I_x$ with speed $1/a(\varepsilon)$ uniformly over $\mathscr{V}$ , if
	\begin{itemize}
	\item For any $A\in\mathscr{V}$, $s_0>0$, and $\zeta>0$,
	$$\liminf_{\varepsilon\to0}\inf_{x\in A}
	\inf_{\varphi\in\Phi_x(s_0)}
	\big(a(\varepsilon)\log\mathbb{P}\big(\rho(X^\varepsilon_x,\varphi)<\zeta\big)+I_x(\varphi)\big)\geq 0.$$
	\item For any $A\in\mathscr{V}$, $s_0>0$, and $\zeta>0$,
	$$\limsup_{\varepsilon\to0}\sup_{x\in A}
	\sup_{s\in[0,s_0]}
	\big(a(\varepsilon)\log\mathbb{P}\big(
	\inf_{\psi\in \Phi_x(s)} \rho(X^\varepsilon_x,\psi)\geq \zeta\big)
	+s\big)\leq 0.$$
\end{itemize}
\end{definition}

\begin{definition}[Equicontinuous uniform Laplace principle]\label{EULP}
 Let $\mathscr{V}$ be a collection of subsets of $\mathcal{V}_0$ and $a(\varepsilon)$ be a function converging to zero as $\varepsilon$ converges to zero. The random variables $\{X^\varepsilon_x\}$ are said to satisfy an equicontinuous uniform Laplace principle with respect to the rate functions $I_x$ with speed $1/a(\varepsilon)$ uniformly over $\mathscr{V}$, if for any $A \in \mathscr{V}$ and any collection $L \subset C_b(\mathcal{V})$ of equibounded and equicontinuous functions from $\mathcal{V}$ to $\mathbb{R}$,
 \begin{align*}
\lim_{\varepsilon\to 0}\sup_{x\in A}\sup_{h\in L}
\Big|
a(\varepsilon)\log\mathbb{E}\exp\Big(-\frac{h(X^\varepsilon_x)}{a(\varepsilon)}\Big)
+\inf_{\varphi\in\mathcal{V}}\big\{h(\varphi)+I_x(\varphi)\big\}
\Big|=0.
 \end{align*}
\end{definition}


We recall the following results from Salins \cite{salins2019equivalences}.
\begin{theorem}\label{th-6}
The equicontinuous uniform Laplace principle and the Freidlin--Wentzell uniform LDP are equivalent with no extra assumptions.
\end{theorem}
	
\bigskip
\noindent
{\large {\bf Acknowlegement:}}
The author gratefully acknowledges the financial support provided by the China Scholarship Council.

\bigskip
\noindent
{\large {\bf Declaration of Generative AI and AI-Assisted Technologies}}
\\
During the preparation of this work, the author used OpenAI ChatGPT solely for language editing. All mathematical formulas, proofs, derivations, and numerical results presented in this manuscript were independently produced and verified by the author, who takes full responsibility for their accuracy and for the content of the publication.

\bibliographystyle{plain}
\bibliography{tempbib}
\end{document}